\documentclass[leqno]{amsart}

\usepackage[leqno]{amsmath}
\usepackage{amssymb,amsthm}
\usepackage{enumerate, latexsym, mathrsfs, mdwlist}
\usepackage{booktabs}

\swapnumbers
\theoremstyle{definition}
\newtheorem{nul}{}[section]
\newtheorem{dfn}[nul]{Definition}

\newtheorem{cnstr}[nul]{Construction}
\newtheorem{ntn}[nul]{Notation}
\newtheorem{exm}[nul]{Example}

\newtheorem{wrn}[nul]{Warning}
\newtheorem*{dfn*}{Definition}
\newtheorem*{axm*}{Axiom}
\newtheorem*{ntn*}{Notation}
\newtheorem*{exm*}{Example}
\newtheorem*{exr*}{Exercise}
\newtheorem*{int*}{Intuition}
\newtheorem*{qst*}{Question}

\theoremstyle{plain}

\newtheorem{thm}[nul]{Theorem}
\newtheorem{prp}[nul]{Proposition}
\newtheorem{lem}[nul]{Lemma}

\newtheorem{cnj}[nul]{Conjecture}
\newtheorem{cor}{Corollary}[nul]
\newtheorem*{thm*}{Theorem}
\newtheorem*{prp*}{Proposition}
\newtheorem*{cor*}{Corollary}
\newtheorem*{lem*}{Lemma}
\newtheorem*{cnj*}{Conjecture}

\numberwithin{equation}{nul}

\DeclareMathOperator{\colim}{colim}

\DeclareMathOperator{\fib}{fib}
\DeclareMathOperator{\Fun}{Fun}

\DeclareMathOperator{\id}{id}

\DeclareMathOperator{\Map}{Map}
\DeclareMathOperator{\Mor}{Mor}
\DeclareMathOperator{\Obj}{Obj}
\DeclareMathOperator{\Pair}{Pair}

\newcommand{\CC}{\mathbf{C}}

\newcommand{\FF}{\mathbf{F}}
\newcommand{\GG}{\mathbf{G}}

\newcommand{\II}{\mathbf{I}}

\newcommand{\NN}{\mathbf{N}}
\newcommand{\OO}{\mathbf{O}}
\newcommand{\PP}{\mathbf{P}}

\newcommand{\RR}{\mathbf{R}}

\newcommand{\XX}{\mathbf{X}}

\newcommand{\Adm}{\mathbf{Adm}}
\newcommand{\Alg}{\mathbf{Alg}}
\newcommand{\Cat}{\mathbf{Cat}}
\newcommand{\Cor}{\mathbf{Cor}}

\newcommand{\Kan}{\mathbf{Kan}}
\newcommand{\Mnd}{\mathbf{Mnd}}

\newcommand{\Mon}{\mathbf{Mon}}
\newcommand{\Op}{\mathbf{Op}}
\newcommand{\Operad}{\mathbf{Operad}}

\newcommand{\Set}{\mathbf{Set}}

\newcommand{\cart}{\mathrm{cart}}

\newcommand{\op}{\mathrm{op}}

\newcommand{\coloneq}{\mathrel{\mathop:}=}

\makeatletter 
\def\revddots{\mathinner{\mkern1mu\raise\p@ 
\vbox{\kern7\p@\hbox{.}}\mkern2mu 
\raise4\p@\hbox{.}\mkern2mu\raise7\p@\hbox{.}\mkern1mu}} 
\makeatother

\usepackage{tikz}
\usetikzlibrary{matrix,arrows,decorations}

\pgfdeclaredecoration{single line}{initial}{\state{initial}[width=\pgfdecoratedpathlength-1sp]{\pgfpathmoveto{\pgfpointorigin}}\state{final}{\pgfpathlineto{\pgfpointorigin}}} 

\newcommand{\fromto}[2]{{#1}\ \tikz[baseline]\draw[>=stealth,->](0,0.5ex)--(0.5,0.5ex);\ {#2}}
\newcommand{\into}[2]{{#1}\ \tikz[baseline]\draw[>=stealth,right hook->](0,0.5ex)--(0.5,0.5ex);\ {#2}}

\newcommand{\equivto}[2]{{#1}\ \tikz[baseline]\draw[>=stealth,->,font=\scriptsize,inner sep=0.5pt](0,0.5ex)--node[above]{$\sim$}(0.5,0.5ex);\ {#2}}
\newcommand{\goesto}[2]{{#1}\ \tikz[baseline]\draw[|->](0,0.5ex)--(0.5,0.5ex);\ {#2}}
\newcommand{\adjunct}[4]{{#1}\colon{#2}\ \begin{tikzpicture}[baseline] \draw[>=stealth,->] (0,1ex) -- (0.75,1ex); \draw[>=stealth,->] (0.75,0.25ex) -- (0,0.25ex); \end{tikzpicture}\ {#3}\colon{#4}}
\renewcommand{\to}{\ \tikz[baseline]\draw[>=stealth,->](0,0.5ex)--(0.5,0.5ex);\ }
\newcommand{\ot}{\ \tikz[baseline]\draw[>=stealth,<-](0,0.5ex)--(0.5,0.5ex);\ }

\title{From operator categories to higher operads}

\author{Clark Barwick}
\address{School of Mathematics, University of Edinburgh, James Clerk Maxwell Building, Peter Guthrie Tait Road, Edinburgh EH9 3FD, UNITED KINGDOM}
\email{clarkbar@gmail.com}

\begin{document}

\begin{abstract} In this paper we introduce the notion of an \emph{operator category} and two different models for homotopy theory of $\infty$-operads over an operator category -- one of which extends Lurie's theory of $\infty$-operads, the other of which is completely new, even in the commutative setting. We define perfect operator categories, and we describe a category $\Lambda(\Phi)$ attached to a perfect operator category $\Phi$ that provides Segal maps. We define a wreath product of operator categories and a form of the Boardman--Vogt tensor product that lies over it. We then give examples of operator categories that provide universal properties for the operads $A_n$ and $E_n$ ($1\leq n\leq +\infty$), as well as a collection of new examples.
\end{abstract}

\maketitle

\setcounter{tocdepth}{1}
\tableofcontents

\addtocounter{section}{-1}

\section{Introduction} A monoid structure on a set $X$ is the data of a product $\prod_{j\in J}x_j\in X$ of a collection of elements $\{x_j\}_{j\in J}$ indexed on a totally ordered finite set $J$. These multiplications are compatible with each other in the following sense: if $\phi\colon\fromto{J}{I}$ is an order-preserving map, then one has
\[
\prod_{j\in J}x_j=\prod_{i\in I}\prod_{j\in J_i}x_j.
\]
(Here $J_i$ is the fiber of $\phi$ over $i\in I$.) When $\phi$ is the inclusion $\into{\{1\}}{\{1,2\}}$ or the inclusion $\into{\{2\}}{\{1,2\}}$, this expresses the existence of the (right and left) unit in $X$. One extracts the associativity from the consideration of the two surjective, order-preserving maps $\fromto{\{1,2,3\}}{\{1,2\}}$. If one drops the ordering on these sets, then one has access to an involution on $\{1,2\}$ that expresses the commutativity of the monoid structure. The suggestion, therefore, is that the category $\OO$ of totally ordered finite sets and order-preserving maps ``controls'' the theory of monoids, and the category $\FF$ of finite sets ``controls'' the theory of commutative monoids.

If $X$ is instead a space (or a category, or a higher category), one may turn the identities in the description above into specified homotopies to obtain the structures observed on loopspaces and infinite loop spaces, respectively. There are (at least) two ways to make this idea precise:
\begin{itemize}
\item One selects a suitably cofibrant model $E_1$ for the unit non-symmetric operad and a suitably cofibrant model $E_{\infty}$ of the unit symmetric operad. Peter May showed that a grouplike $E_1$ (respectively, $E_{\infty}$) algebra structure on $X$ is essentially equivalent to a single delooping (resp., an infinite delooping) structure thereupon.
\item Alternately, one considers the category $\Delta$ of nonempty totally ordered finite sets and the category $\Gamma$ opposite to that of pointed finite sets. Graeme Segal observed that an $E_1$ algebra structure on a space $X$ is equivalent to a functor $\fromto{\Delta^{\op}}{\mathbf{Top}}$ that carries $1$ to $X$ and satisfies the so-called Segal condition, and an $E_{\infty}$ algebra structure on a space $X$ is equivalent to a functor $\fromto{\Gamma^{\op}}{\mathbf{Top}}$ that carries $1$ to $X$ and satisfies the Segal condition.
\end{itemize}
Our motivating -- but inchoate -- observation is that the category $\OO$ completely determines both the operad $E_1$ and Segal's theory of special $\Delta$-spaces, while the category $\FF$ completely determines both the operad $E_{\infty}$ and Segal's theory of special $\Gamma$-spaces.

To make precise what one means by ``determines,'' we intoduce here the notion of an \emph{operator category}. An operator category is a locally finite category that admits a terminal object $1$ and, for any object $I$, any morphism $i\colon\fromto{1}{I}$, and any map $\phi\colon\fromto{J}{I}$, a \emph{fiber} $J_i$ (i.e., a pullback of $\phi$ along $i$). Roughly speaking, an object $I$ of an operator category can be regarded as a finite set $|I|=\Mor(1,I)$ equipped with a suitable additional structure, and morphisms are structure-preserving maps. Of course $\OO$ and $\FF$ are operator categories. Attached to any operator category $\Phi$ is a theory of a \emph{$\Phi$-operads}. In effect, a $\Phi$-operad $P$ consists of the following: for each object $I\in\Phi$, a space $P(I)$, a point $e\in P(1)$ and for each morphism $\phi\colon\fromto{J}{I}$, a composition map
\[
\fromto{P(I)\times\prod_{i\in I}P(J_i)}{P(J)}.
\]
These data are required to satisfy an associativity condition for any composable pair of maps $K\to J\to I$ and a condition exhibiting $e$ as the unit. (In praxis, it is far simpler for us to work directly with a suitable theory of $\infty$-operads.) When $\Phi=\OO$, one obtains the theory of nonsymmetric operads, and when $\Phi=\FF$, one obtains the theory of symmetric operads. Of course $E_1$ is a suitably cofibrant replacement of the unit nonsymmetric operad, and $E_{\infty}$ is a suitably cofibrant replacement of the unit symmetric operad.

At the same time, the categories $\OO$ and $\FF$ themselves support canonically defined monads, which we shall denote $T_{\OO}$ and $T_{\FF}$. (For this, $\OO$ and $\FF$ are required to enjoy an additional technical property, which we call \emph{perfection}, but we emphasize that under this hypothesis the monad is \emph{not} additional structure.) In effect, these monads add points to any object in all the ways one can do so functorially. So $T_{\OO}$ adds both a new minimal point and a new maximal point to any totally ordered finite set, whereas $T_{\FF}$ simply adds a point to a finite set. An old observation of Joyal then shows that the \emph{Kleisli category} of the monad (i.e., the full subcategory of algebras spanned by the free algebras) on $\OO$ is naturally equivalent to $\Delta^{\op}$, and it is obvious that the Kleisli category of the monad on $\FF$ is $\Gamma^{\op}$. The \emph{Segal condition} for a functor from one of these categories to spaces ensures that $X(I)$ is equivalent to $\prod_{i\in|I|}X(\{i\})$.

For each operator category $\Phi$, one obtains an associated theory of $\Phi$-operads, and one can form a suitably cofibrant replacement of the unit $\Phi$-operad, whose algebras we will just call \emph{$\Phi$-algebras}. (Again, we will actually work with models of $\infty$-operads, which permit us to dodge the delicate cofibrancy issues for strict operads.) If $\Phi$ happens to be perfect, then it supports a monad $T$, and the Kleisli category $\Lambda(T)$ of $T$ contains a class of Segal maps $\chi_i$, and we obtain an equivalence of homotopy theories between $\Phi$-algebras in spaces and functors $\fromto{\Lambda(\Phi)}{\mathbf{Top}}$ that satisfy the Segal condition as above.

All of this is determined the instant one has a tiny quantity of combinatorial data -- the operator category $\Phi$. Moreover, the axioms for an operator category are so invitingly uncomplicated that one cannot help constructing other examples.  The surprise -- the biggest of this paper -- is that many operads arise in this manner.

For example, there is a wreath product of operator categories, which we describe in \S \ref{sect:wreath}: if $\Phi$ and $\Psi$ are operator categories, then $\Psi\wr\Phi$ is the category whose objects are pairs $(M_I,I)$ consisting of an object $I\in\Phi$ and an object $M_I\in\Psi^{\times|I|}$. (Morphisms are defined in the obvious way.) If $\Phi$ and $\Psi$ are perfect, then so is $\Psi\wr\Phi$. We also construct an \emph{external} Boardman--Vogt tensor product \cite{MR524181}, which takes a $\Psi$-operad $Q$ and a $\Phi$-operad $P$, and constructs a $(\Psi\wr\Phi)$-operad $Q\otimes P$. In effect, $(Q\otimes P)$-algebras are $Q$-algebras in $P$-algebras. At the same time, we show that if $P$ is the unit $\Phi$-operad and $Q$ is the unit $\Psi$-operad, then $Q\otimes P$ is the unit $(\Psi\wr\Phi)$-operad.

So if we form the iterated wreath product
\[
\OO^{(n)}\coloneq\OO\wr\OO\wr\cdots\wr\OO,
\]
then we get a perfect operator category. One sees readily that the theory of (strict) $\OO^{(n)}$-operads coincides with Michael Batanin's theory of $(n-1)$-terminal $n$-operads \cite{MR2365200}. Using our theory of $\infty$-operads, we prove the following.
\begin{thm*} The homotopy theory of $\OO^{(n)}$-algebras is equivalent to that of $E_n$-algebras, and the Kleisli category of the monad on $\OO^{(n)}$ is precisely Joyal's disk category $\Theta_n^{\op}$.
\end{thm*}
\noindent This can be regarded as an $\infty$-categorical analogue of Batanin's main result \cite{MR2365200}; the flexibility of the present context allows for a dramatically simpler and more conceptual proof. Moreover, since the tensor product of unit operads is again a unit operad, this yields a description of $E_{k_1+\cdots+k_m}$ as an iterated Boardman--Vogt tensor product $E_{k_1}\otimes E_{k_2}\otimes\cdots\otimes E_{k_m}$. This is not in itself new. In 1988, Gerald Dunn proved a very strict version of this fact using ordinary operads, and Michael Brinkmeier extended this result in 2000. Later, Zbigniew Fiedorowicz and Rainer Vogt proved a more general Additivity Theorem, stating that any cofibrant $E_k$ operad tensored with any cofibrant $E_{\ell}$ operad is $E_{k+\ell}$. However, because the Boardman--Vogt tensor product of strict operads does not preserve weak equivalences, these results are homotopically delicate (and some early attempts to prove a result of this kind appear to be incorrect). For this reason, Lurie's approach via $\infty$-operads is far better adapted to our work, and our expression of Additivity using $\OO^{(n)}$-operads is completely combinatorial.

Another, even simpler, construction of operator categories takes an operator category $\Phi$ and forms the full subcategory $\Phi_{\leq m}$ spanned by those objects $I$ such that $|I|$ is of cardinality $\leq m$. This is again an operator category, but it is not perfect, so only the operad story is accessible. A $\Phi_{\leq m}$-algebra is then much like a $\Phi$-algebra, except that one has only operations of arity $\leq m$. Consequently, one sees that $\OO_{\leq m}$-algebras are precisely $A_m$-algebras. We do not know a standard name for $\FF_{\leq m}$-algebras; we write $F_m$ for the operad freely generated by the operations of the $E_{\infty}$ operad of arity $\leq m$, so that $\FF_{\leq m}$-algebras are $F_m$-algebras. More generally still, one has an orbital category $\OO^{(n)}_{\leq m}$, providing a bifiltration of the $E_{\infty}$ operad (Ex. \ref{exm:bifiltrobinson}).

Here is a table summarizing the situation:

\begin{center}
\begin{tabular}{cccc}
\toprule
$\Phi$ & $\Phi$-algebras in $V$ &  perfect? & $\Lambda(\Phi)$ \\
\midrule
{$\{1\}$} & $\Alg_{E_0}(V)$ & yes & {$\{1\}$} \\
{$\OO$} & $\Alg_{E_1}(V)$ & yes & {$\Delta^{\op}$} \\
{$\OO^{(n)}$} & $\Alg_{E_n}(V)$ & yes & {$\Theta_n^{\op}$} \\
{$\FF$} & $\Alg_{E_{\infty}}(V)$ & yes & {$\Gamma^{\op}$} \\
$\OO_{\leq n}$ & $\Alg_{A_n}(V)$ & no & \\
$\FF_{\leq n}$ & $\Alg_{F_n}(V)$ & no & \\
\bottomrule
\end{tabular}
\end{center}

We introduce the $2$-category of operator categories in \S \ref{sect:opcats}. Given an operator category $\Phi$, we develop the homotopy theory of weak operads in two ways. The first and simplest of these models (\S \ref{sect:wkoperads}) is as suitable families of spaces over the nerve of a tree-like category $\Delta^{\op}_{\Phi}$ of finite sequences of objects of $\Phi$ that is related to a version of the dendroidal category with level structure when $\Phi=\FF$; this homotopy theory --- that of \emph{complete Segal $\Phi$-operads} --- can be described for any operator category and is easy to describe without many additional combinatorial complications. The second, a generalization of Lurie's theory of $\infty$-operads, requires more technology, and will come later (\S \ref{sect:Luriestyle}).

To prove this, we pass to a different way of describing weak operads and their algebras, which will only work for a special class of operator categories that we call \emph{perfect} (\S \ref{sect:perfect}). In effect, a perfect operator category $\Phi$ admits a universal way of adding elements to any object. This defines a monad on $\Phi$ (\S \ref{sect:canonicalmonad}), and the free algebras for this monad form a supplementary category $\Lambda(\Phi)$, called the \emph{Leinster category} (\S \ref{sect:Leinster}), which provides an alternative way of parametrizing the operations of an operad over $\Phi$.

Our second model of the homotopy theory of weak operads over a perfect operator category --- that of \emph{$\Phi$-quasioperads} --- is a natural generalization of Lurie's theory of $\infty$-operads (\S \ref{sect:Luriestyle}). In effect, a $\Phi$-quasioperad is an inner fibration $\fromto{X^{\otimes}}{N\Lambda(\Phi)}$ enjoying certain properties analogous to the ones developed by Lurie. We show that this and our original model of the homotopy theory of weak operads over $\Phi$ are equivalent in manner compatible with changes of operator category and the formation of $\infty$-categories of algebras (\S \ref{sect:equivalence}).

There is a monoidal structure on the $\infty$-category of all operator categories called the \emph{wreath product}, which we describe in \S \ref{sect:wreath}. An object of the wreath product $\Psi\wr\Phi$ is a pair consisting of an object $I\in\Phi$ and a collection of objects $\{J_i\ |\ i\in I\}$ indexed by the elements of $I$. The Boardman--Vogt tensor product of operads is externalized relative to this wreath product in \S \ref{sect:BV}, and we deduce our identification of $\OO^{(n)}$-algebras with $E_n$-algebras.

\subsection*{An apology} The results of this paper largely date from 2005--2006. Since that time, there has been a series of dramatic advances in the study of homotopy coherent algebraic structures, spearheaded by Moerdijk--Weiss and Lurie. As soon as I managed to record some of these results, the techniques I employed had, discouragingly, become outmoded. It was some time before it became clear to me how the work here interacts with some of these new advances --- particularly with Lurie's framework. Ultimately, these advances have simplified the work here greatly; the main result of this paper is now an immediate consequence of the comparison between our theory of weak operads and Lurie's. Nevertheless, my efforts to bring this work in line with current technology have led to an embarrassing delay in the publication of this work. I apologize for this, especially to the students who have sought to employ aspects of the theory introduced here.

\subsection*{Acknowledgments} All of this was inspired by Bertrand To\"en's visionary preprint \cite{dualite}. Early conversations with Markus Spitzweck were instrumental to my understanding. It was an offhand remark in a paper of Tom Leinster \cite[pp. 40--43]{leinsterhaops} that led me to formulate the notion of perfection for operator categories; the \emph{Leinster category} of \S \ref{sect:Leinster} is thus named after him. Haynes Miller has always been very kind in his support and encouragement, and he has asked a number of questions that helped refine my understanding of these objects. More recently, conversations with Chris Schommer-Pries revitalized my interest in this material; note that Lm. \ref{lem:actinert} in this paper is due to him. A visit from Clemens Berger helped me understand better how the work here interacts with concepts he's been developing since the dawn of the new millennium. In addition, I've benefitted from conversations with David Ayala, Ezra Getzler, John Rognes, and Sarah Whitehouse.


\section{Operator categories}\label{sect:opcats} The objects of an operator category are finite sets equipped with some additional structure. Such an object will be regarded as an indexing set for some multiplication law. The structure of the operator category can thus be thought of determining the associativity and commutativity constraints on that law.

\begin{ntn} For any ordinary category $\Phi$ with a terminal object $1$ and for any object $K\in\Phi$, we write $|K|\coloneq\Mor_{\Phi}(1,K)$. For any $i\in|K|$, it will be convenient to denote the morphism $i\colon\fromto{1}{K}$ as $\into{\{i\}}{K}$. We call the elements of $|K|$ \emph{points of $K$}, but a more familiar name (to some) might be \emph{global element}.
\end{ntn}

\begin{dfn}\label{df:opcat} An \textbf{\emph{operator category}} $\Phi$ is an essentially small category that satisfies the following three conditions.
\begin{enumerate}[(\ref{df:opcat}.1)]
\item The category $\Phi$ has a terminal object.
\item For any morphism $\fromto{J}{I}$ of $\Phi$ and for any point $i\in|I|$, there exists a fiber $J_i\coloneq\{i\}\times_IJ$.
\item For any pair of objects $I,J\in\Phi$, the set $\Mor_{\Phi}(J,I)$ is finite.
\end{enumerate}
\end{dfn}

\begin{nul}[A note on terminology] The notion we have defined here is distinct from the notion of a \emph{category of operators} used in the brilliant work of May and Thomason \cite{MR508885}, and it serves a distinct mathematical role. We hope that the obvious similarity in nomenclature will not lead to confusion.
\end{nul}

There are very many interesting examples of operator categories, but for now, let us focus on a small number of these.

\begin{exm}\label{exm:firstopcats} The following categories are operator categories:
\begin{enumerate}[(\ref{exm:firstopcats}.1)]
\item the trivial category $\{1\}$;
\item the category $\OO$ of ordered finite sets; and
\item the category $\FF$ of finite sets.
\end{enumerate}
\end{exm}

\begin{exm}\label{exm:truncation} For any operator category $\Phi$ and for any integer $n\geq1$, write $\Phi_{\leq n}$ for the full subcategory of $\Phi$ spanned by those objects $I\in\Phi$ such that $\#|I|\leq n$. Then the category $\Phi_{\leq n}$ is an operator category as well.
\end{exm}

\begin{exm} Suppose $\Psi$ and $\Phi$ two operator categories. Then we may define a category $\Psi\wr\Phi$ as follows. An object of $\Psi\wr\Phi$ will be a pair $(I,M)$ consisting of an object $I\in\Phi$ and a collection $M=\{M_i\}_{i\in|I|}$ of objects of $\Psi$, indexed by the points of $I$. A morphism $(\eta,\omega)\colon\fromto{(J,N)}{(I,M)}$ of $\Psi\wr\Phi$ consists of a morphism $\eta\colon\fromto{J}{I}$ of $\Phi$ and a collection
\begin{equation*}
\{\omega_j\colon\fromto{N_j}{M_{\eta(j)}}\}_{j\in|J|}
\end{equation*}
of morphisms of $\Psi$, indexed by the points of $J$. Then $\Psi\wr\Phi$ is an operator category. In \S \ref{sect:wreath} we will give a more systematic discussion of this story, and we will show that this \emph{wreath product} of operator categories in fact determines a monoidal structure on the collection of all operator categories.
\end{exm}

\begin{exm}\label{exm:semidir} Suppose $\Phi$ any operator category. Then we may define a \textbf{\emph{semidirect product}} category $\Phi\rtimes\OO$ as follows. An object of $\Phi\rtimes\OO$ will be a pair $(I,M)$ consisting of a totally ordered finite set $I\in\OO$ and a functor $M\colon\fromto{I}{\Phi}$. A morphism $(\eta,\omega)\colon\fromto{(J,N)}{(I,M)}$ of $\Phi\rtimes\OO$ consists of a morphism $\eta\colon\fromto{J}{I}$ of $\OO$ and a natural transformation
\begin{equation*}
\omega\colon\fromto{N}{\eta^{\star}M}
\end{equation*}
of functors $\fromto{J}{\Phi}$. Then $\Phi\rtimes\OO$ is an operator category.
\end{exm}

\begin{exm}\label{exm:cyclic} Denote by $\CC$ the category of finite cyclically ordered sets and monotone maps, i.e., maps $\phi\colon\fromto{J}{I}$ such that for any $r,s,t\in J$, if $[\phi(r),\phi(s),\phi(t)]$ in $I$, then $[r,s,t]$ in $J$. Then $\CC$ is an operator category. 
\end{exm}

\begin{exm}\label{exm:graph} Denote by $\GG$ the category of finite simple graphs. Then $\GG$ is an operator category.
\end{exm}

\begin{dfn} A functor $G\colon\fromto{\Psi}{\Phi}$ between operator categories will be said to be \textbf{\emph{admissible}} if it preserves terminal objects and the formation of fibers. An admissible functor $G\colon\fromto{\Psi}{\Phi}$ will be said to be an \textbf{\emph{operator morphism}} if in addition, for any object $I$ of $\Psi$, the induced morphism $\fromto{|I|}{|GI|}$ is a surjection. 
\end{dfn}

\begin{exm}\label{exm:admissopmorph}
\begin{enumerate}[(\ref{exm:admissopmorph}.1)]
\item Any equivalence $\equivto{\Psi}{\Phi}$ between operator categories is an operator morphism.
\item For any operator category $\Phi$, the assignment $\goesto{I}{|I|}$ is an operator morphism $u\colon\fromto{\Phi}{\FF}$ is an operator morphism.
\item For any operator category $\Phi$, the inclusion $\into{\{1\}}{\Phi}$ of the terminal object is an operator morphism.
\item For any operator category $\Phi$ and for any positive integer $n$, the inclusion $\into{\Phi_{\leq n}}{\Phi}$ is an operator morphism.
\item For any operator category $\Phi$, the unique functor $\fromto{\Phi}{\{1\}}$ is an admissible functor, but it is generally not an operator morphism.
\end{enumerate}
\end{exm}

\begin{prp}\label{lem:operatormorphism} Suppose $G\colon\fromto{\Psi}{\Phi}$ an operator morphism. Then for any object $I$ of $\Psi$, the induced map $\fromto{|I|}{|GI|}$ is a bijection.
\begin{proof} Two points $i,j\in|I|$ are distinct if and only if the underlying set $|j^{-1}\{i\}|\cong|j|^{-1}\{i\}$ of the fiber of one point along the other is empty -- or, equivalently, if and only if the intersection $i\cap j$ has no points, that is, $|i\cap j|$ is empty. Since $G$ and $\goesto{I}{|I|}$ are admissible, one also has a bijection $|G(j^{-1}\{i\})|\cong|G(j)|^{-1}\{G(i)\}$. Since $G$ is an operator morphism, the map
\[
\fromto{|j^{-1}\{i\}|}{|G(j^{-1}\{i\})|\cong|G(j)|^{-1}\{G(i)\}}
\]
is surjective, so if the source is empty, then the target is as well.
\end{proof}
\end{prp}

\begin{cor}\label{lem:2outof3operator} Suppose 
\begin{equation*}
\begin{tikzpicture} 
\matrix(m)[matrix of math nodes, 
row sep=6ex, column sep=4ex, 
text height=1.5ex, text depth=0.25ex] 
{&\Psi&\\ 
X&&\Phi\\}; 
\path[>=stealth,->,font=\scriptsize] 
(m-2-1) edge node[above left]{$G$} (m-1-2) 
edge node[below]{$H$} (m-2-3)
(m-1-2) edge node[above right]{$F$} (m-2-3); 
\end{tikzpicture}
\end{equation*}
a commutative triangle of admissible functors in which $F$ is an operator morphism; then $G$ is an operator morphism if and only if $H$ is.
\end{cor}

We organize the collection of operator categories into an $\infty$-category.
\begin{ntn} Denote by $\mathrm{Adm}$ the (strict) $2$-category in which the objects are small operator categories, the $1$-morphisms are admissible functors, and the $2$-morphisms are isomorphisms of functors. Denote by $\mathrm{Op}$ the sub-$2$-category of $\mathrm{Adm}$ in which the objects are small operator categories, the $1$-morphisms are operator morphisms, and the $2$-morphisms are isomorphisms of functors.

Applying the nerve to each $\Mor$-groupoid in $\mathrm{Adm}$ and $\mathrm{Op}$, we obtain categories enriched in fibrant simplicial sets, and we may apply the simplicial nerve to obtain $\infty$-categories that are $2$-categories in the sense of \cite[\S 2.3.4]{HTT} (which could perhaps more precisely be called ``$(2,1)$-categories''). We will refer to these $\infty$-categories as $\Adm$ and $\Op$. 
\end{ntn}

As a result of Pr. \ref{lem:operatormorphism}, we have the following.
\begin{prp} The trivial operator category $\{1\}$ is initial in both $\Adm$ and $\Op$, and it is terminal in $\Adm$. The operator category $\FF$ is terminal in $\Op$.
\end{prp}


\section{Complete Segal operads}\label{sect:wkoperads} Any operator category gives rise to a theory of operads (elsewhere called a colored operad or multicategory). Here we define a \emph{weak} version of this theory, as well as its theory of algebras. To this end, we first single out an important class of morphisms of an operator category.

\begin{dfn}\label{dfn:Miller} A morphism $\fromto{K}{J}$ of an operator category $\Phi$ is a \textbf{\emph{fiber inclusion}} if there exists a morphism $\fromto{J}{I}$ and a point $i\in|I|$ such that the square
\begin{equation*}
\begin{tikzpicture} 
\matrix(m)[matrix of math nodes, 
row sep=4ex, column sep=4ex, 
text height=1.5ex, text depth=0.25ex] 
{K&J\\ 
\{i\}&I\\}; 
\path[>=stealth,->,font=\scriptsize] 
(m-1-1) edge (m-1-2) 
edge (m-2-1) 
(m-1-2) edge (m-2-2) 
(m-2-1) edge[right hook->] (m-2-2); 
\end{tikzpicture}
\end{equation*}
is a pullback square. A morphism $\fromto{K}{J}$ is an \textbf{\emph{interval inclusion}} if it can be written as the composite of a finite sequence of fiber inclusions. An interval inclusion will be denoted by a hooked arrow $\into{K}{J}$.
\end{dfn}

\begin{exm}\label{exm:intinc}
\begin{enumerate}[(\ref{exm:intinc}.1)]
\item Any isomorphism of an operator category $\Phi$ is a fiber inclusion, and any point of any object of an operator category $\Phi$ is a fiber inclusion.
\item Suppose $I$ an object of $\mathbf{O}$. Then for any elements $i_0,i_1\in|I|$, write $[i_0,i_1]$ for the subset of elements $i\in|I|$ such that $i_0\leq i\leq i_1$, with the induced ordering. A morphism $\fromto{J}{I}$ is a fiber inclusion if and only if it is a monomorphism whose image is precisely $[i_0,i_1]$ for some elements $i_0,i_1\in|I|$. Interval inclusions and fibers inclusions coincide in $\mathbf{O}$.
\item In the category $\mathbf{F}$, fiber inclusions are precisely monomorphisms. Again interval inclusions and fibers inclusions coincide.
\end{enumerate}
\end{exm}

\begin{nul}\label{nul:intervinc} We make the following observations.
\begin{enumerate}[(\ref{nul:intervinc}.1)]
\item Interval inclusions in an operator category are monomorphisms.
\item\label{lem:intincclosedpullback} If $\Phi$ is an operator category, then the pullback of any morphism $\fromto{L}{J}$ along any interval inclusion $\into{K}{J}$ exists, and the canonical morphism $\fromto{K\times_JL}{L}$ is an interval inclusion.
\item Suppose $\Phi$ an operator category; suppose $\phi\colon\into{K}{J}$ an interval inclusion in $\Phi$; and suppose $\psi\colon\fromto{L}{K}$ a morphism of $\Phi$. Then $\psi$ is an interval inclusion if and only if $\phi\circ\psi$ is.
\item\label{lem:emptytermdilemma} Suppose $\Phi$ an operator category, $\into{K}{J}$ an interval inclusion, and $j\in|J|$ a point. Then either the fiber $K_j$ is a terminal object or else $|K_j|=\varnothing$.
\end{enumerate}
\end{nul}

\begin{dfn} Suppose $\Phi$ an operator category. Then a \textbf{\emph{$\Phi$-sequence}} is a pair $(\mathbf{m},I)$ consisting of an object $m\in\Delta$ and a functor $I\colon\fromto{\mathbf{m}}{\Phi}$. We will denote such an object by
\begin{equation*}
[I_0\to I_1\to\cdots\to I_m].
\end{equation*}

A \textbf{\emph{morphism $(\eta,\phi)\colon\fromto{(\mathbf{n},J)}{(\mathbf{m},I)}$ of $\Phi$-sequences}} consists of a morphism $\eta\colon\fromto{\mathbf{n}}{\mathbf{m}}$ of $\Delta$ and a natural transformation $\phi\colon\fromto{J}{I\circ\eta}$ such that for any integer $0\leq k\leq m$, the morphism $\phi_k\colon\fromto{J_k}{I_{\eta(k)}}$ is an interval inclusion, and for any pair of integers $0\leq k\leq\ell\leq m$, the square
\begin{equation*}
\begin{tikzpicture} 
\matrix(m)[matrix of math nodes, 
row sep=4ex, column sep=4ex, 
text height=1.5ex, text depth=0.25ex] 
{J_{k}&I_{{\eta(k)}}\\ 
J_{\ell}&I_{{\eta(\ell)}}\\}; 
\path[>=stealth,->,font=\scriptsize] 
(m-1-1) edge[right hook->] node[above]{$\phi_{k}$} (m-1-2) 
edge (m-2-1) 
(m-1-2) edge (m-2-2) 
(m-2-1) edge[right hook->] node[below]{$\phi_{\ell}$} (m-2-2); 
\end{tikzpicture}
\end{equation*}
is a pullback of $\Phi$. Denote by $\Delta_{\Phi}$ the category of $\Phi$-sequences.

Any admissible functor $G\colon\fromto{\Psi}{\Phi}$ induces a functor $\fromto{\Delta_{\Psi}}{\Delta_{\Phi}}$ given by $\goesto{(\mathbf{m},I)}{(\mathbf{m},G\circ I)}$.
\end{dfn}

\begin{ntn}\label{ntn:Segcondsonseqs} Suppose $\Phi$ an operator category, and suppose $X\colon\fromto{\Delta_{\Phi}^{\op}}{\Kan}$ a functor. We study three classes of maps.
\begin{enumerate}[(\ref{ntn:Segcondsonseqs}.1)]
\item For any $\Phi$-sequence $(m,I)$ and for any point $i\in|I_m|$, one obtains a map
\begin{equation*}
\fromto{X[I_0\to I_1\to\cdots\to I_m]}{X[I_{0,i}\to I_{1,i}\to\cdots\to\{i\}]}.
\end{equation*}
Consequently, one obtains a map
\begin{equation*}
p_{(\mathbf{m},I)}\colon\fromto{X[I_0\to I_1\to\cdots\to I_m]}{\prod_{i\in|I_m|}X[I_{0,i}\to I_{1,i}\to\cdots\to\{i\}]}.
\end{equation*}
\item For any $\Phi$-sequence $(\mathbf{m},I)$, one obtains a map
\begin{equation*}
s_{\mathbf{m},I}\colon\fromto{X[I_0\to\cdots\to I_m]}{X[I_0\to I_1]\times^{h}_{X[I_1]}\cdots\times^{h}_{[I_{m-1}]}X[I_{m-1}\to I_m]},
\end{equation*}
where the target is the homotopy fiber product.
\item Lastly, the inclusion $\into{\{1\}}{\Phi}$ induces an inclusion $\into{\Delta\cong\Delta_{\{1\}}}{\Delta_{\Phi}}$; hence $X$ restricts to a simplicial space $(X|\Delta^{\op})$, and one can define the map 
\begin{equation*}
r\colon\fromto{(X|\Delta^{\op})_0}{(X|\Delta^{\op})_K},
\end{equation*}
where $K=\Delta^3/(\Delta^{\{0,2\}}\sqcup\Delta^{\{1,3\}})$, and $(X|\Delta^{\op})_K$ is the homotopy limit of the diagram
\begin{equation*}
\Delta_{/K}^{\op}\to\Delta^{\op}\ \tikz[baseline]\draw[>=stealth,->,font=\scriptsize](0,0.5ex)--node[above]{$(X|\Delta^{\op})$}(1,0.5ex);\ s\Set.
\end{equation*}
\end{enumerate}
\end{ntn}

\begin{dfn} Suppose $\Phi$ an operator category. A \textbf{\emph{complete Segal $\Phi$-operad}} is a left fibration $q\colon\fromto{X}{N\Delta_{\Phi}^{\op}}$ such that any functor $\fromto{\Delta_{\Phi}^{\op}}{\Kan}$ that classifies $q$ has the property that the maps $p_{(\mathbf{m},I)}$, $s_{(\mathbf{m},I)}$, and $r$ (Nt. \ref{ntn:Segcondsonseqs}) are all equivalences. A \textbf{\emph{morphism $g\colon\fromto{X}{Y}$ of complete Segal $\Phi$-operads}} is a commutative diagram
\begin{equation*}
\begin{tikzpicture} 
\matrix(m)[matrix of math nodes, 
row sep=4ex, column sep=4ex, 
text height=1.5ex, text depth=0.25ex] 
{X&&Y\\ 
&N\Delta_{\Phi}^{\op}&\\}; 
\path[>=stealth,->,font=\scriptsize] 
(m-1-1) edge node[above]{$g$} (m-1-3) 
edge (m-2-2) 
(m-1-3) edge (m-2-2); 
\end{tikzpicture}
\end{equation*}
A morphism $g\colon\fromto{X}{Y}$ will be said to be a \textbf{\emph{equivalence of complete Segal $\Phi$-operads}} if it is a covariant weak equivalence.

Denote by $\Operad_{\mathrm{CSS}}^{\Phi,\Delta}$ the full simplicial subcategory of $s\Set_{/N\Delta^{\op}_{\Phi}}$ spanned by the complete Segal $\Phi$-operads. This is a fibrant simplicial category, so we can define an $\infty$-category $\Operad_{\mathrm{CSS}}^{\Phi}$ as the simplicial nerve of $\Operad_{\mathrm{CSS}}^{\Phi}$.
\end{dfn}

\begin{nul} Suppose $\Phi$ an operator category. We refer to the fiber of a complete Segal $\Phi$-operad $X$ over an object
\begin{equation*}
[I_0\to I_1\to\cdots\to I_m]\in\Delta^{\op}_{\Phi}
\end{equation*}
as the \textbf{\emph{space of operations of type}} $[I_0\to I_1\to\cdots\to I_m]$. The space of operations of type $\{1\}$ is the \textbf{\emph{space of objects}} of $X$ or \emph{interior $\infty$-groupoid} of $X$, and the space of operations of type $I$ is the \textbf{\emph{space of $I$-tuples of objects}} of $X$. Let us denote by
\begin{equation*}
[I_0/I_1/\cdots/I_m]\Map_X((x^0_{i_0})_{i_0\in|I_0|},\dots,(x^m_{i_m})_{i_m\in|I_m|})
\end{equation*}
the fiber of the map $X[I_0\to I_1\to\cdots\to I_m]\to X(I_0)\times\cdots\times X(I_m)$ over the vertex $((x^0_{i_0})_{i_0\in|I_0|},\dots,(x^m_{i_m})_{i_m\in|I_m|})$. The condition that $X$ is a complete Segal $\Phi$-operad then gives equivalences
\begin{equation*}\notag\begin{split}
[J/I/\{1\}]\Map_X((x_j)_{j\in|J|},(y_i)_{i\in|I|},z)&\\
&\!\!\!\!\!\!\!\!\!\!\!\!\!\!\!\!\!\!\!\!\!\!\!\!\!\!\!\!\!\!\!\!\!\!\!\!\!\!\!\!\!\!\!\!\!\!\!\!\simeq[I/\{1\}]\Map_X((y_i)_{i\in|I|},z)\times[J/I]\Map_X((x_j)_{j\in|J|},(y_i)_{i\in|I|})\\
&\!\!\!\!\!\!\!\!\!\!\!\!\!\!\!\!\!\!\!\!\!\!\!\!\!\!\!\!\!\!\!\!\!\!\!\!\!\!\!\!\!\!\!\!\!\!\!\!\!\!\simeq[I/\{1\}]\Map_X((y_i)_{i\in|I|},z)\times\prod_{i\in|I|}[J_i/\{i\}]\Map_X((x_j)_{j\in|J_i|},y_i),
\end{split}
\end{equation*}
and the map
\begin{equation*}
\fromto{[J/I/\{1\}]\Map_X((x_j)_{j\in|J|},(y_i)_{i\in|I|},z)}{[J/\{1\}]\Map_X((x_j)_{j\in|J|},z)}
\end{equation*}
induced by the map
\begin{equation*}
\begin{tikzpicture} 
\matrix(m)[matrix of math nodes, 
row sep=4ex, column sep=4ex, 
text height=1.5ex, text depth=0.25ex] 
{J&&\{1\}\\ 
J&I&\{1\}\\}; 
\path[>=stealth,->,font=\scriptsize] 
(m-2-1) edge (m-2-2) 
edge[-,double distance=1.5pt] (m-1-1) 
(m-2-2) edge (m-2-3)
(m-2-3) edge[-,double distance=1.5pt] (m-1-3) 
(m-1-1) edge (m-1-3); 
\end{tikzpicture}
\end{equation*}
amounts to a \textbf{\emph{polycomposition map}}, which is defined up to coherent homotopy. The functoriality in $N\Delta^{\op}_{\Phi}$ amounts to a coherent associativity condition.
\end{nul}

We immediately obtain the following characterization of equivalences between complete Segal $\Phi$-operads.
\begin{prp}\label{prp:DKequivs} Suppose $\Phi$ an operator category. Then a morphism $g\colon\fromto{X}{Y}$ of complete Segal $\Phi$-operads is an equivalence just in case the following conditions are satisfied.
\begin{enumerate}[(\ref{prp:DKequivs}.1)]
\item \emph{Essential surjectivity}. The map $\fromto{\pi_0X[\{1\}]}{\pi_0Y[\{1\}]}$ is surjective.
\item \emph{Full faithfulness}. For any object $I\in\Phi$, any vertex $x\in X[I]_{0}$, and any vertex $y\in X[\{1\}]_{0}$, the induced map
\begin{equation*}
\fromto{[I/\{1\}]\Map_{X}(x,y)}{[I/\{1\}]\Map_{Y}(g(x),g(y))}
\end{equation*}
is an equivalence.
\end{enumerate}
\end{prp}

\begin{exm} For any operator category $\Phi$, the identity functor on the simplicial set $N\Delta_{\Phi}^{\op}$ is a complete Segal $\Phi$-operad --- the \textbf{\emph{terminal complete Segal $\Phi$-operad}}, which we denote $U_{\Phi}$. These complete Segal $\Phi$-operads, for suitable choice of $\Phi$, give rise to all the operads discussed in the introduction. In particular, when $\Phi=\FF$, we show that the terminal complete Segal $\FF$-operad $U_{\FF}$ is equivalent to the operad $E_{\infty}$ (\S \ref{sect:examples}).
\end{exm}

\begin{exm} When $\Phi=\{1\}$, we find that a complete Segal $\{1\}$-operad is a left fibration $\fromto{X}{N\Delta^{\op}}$ classified by a complete Segal \emph{space} in the sense of Rezk \cite{MR1804411}. (In particular, $\Operad_{\mathrm{CSS}}^{\{1\}}$ is a homotopy theory of $(\infty,1)$-categories, in the sense of \cite{BSP}.) For simplicity, we will, by a small abuse, simply call such left fibrations \textbf{\emph{complete Segal spaces}}.
\end{exm}

\begin{exm} When $\Phi=\FF$, we obtain a new homotopy theory of weak symmetric operads, which we will show is equivalent to Lurie's in \S \ref{sect:equivalence}.
\end{exm}

Using the left Bousfield localization \cite{MR2003j:18018, MR2771591} of the covariant model structure \cite[Pr. 2.1.4.7]{HTT} with respect to the set $S_{\Phi}$ of those morphisms that represent the morphisms of \ref{ntn:Segcondsonseqs}, one obtains the following.
\begin{prp}\label{prp:wkoperadmodstruct} The category $s\Set_{/N\Delta_{\Phi}^{\op}}$ admits a left proper, tractable, simplicial model structure --- called the \textbf{\emph{operadic model structure}} --- with the following properties.
\begin{enumerate}[(\ref{prp:wkoperadmodstruct}.1)]
\item A map $\fromto{X}{Y}$ over $N\Delta_{\Phi}^{\op}$ is a cofibration just in case it is a monomorphism.
\item An object $\fromto{X}{N\Delta_{\Phi}^{\op}}$ is fibrant just in case it is a complete Segal $\Phi$-operad.
\item A map $\fromto{X}{Y}$ of simplicial sets over $N\Delta^{\op}_{\Phi}$ is a weak equivalence if and only if, for any complete Segal $\Phi$-operad $Z$, the induced map
\begin{equation*}
\fromto{\Map_{/N\Delta^{\op}_{\Phi}}(Y,Z)}{\Map_{/N\Delta^{\op}_{\Phi}}(X,Z)}
\end{equation*}
is a weak equivalence.
\item A map $\fromto{X}{Y}$ between complete Segal $\Phi$-operads is a weak equivalence just in case it is an equivalence of complete Segal $\Phi$-operads.
\end{enumerate}
\end{prp}
\begin{cor} The $\infty$-category $\Operad_{\mathrm{CSS}}^{\Phi}$ of complete Segal $\Phi$-operads is an accessible localization of the functor $\infty$-category $\Fun(N\Delta^{\op}_{\Phi},\Kan)$; in particular, it is a presentable $\infty$-category.
\end{cor}

Critically, the homotopy theory of weak operads over operator categories is functorial with respect to operator morphisms.
\begin{prp}\label{prp:operadicfunctorial} For any operator morphism $G\colon\fromto{\Psi}{\Phi}$, the adjunction
\begin{equation*}
\adjunct{G_!}{s\Set_{/N\Delta_{\Psi}^{\op}}}{s\Set_{/N\Delta_{\Phi}^{\op}}}{G^{\star}}
\end{equation*}
is a Quillen adjunction for the operadic model structure.
\begin{proof} It is enough to note that the functor $G^{\star}$ is a right adjoint for the covariant model structure, and a left fibration $\fromto{X}{N\Delta_{\Phi}^{\op}}$ is a complete Segal $\Phi$-operad only if $\fromto{G^{\star}X\cong X\times_{N\Delta_{\Phi}^{\op}}N\Delta_{\Psi}^{\op}}{N\Delta_{\Psi}^{\op}}$ is a complete Segal $\Psi$-operad.
\end{proof}
\end{prp}

\begin{cor} An operator morphism $G\colon\fromto{\Psi}{\Phi}$ induces an adjunction
\begin{equation*}
\adjunct{G_!}{\Operad_{\mathrm{CSS}}^{\Psi}}{\Operad_{\mathrm{CSS}}^{\Phi}}{G^{\star}}
\end{equation*}
of $\infty$-categories.
\end{cor}

\begin{exm} For any operator category $\Psi$, the right adjoint
\begin{equation*}
p^{\star}\colon\fromto{\Operad^{\Psi}}{\Operad^{\{1\}}}
\end{equation*}
induced by the inclusion $p\colon\into{\{1\}}{\Psi}$ carries any complete Segal $\Psi$-operad to its \textbf{\emph{underlying complete Segal space}}.
\end{exm}

\begin{exm} When $\Phi=\FF$, the notion of a complete Segal $\Phi$-operad is closely related to the usual notion of a symmetric operad in simplicial sets. For any operator category $\Psi$, the right adjoint $u^{\star}\colon\fromto{\Operad_{\mathrm{CSS}}^{\FF}}{\Operad_{\mathrm{CSS}}^{\Psi}}$ induced by the essentially unique operator morphism $u\colon\fromto{\Psi}{\FF}$ carries any complete Segal $\FF$-operad to its \textbf{\emph{underlying complete Segal $\Psi$-operad}}. The left adjoint $\mathbf{Symm}\colon u_!$ carries any complete Segal $\Psi$-operad to its \textbf{\emph{symmetrization}}.
\end{exm}

\begin{exm} When $\Psi=\OO$, we show below that the symmetrization of the terminal complete Segal $\OO$-operad $U_{\OO}$ is equivalent to the operad $E_1$, or, equivalently, the operad $A_{\infty}$. When $\Psi=\OO_{\leq n}$ for some integer $n\geq 1$, we show below that the symmetrization of the terminal complete Segal $\OO_{\leq n}$-operad $U_{\OO_{\leq n}}$ is equivalent to the operad $A_n$ (\S \ref{sect:examples}).
\end{exm}

We note that it will sometimes be convenient to work with a straightened variant of the operadic model structure of Pr. \ref{prp:wkoperadmodstruct}, which is provided by the following left Bousfield localization of the injective model structure.
\begin{prp}\label{prp:straightenedoperads} The category $\Fun(\Delta_{\Phi}^{\op},s\Set)$ admits a left proper tractable model structure --- called the \textbf{\emph{operadic model structure}} --- with the following properties.
\begin{enumerate}[(\ref{prp:straightenedoperads}.1)]
\item A natural transformation $\fromto{X}{Y}$ is a cofibration just in case it is a monomorphism.
\item An object $X$ is fibrant just in case it is valued in Kan complexes, and it classifies a complete Segal $\Phi$-operad.
\item A natural transformation $\fromto{X}{Y}$ is a weak equivalence if and only if, for any fibrant object $Z$, the induced map
\begin{equation*}
\fromto{\Map(Y,Z)}{\Map(X,Z)}
\end{equation*}
is a weak equivalence.
\end{enumerate}
\end{prp}

\begin{nul} For any operator category $\Phi$, it is clear that the straightening/unstraightening Quillen equivalence of \cite[\S 2.2.1]{HTT} localizes to a Quillen equivalence
\begin{equation*}
\adjunct{\mathrm{St}}{s\Set_{/N\Delta^{\op}}}{\Fun(\Delta_{\Phi}^{\op},s\Set)}{\mathrm{Un}}
\end{equation*}
between the operadic model structures on $\Fun(\Delta_{\Phi}^{\op},s\Set)$ and $s\Set_{/N\Delta_{\Phi}^{\op}}$.
\end{nul}


\section{Wreath products}\label{sect:wreath} Roughly speaking, if the objects of an operator category index a certain sort of multiplications, then the objects of a wreath product of two operator categories $\Phi$ and $\Psi$ index $\Phi$-multiplications in objects that already possess $\Psi$-multiplications. This provides a wealth of new examples of operator categories, and with this insight, we can introduce the operator categories $\OO^{(n)}$ ($1\leq n<\infty$), which are key to the combinatorial gadgets that characterize the operads $E_n$. 

\begin{dfn} Suppose $\Phi$ an operator category. Then a \textbf{\emph{coronal fibration}}
\begin{equation*}
p\colon\fromto{X}{N\Phi}
\end{equation*}
is a cartesian fibration such that, for any object $I\in\Phi$, the functors
\begin{equation*}
\{\fromto{X_I}{X_{\{i\}}}\ |\ i\in|I|\}
\end{equation*}
together exhibit the fiber $X_I$ as a product of the fibers $X_{\{i\}}$. In this situation, $p$ will be said to \textbf{\emph{exhibit $X$ as a wreath product of $X_1$ with $\Phi$}}.
\end{dfn}

\begin{ntn} Denote by $\mathscr{O}^{\cart}(\Cat_{\infty})$ the subcategory \cite[\S 1.2.11]{HTT} of the $\infty$-category
\begin{equation*}
\mathscr{O}(\Cat_{\infty})\coloneq\Fun(\Delta^1,\Cat_{\infty})
\end{equation*}
whose objects are cartesian fibrations and whose morphisms carry cartesian morphisms to cartesian morphisms. Let $\Cat_{\infty/S}^{\cart}$ denote the fiber of the target functor
\begin{equation*}
t\colon\fromto{\mathscr{O}^{\cart}(\Cat_{\infty})}{\Cat_{\infty}}
\end{equation*}
over an object $\{S\}\subset\Cat_{\infty}$. Note that $\Cat_{\infty/S}^{\cart}$ may be identified with the nerve of the cartesian simplicial model category of marked simplicial sets over $S$ \cite[Pr. 3.1.3.7]{HTT}, whence it is the relative nerve of the category of cartesian fibrations over $S$, equipped with the cartesian equivalences.

Now for any operator category $\Phi$, denote by
\begin{equation*}
\Cor_{/N\Phi}\subset\Cat_{\infty,/N\Phi}^{\cart}
\end{equation*}
the full subcategory spanned by the coronal fibrations.
\end{ntn}

\begin{lem}\label{lem:CoroverNPhi} Suppose $\Phi$ an operator category with a terminal object $1$. Then the functor $\fromto{\Cor_{/N\Phi}}{\Cat_{\infty}}$ given by the extraction $\goesto{X}{X_1}$ of the fiber over $1$ is a trivial fibration.
\begin{proof} It is a direct consequence of straightening \cite[\S 3.2]{HTT} that the $\infty$-category $\Cat_{\infty,/N\Phi}^{\cart}$ is equivalent to the functor $\infty$-category $\Fun(N\Phi^{\op},\Cat_{\infty})$, and the functor $\fromto{\Cat_{\infty,/N\Phi}^{\cart}}{\Cat_{\infty}}$ given by the assignment $\goesto{X}{X_{1}}$ corresponds to evaluation at $1\in N\Phi^{\op}$. Now the result follows from the observation that a cartesian fibration $p\colon\fromto{X}{N\Phi}$ is a coronal fibration just in case the corresponding functor $F_p\colon\fromto{N\Phi^{\op}}{\Cat_{\infty}}$ is a right Kan extension of the restriction of $F_p$ to $\{1\}$.
\end{proof}
\end{lem}

\begin{ntn} For any $\infty$-category $C$ and any operator category $\Phi$, we write $C\wr\Phi$ for the total space $X$ of a coronal fibration $p\colon\fromto{X}{N\Phi}$ such that $C$ appears as the fiber over a terminal object. One can make the assignment $\fromto{C}{C\wr\Phi}$ into a functor by choosing a terminal object $1\in\Phi$ and a section of the trivial fibration $\fromto{\Cor_{/N\Phi}}{\Cat_{\infty}}$ informally described as the assignment $\goesto{X}{X_{1}}$. The space of these choices is contractible.
\end{ntn}

\begin{lem} If $F\colon\fromto{\Psi}{\Phi}$ is an operator morphism, then the pullback
\begin{equation*}
F^{\star}p\colon\fromto{X\times_{N\Phi}N\Psi}{N\Psi}
\end{equation*}
of a coronal fibration $p\colon\fromto{X}{N\Phi}$ is a coronal fibration.
\begin{proof} Suppose $I\in\Psi$. The maps $\fromto{(X\times_{N\Phi}N\Psi)_I}{(X\times_{N\Phi}N\Psi)_{\{i\}}}$ for $i\in|I|$ can be identified with the maps $\fromto{X_{FI}}{X_{F\{i\}}}$. Since the map $\fromto{|I|}{|FI|}$ is a bijection, these maps exhibit the fiber $(X\times_{N\Phi}N\Psi)_I$ as a product of the fibers $(X\times_{N\Phi}N\Psi)_{\{i\}}$.
\end{proof}
\end{lem}

This lemma, combined with the previous one, implies that pullback along an operator morphism $F\colon\fromto{\Psi}{\Phi}$ induces an equivalence of $\infty$-categories
\begin{equation*}
\equivto{\Cor_{/N\Phi}}{\Cor_{/N\Psi}};
\end{equation*}
they also imply that the whole $\infty$-category of coronal fibrations is equivalent to the $\infty$-category $\Op\times\Cat_{\infty}$. We will return momentarily to a special case of this observation.

\begin{lem} Suppose $\Phi$ an operator category, and suppose $p\colon\fromto{X}{N\Phi}$ a coronal fibration such that the fiber $X_{1}$ is the nerve of a operator category. Then $X$ itself is also the nerve of an operator category, and $p$ is the nerve of admissible functor.
\begin{proof} Since $X_{1}$ is a $1$-category, every fiber $X_I\simeq X_{1}^{\times|I|}$ is a $1$-category. It follows that $X$ itself is a $1$-category. The fact that $X$ has a terminal object and all fibers follows directly from \cite[Pr. 4.3.1.10]{HTT}.
\end{proof}
\end{lem}
\noindent Informally, we conclude that the wreath product of two operator categories is again an operator category.

\begin{lem} Suppose $\Phi$ and $\Phi'$ operator categories, and suppose $p\colon\fromto{\Phi'}{\Phi}$ an admissible functor whose nerve is a coronal fibration. Then for any operator morphism $F\colon\fromto{\Psi}{\Phi}$ and any pullback diagram
\begin{equation*}
\begin{tikzpicture} 
\matrix(m)[matrix of math nodes, 
row sep=4ex, column sep=4ex, 
text height=1.5ex, text depth=0.25ex] 
{\Psi'&\Phi'\\ 
\Psi&\Phi,\\}; 
\path[>=stealth,->,font=\scriptsize] 
(m-1-1) edge node[above]{$F'$} (m-1-2) 
edge node[left]{$q$} (m-2-1) 
(m-1-2) edge node[right]{$p$} (m-2-2) 
(m-2-1) edge node[below]{$F$} (m-2-2); 
\end{tikzpicture}
\end{equation*}
the functor $q$ is admissible, and the functor $F'$ is an operator morphism.
\begin{proof} It is a simple matter to see that $q$ and $F'$ are admissible functors. To check that $F'$ is an operator morphism, let us note first that $F'$ induces an equivalence $\equivto{\Psi'_{\ast}}{\Phi'_{\ast}}$; hence for any object $K\in\Psi'_{\ast}$, the natural map $\fromto{|K|}{|F'K|}$ is a bijection. Now for any object $J\in\Psi$, the decomposition $\Psi'_J\simeq(\Psi'_{\star})^{\times|J|}$ gives, for any object $K$ of $\Psi'$ lying over $J$, a decomposition of finite sets
\begin{equation*}
|K|\cong\coprod_{j\in|J|}|K_{\{j\}}|.
\end{equation*}
After applying $F'$, we similarly obtain a decomposition of finite sets
\begin{equation*}
|F'K|\cong\coprod_{i\in|FJ|}|(F'K)_{\{i\}}|.
\end{equation*}
We thus conclude that the map $\fromto{|K|}{|F'K|}$ is a bijection since both the map $\fromto{|J|}{|FJ|}$ and all the maps $\fromto{|K_{\{j\}}|}{|K_{F\{j\}}|}$ are bijections.
\end{proof}
\end{lem}

We now set about showing that the wreath product defines a monoidal structure on the $\infty$-category $\Op$.
\begin{ntn} Denote by $\mathrm{M}$ the ordinary category whose objects are pairs $(m,i)$ consisting of an integer $m\geq 0$ and an integer $0\leq i\leq m$ and whose morphisms $\fromto{(n,j)}{(m,i)}$ are maps $\phi\colon\fromto{\mathbf{m}}{\mathbf{n}}$ of $\Delta$ such that $j\leq\phi(i)$. This category comes equipped with a natural projection $\fromto{\mathrm{M}}{\Delta^{\op}}$.

Denote by $E(\Adm)$ the simplicial set over $N\Delta$ specified by the following universal property. We require, for any simplicial set $K$ and any map $\sigma\colon\fromto{K}{N\Delta^{\op}}$, a bijection
\begin{equation*}
\Mor_{/N\Delta^{\op}}(K,E(\Adm)^{\op})\cong\Mor(K\times_{N\Delta^{\op}}N\mathrm{M},\Adm^{\op}),
\end{equation*}
functorial in $\sigma$. Now by \cite[Cor. 3.2.2.13]{HTT}, the map $\fromto{E(\Adm)}{N\Delta}$ is a cartesian fibration, and $E(\Adm)$ is an $\infty$-category whose objects are pairs $(m,X)$ consisting of an integer $m\geq 0$ and a functor $X\colon\fromto{(\Delta^m)^{\op}}{\Adm}$. A morphism $(\phi,g)\colon\fromto{(n,Y)}{(m,X)}$ can be regarded as a map $\phi\colon\fromto{\mathbf{n}}{\mathbf{m}}$ of $\Delta$ and an edge $g\colon\fromto{Y}{\phi^{\star}X}$ in $\Fun((\Delta^n)^{\op},\Adm)$.

Denote by $\Op^{\wr}$ the following subcategory of $E(\Adm)$. The objects of $\Op^{\wr}$ are those pairs $(m,X)\in E(\Adm)$ such that for any integer $1\leq i\leq m$, the nerve of the admissible functor $\fromto{X_i}{X_{i-1}}$ is a coronal fibration, and the operator category $X_0$ is equivalent to $\{\ast\}$. The morphisms $\fromto{(n,Y)}{(m,X)}$ are those pairs $(\phi,g)$ such that for any integer $0\leq i\leq n$, the admissible functor $\fromto{Y_i}{X_{\phi(i)}}$ is an operator morphism.
\end{ntn}

\begin{prp}\label{prp:Opmonoidal} The inner fibration $\fromto{(\Op^{\wr})^{\op}}{N\Delta^{\op}}$ is a monoidal $\infty$-category.
\begin{proof} We first show that $\fromto{\Op^{\wr}}{N\Delta}$ is a cartesian fibration. Indeed, for any object $(m,X)$ of $\Op^{\wr}$ and any edge $\phi\colon\fromto{\mathbf{n}}{\mathbf{m}}$ of $N\Delta$, there exists a morphism $\fromto{Z}{\phi^{\star}X}$ in which $Z\colon\fromto{(\Delta^n)^{\op}}{\Adm}$ is a diagram such that $Z_0$ is equivalent to $\{1\}$ and, for any integers $0\leq i\leq j\leq n$, the square
\begin{equation*}
\begin{tikzpicture} 
\matrix(m)[matrix of math nodes, 
row sep=4ex, column sep=4ex, 
text height=1.5ex, text depth=0.25ex] 
{Z_{j}&X_{\phi(j)}\\ 
Z_{i}&X_{\phi(i)}\\}; 
\path[>=stealth,->,font=\scriptsize] 
(m-1-1) edge (m-1-2) 
edge (m-2-1) 
(m-1-2) edge (m-2-2) 
(m-2-1) edge (m-2-2); 
\end{tikzpicture}
\end{equation*}
is a pullback square. It is straightforward now to check that the resulting edge $\fromto{Z}{X}$ is cartesian over $\phi$.

It now remains to show that for any integer $m\geq 0$, the natural map
\begin{equation*}
\sigma_m\colon\fromto{\Op^{\wr}_{\mathbf{m}}}{\prod_{i=1}^m\Op^{\wr}_{\{i-1,i\}}}
\end{equation*}
is an equivalence of $\infty$-categories. This is the functor that assigns to any object $(m,X)$ of $\Op^{\wr}_{\mathbf{m}}$ the tuple
\begin{equation*}
(X(m)_1,X(m-1)_1,\dots,X(2)_1,X(1)),
\end{equation*}
where $X(k)_1$ denotes the fiber of $\fromto{X(k)}{X(k-1)}$ over a terminal object of $X(k-1)$. It follows from Lm. \ref{lem:CoroverNPhi} that this functor is indeed an equivalence.
\end{proof}
\end{prp}

\begin{exm} Of course we may form the wreath product of \emph{any} two operator categories, but of particular import are the operator categories obtained by forming the iterated wreath product of $\OO$ with itself:
\begin{equation*}
\OO^{(n)}\coloneq\OO\wr\OO\wr\cdots\wr\OO.
\end{equation*}
As Clemens Berger has observed, the homotopy theory of complete Segal $\OO^{(n)}$-operads is a homotopy-coherent variant of Michael Batanin's notion of an \emph{$(n-1)$-terminal $n$-operad} \cite{MR2365200}. We prove below that the little $n$-cubes operad $E_n$ is equivalent to the symmetrization $\mathbf{Symm}(U_{\OO^{(n)}})$ of the terminal complete Segal $\OO^{(n)}$-operad (\S \ref{sect:examples}). This is a variant of Batanin's result in \cite{MR2365200}.
\end{exm}


\section{Perfect operator categories}\label{sect:perfect}

The May--Thomason category of a symmetric operad \cite{MR508885} is a category that lives over Segal's category $\Gamma^{\op}$ of pointed finite sets. Lurie's theory of $\infty$-operads \cite[Ch. 2]{HA} is built upon a generalization of this picture. How do we understand the relationship between the operator category $\FF$ and the category $\Gamma^{\op}$? For nonsymmetric operads, one has a similar May--Thomason construction over $\Delta^{\op}$ \cite{MR534053}. How do we understand the relationship between the operator category $\OO$ and the category $\Delta^{\op}$? Is it analogous to the relationship between the operator category $\FF$ and the category $\Gamma^{\op}$?

At first blush, the answer appears to be no: $\Gamma^{\op}$ is the category of pointed objects of $\FF$, and $\Delta^{\op}$ is opposite of the category of nonempty objects of $\OO$. Expressed this way, $\Gamma^{\op}$ and $\Delta^{\op}$ don't look the least bit similar. However, an insight that we inherited from Tom Leinster \cite[pp. 40--43]{leinsterhaops} shows us how to think of them each as special cases of a general construction. This insight leads us to study a special class of operator categories, which we call \emph{perfect}. We discuss the key properties of these operator categories here. In the next section we find that perfect operator categories admit a canonical monad, and in the section after that, we use that monad to define what we call the \emph{Leinster category} $\Lambda(\Phi)$ of a perfect operator category, and we show that
\begin{equation*}
\Lambda(\FF)\simeq\Gamma^{\op}\textrm{\quad and\quad}\Lambda(\OO)\simeq\Delta^{\op}.
\end{equation*}
These Leinster categories will be the foundation upon which our analogue of Lurie's theory of $\infty$-operads over operator categories is built.

The first property enjoyed by perfect operator categories is the existence of a \emph{point classifier}. Point classifiers play a role in the theory of perfect operator categories that is in many respects analogous to the role played by the subobject classifier in topos theory.

\begin{dfn}\label{dfn:ptclassify} Suppose $\Phi$ an operator category. Consider the category $\Phi^{\mathrm{cons}}$ whose objects are \emph{points}, -- that is, maps $\fromto{\{i\}}{I}$, which we shall denote $(I,i)$ -- and whose morphisms $\fromto{(I,i)}{(J,j)}$ are pullback squares
\begin{equation*}
\begin{tikzpicture}[baseline]
\matrix(m)[matrix of math nodes, 
row sep=4ex, column sep=4ex, 
text height=1.5ex, text depth=0.25ex] 
{\{i\}&I\\
\{j\}&J.\\}; 
\path[>=stealth,->,font=\scriptsize] 
(m-1-1) edge (m-1-2)
(m-1-1) edge[-,double distance=1.5pt] (m-2-1)
(m-1-2) edge (m-2-2) 
(m-2-1) edge (m-2-2); 
\end{tikzpicture}
\end{equation*}
A \textbf{\emph{point classifier}} for $\Phi$ is a terminal object $T$ of $\Phi^{\mathrm{cons}}$.
\end{dfn}

\noindent If a point classifier in a pointed category exists, then it is essentially unique.

\begin{prp} Suppose $\Phi$ an operator category, and suppose $(T,t)$ a point classifier for $\Phi$. Then there is a functor $\chi\colon\fromto{\Phi^{\mathrm{cons}}}{(\Phi/T)}$ such that the square
\begin{equation*}
\begin{tikzpicture}[baseline]
\matrix(m)[matrix of math nodes, 
row sep=4ex, column sep=4ex, 
text height=1.5ex, text depth=0.25ex] 
{\Phi^{\mathrm{cons}}&(\Phi/T)\\
\{1\}&\Phi\\}; 
\path[>=stealth,->,font=\scriptsize] 
(m-1-1) edge node[above]{$\chi$} (m-1-2)
(m-1-1) edge (m-2-1)
(m-1-2) edge node[right]{$\fib$} (m-2-2) 
(m-2-1) edge (m-2-2); 
\end{tikzpicture}
\end{equation*}
is a pullback square, where $\fib\colon\fromto{(\Phi/T)}{\Phi}$ is the functor that assigns to any morphism $\fromto{I}{T}$ its fiber over $t$.
\begin{proof} Compose the equivalence of categories $\equivto{\Phi^{\mathrm{cons}}}{(\Phi^{\mathrm{cons}}/T)}$ with the inclusion $\into{(\Phi^{\mathrm{cons}}/T)}{(\Phi/T)}$ to construct the functor $\chi$. Now the result follows from the universal property of $T$.
\end{proof}
\end{prp}

\begin{dfn}\label{dfn:specialgeneric} Suppose $\Phi$ an operator category, and suppose $(T,t)$ a point classifier for $\Phi$. Then for any object $(I,i)$ of $\Phi_{\ast}$, the unique conservative morphism $\fromto{(I,i)}{(T,t)}$ will be called the \textbf{\emph{classifying morphism}} for $i$, and will be denoted $\chi_i$. We shall call the point $t\in|T|$ the \textbf{\emph{special point}} of $T$, and for any morphism $\fromto{I}{T}$ of $\Phi$, the fiber $I_t$ will be called the \textbf{\emph{special fiber}}. Write $\fib\colon\fromto{(\Phi/T)}{\Phi}$ for the \textbf{\emph{special fiber functor}} $\goesto{I}{I_t}$.
\end{dfn}

\begin{exm}\label{exm:ptclassifiers}
\begin{enumerate}[(\ref{exm:ptclassifiers}.1)]
\item The category $\{1\}$ trivially has a point classifier.
\item The category $\mathbf{O}_{\ast}$ of pointed ordered finite sets has a point classifier, namely, $T_{\OO}\coloneq\{0,1,2\}$, wherein $1$ is the special point. Indeed, suppose $(J,j)$ a pointed ordered finite set. Then there is a morphism $\chi_j\colon\fromto{J}{T_{\OO}}$ defined by the formula
\begin{equation*}
\chi_j(k)\coloneq \begin{cases}0&\textrm{if }k<j;\\
1&\textrm{if }k=j;\\
2&\textrm{else.}
\end{cases}
\end{equation*}
The fiber of $\chi_j$ is of course the point $j$ itself, and $\chi_j$ is moreover unique with this property.
\item The category $\FF_{\ast}$ of pointed finite sets has a point classifier, namely, the set $T_{\FF}\coloneq\{0,1\}$, wherein $1$ is the special point. Indeed, this is a consequence of the observation that $T_{\FF}$ is the subobject classifier for the topos of sets.
\item The category $(\OO\wr\OO)_{\ast}$ has a point classifier. This is the object
\begin{equation*}
(T_{\OO},\{\ast,T_{\OO},\ast\}),
\end{equation*}
which may be pictured thus:
\begin{equation*}
\begin{tikzpicture}[baseline=1.25ex]
\matrix(m)[matrix of math nodes, 
row sep={2ex,between origins}, column sep={4ex,between origins}, 
text height=0ex, text depth=-1ex] 
{&&\ast&T_{\OO}&\ast&\\
&&&&&\\
&&&\bullet&&\\
&&&&&\\
&&\bullet&\circ&\bullet&\\
&&&&&\\
&&&\bullet&&\\
&{}&&&&{}\\
T_{\OO}:&&\bullet&\bullet&\bullet&\\}; 
\path[>=stealth] 
(m-8-2) edge (m-8-6); 
\end{tikzpicture}
\end{equation*}
This trend continues: the category $(\OO^{(n)})_{\ast}$ has a point classifier, which may be represented in $\RR^{n}$ as the special point at the origin and $2n$ points at the intersection of the unit $(n-1)$-sphere and the coordinate axes.
\end{enumerate}
\end{exm}

\begin{exm} Suppose $F\colon\fromto{\Psi}{\Phi}$ a fully faithful operator morphism, and suppose $(T,t)$ an object of $\Psi_{\ast}$ such that $F(T,t)$ is a point classifier for $\Phi_{\ast}$. Then $(T,t)$ is a point classifier for $\Psi_{\ast}$ as well.

It follows that for any integer $m\geq 3$, the category $\OO_{\leq m,\ast}$ has a point classifier, and for any integer $n\geq 2$, the category $\mathbf{F}_{\leq n,\ast}$ has a point classifier.
\end{exm}

\begin{dfn}\label{dfn:perfopcat} An operator category $\Phi$ is \textbf{\emph{perfect}} if the following conditions are satisfied.
\begin{enumerate}[(\ref{dfn:perfopcat}.1)]
\item The category $\Phi_{\ast}$ contains a point classifier $(T_{\Phi},t_{\Phi})$.
\item The special fiber functor $\fib\colon\fromto{(\Phi/T_{\Phi})}{\Phi}$ admits a right adjoint $E_{\Phi}$.
\end{enumerate}
One denotes the full subcategory of $\Adm$ (respectively, of $\Op$) spanned by the perfect operator categories by $\Adm^{\mathrm{perf}}$ (resp., $\Op^{\mathrm{perf}}$).
\end{dfn}

\begin{ntn} As a rule, we drop the subscripts from the notation for the point classifier and the right adjoint of $\fib$ if it is clear from the context which operator category is under consideration.

If $\Phi$ is a perfect operator category, then let us abuse notation by writing $T$ for the endofunctor of $\Phi$ obtained by composing $E\colon\fromto{\Phi}{(\Phi/T_{\Phi})}$ with the forgetful functor $\fromto{(\Phi/T_{\Phi})}{\Phi}$. Hence for some object $I$ of $\Phi$, the notation $EI$ will denote $TI$ along with the \emph{structure morphism} $e_I\colon\fromto{TI}{T}$. This abuse is partly justified by the following observation.
\end{ntn}

\begin{prp} If $\Phi$ is a perfect operator category, then the structural morphism $e_{1}\colon\fromto{T(1)}{T}$ is an isomorphism.
\end{prp}

\begin{cor} If $\Phi$ is a perfect operator category, and $I$ is an object of $\Phi$, then $e_I=T(!)$, where $!\colon\fromto{I}{1}$ is the canonical morphism.
\end{cor}

\begin{exm}\label{exm:perfopcats}
\begin{enumerate}[(\ref{exm:perfopcats}.1)]
\item Of course the operator category $\{1\}$ is perfect.
\item The operator category $\OO$ is perfect; the functor $E$ assigns to any ordered finite set $I$ the ordered finite set $TI$ obtained by adding a single point at the beginning and a single point at the end, along with the unique map $e_I\colon\fromto{TI}{T}$ whose special fiber is precisely $I\subset TI$.
\item The operator category $\FF$ is perfect; the functor $E$ assigns to any finite set $I$ the finite set $TI$ obtained by adding a disjoint basepoint to $I$, along with the unique map $e_I\colon\fromto{TI}{T}$ whose special fiber is precisely $I\subset TI$.
\item For any integer $n\geq 1$, neither $\OO_{\leq n}$ nor $\FF_{\leq n}$ are perfect.
\item The operator category $\OO\wr\OO$ is perfect. The functor $T_{\OO\wr\OO}$ carries an object $(I,M)$ to an object
\begin{equation*}
(T_{\OO}I,\{S_{\eta}\}_{\eta\in|T_{\OO}I|}),
\end{equation*}
where $S_{\eta}=1$ if $\eta$ is either of the endpoints in $T_{\OO}I$ and $S_{\eta}=T_{\OO}M_i$ if $\eta=i\in|I|\subset|T_{\OO}I|$.
\end{enumerate}
\end{exm}

More generally, wreath products of perfect operator categories are perfect.
\begin{prp} Suppose $\Phi$ and $\Psi$ two perfect operator categories. Then the operator category $\Psi\wr\Phi$ is perfect as well.
\begin{proof} Choose point classifiers $(T_{\Phi},t_{\Phi})$ of $\Phi$ and $(T_{\Psi},t_{\Psi})$ of $\Psi$ and a terminal object $1\in\Psi$. Consider an object $T_{\Psi\wr\Phi}=(T_{\Phi},\{S_{\eta}\}_{\eta\in|T_{\Phi}|})$ of $\Psi\wr\Phi$ in which
\begin{equation*}
S_{\eta}=\begin{cases}
1&\textrm{if }\eta\neq t_{\Phi};\\
T_{\Psi}&\textrm{if }\eta=t_{\Phi}.
\end{cases}
\end{equation*}
Consider the point $t_{\Psi\wr\Phi}\coloneq(t_{\Phi},t_{\Psi})\in|T_{\Psi\wr\Phi}|$. Clearly the pair $(T_{\Psi\wr\Phi},t_{\Psi\wr\Phi})$ is a point classifier of $\Psi\wr\Phi$.

The right adjoint $E_{\Psi\wr\Phi}$ of the special fiber functor $\fib$ is defined by carrying any object $(I,\{M_i\}_{i\in|I|})$ to the object $(T_{\Phi}I,\{N_j\}_{j\in|T_{\Phi}I|})$, where
\begin{equation*}
N_j=\begin{cases}
1&\textrm{if }j\notin|I|\subset|T_{\Phi}I|;\\
T_{\Psi}M_j&\textrm{if }j\in|I|\subset|T_{\Phi}I|,
\end{cases}
\end{equation*}
along with the morphism
\begin{equation*}
\fromto{(T_{\Phi}I,\{N_j\}_{j\in|T_{\Phi}I|})}{(T_{\Phi},\{S_{\eta}\}_{\eta\in|T_{\Phi}|})}
\end{equation*}
given by $e_I\colon\fromto{T_{\Phi}I}{T_{\Phi}}$ and, for any $j\in|T_{\Phi}I|$, the morphism $e_{M_{j}}\colon\fromto{T_{\Psi}M_j}{T_{\Psi}}$ when $j\in|I|$, and the identity map on $1$ when $j\notin|I|$.
\end{proof}
\end{prp}
\noindent This allows us to restrict the monoidal structure of Pr. \ref{prp:Opmonoidal} to a monoidal structure on $\Op^{\mathrm{perf}}$.
\begin{cor} Denote by $\Op^{\mathrm{perf},\wr}$ the full subcategory of $\Op^{\wr}$ spanned by those pairs $(m,X)$ such that each $X_i$ is perfect. Then the composite $\fromto{(\Op^{\mathrm{perf},\wr})^{\op}}{N\Delta^{\op}}$ is a monoidal $\infty$-category.
\end{cor}

\noindent Consequently, we obtain the following example.
\begin{exm} The operator categories $\OO^{(n)}$ are all perfect.
\end{exm}


\section{The canonical monad on a perfect operator category}\label{sect:canonicalmonad} A perfect operator category $\Phi$ always comes equipped with a monad. In effect, this monad adds points ``in every direction'' in $\Phi$; these ``directions'' are indexed by the ``non-special'' points of the point classifier.

\begin{nul} Observe that the for any object $I$ of an operator category $\Phi$, and any point $i\in|I|$, the fiber functor $(-)_i$ is already right adjoint to the fully faithful functor $p_i\colon\fromto{\Phi}{(\Phi/I)}$, which sends any object $J$ of $\Phi$ to the composite
\begin{equation*}
J\ \tikz[baseline]\draw[>=stealth,->,font=\scriptsize](0,0.5ex)--(0.5,0.5ex);\ \{i\}\ \tikz[baseline]\draw[>=stealth,right hook->,font=\scriptsize](0,0.5ex)--(0.5,0.5ex);\ I.
\end{equation*}
Observe that $\Phi$ thus has the structure of a localization of $(\Phi/I)$; that is, the unit $\fromto{J}{p_i(J)_i}$ is an isomorphism.

If $\Phi$ is perfect, there is a string of adjoints
\begin{equation*}
\Phi\ \begin{tikzpicture}[baseline] \draw[>=stealth,->,font=\scriptsize] (0,1.25ex) --node[above]{$p_t$} (0.75,1.25ex); \draw[>=stealth,->,font=\scriptsize] (0.75,0.5ex) -- (0,0.5ex); \draw[>=stealth,->,font=\scriptsize] (0,-0.25ex) --node[below]{$E$} (0.75,-0.25ex);\end{tikzpicture}\ (\Phi/T)
\end{equation*}
and as the following result shows, $\Phi$ is both a localization and a colocalization of $(\Phi/T)$.
\end{nul}

\begin{lem}\label{lem:Phicolocal} If $\Phi$ is a perfect operator $\XX$-category, then the adjoint pair $(\fib,E)$ gives $\Phi$ the structure of a colocalization of $(\Phi/T)$.
\begin{proof} The claim is simply that the counit $\kappa\colon\fromto{\fib\circ E}{\id_{\Phi}}$ is an isomorphism. The inverse to the unit $\fromto{\id_{\Phi}}{\fib\circ p_t}$ induces a morphism $\epsilon\colon\fromto{p_t}{E}$ of functors $\fromto{\Phi}{(\Phi/T)}$; this gives a morphism $\iota\colon\fromto{\id_{\Phi}}{T}$ of endofunctors of $\Phi$.

Now for any object $I$ of $\Phi$, the resulting square
\begin{equation*}
\begin{tikzpicture} 
\matrix(m)[matrix of math nodes, 
row sep=4ex, column sep=4ex, 
text height=1.5ex, text depth=0.25ex] 
{I&TI\\ 
\{t\}&T\\}; 
\path[>=stealth,->,font=\scriptsize] 
(m-1-1) edge (m-1-2) 
edge (m-2-1) 
(m-1-2) edge (m-2-2) 
(m-2-1) edge[right hook->] (m-2-2); 
\end{tikzpicture}
\end{equation*}
is a pullback square; indeed, a commutative diagram
\begin{equation*}
\begin{tikzpicture} 
\matrix(m)[matrix of math nodes, 
row sep=4ex, column sep=4ex, 
text height=1.5ex, text depth=0.25ex] 
{I'&TI\\ 
\{t\}&T\\}; 
\path[>=stealth,->,font=\scriptsize] 
(m-1-1) edge (m-1-2) 
edge (m-2-1) 
(m-1-2) edge (m-2-2) 
(m-2-1) edge[right hook->] (m-2-2); 
\end{tikzpicture}
\end{equation*}
is the same data as a morphism $\fromto{p_{t}I'}{EI}$ of $(\Phi/T)$, which is in turn the same data as a $\fromto{I'\cong(p_{t}I')_t}{I}$. 
\end{proof}
\end{lem}

If $\Phi$ is a perfect operator category, then, in effect, the endofunctor $T\colon\fromto{\Phi}{\Phi}$, when applied to an object of $\Phi$, adds as few points as possible in as many directions as possible. It turns out that this endofunctor is a monad; let us now construct a natural transformation $\mu\colon\fromto{T^{\,2}}{T}$ that, together with the natural transformation $\iota\colon\fromto{\id_{\Phi}}{T}$ from the previous proof, will exhibit a monad structure on $T$.

\begin{ntn} Suppose $\Phi$ a perfect operator category. The embedding $\iota_T\colon\into{T}{TT}$ permits us to regard the special point $t\in|T|$ as a point of $TT$ as well. Now consider the classifying morphism $\chi_t\colon\fromto{TT}{T}$. Its special fiber is the point $\iota_T(t)$, so that the following diagram is a pullback:
\begin{equation*}
\begin{tikzpicture}[baseline]
\matrix(m)[matrix of math nodes, 
row sep=4ex, column sep=5ex, 
text height=1.5ex, text depth=0.25ex] 
{\{t\}&T&TT\\
\{t\}&&T\\}; 
\path[>=stealth,->,font=\scriptsize] 
(m-1-1) edge[right hook->] (m-1-2)
(m-1-2) edge[right hook->] node[above]{$\iota_T$} (m-1-3)
(m-1-1) edge[-,double distance=1.5pt] (m-2-1)
(m-1-3) edge node[right]{$\chi_t$} (m-2-3) 
(m-2-1) edge[right hook->] (m-2-3); 
\end{tikzpicture}
\end{equation*}
\end{ntn}

\begin{ntn}\label{nul:EI} For any object $I$ of a perfect operator category $\Phi$, the functor $E$ induces a functor
\begin{equation*}
E_{/I}\colon\fromto{(\Phi/I)}{(\Phi/TI)},
\end{equation*}
which is right adjoint to the functor $\iota_I^{-1}\coloneq -\times_{TI}I\colon\fromto{(\Phi/TI)}{(\Phi/I)}$.
\end{ntn}

\begin{nul} This feature of the endofunctor $T\colon\fromto{\Phi}{\Phi}$ is what makes it a \emph{local} or \emph{parametric right adjoint}; see \cite[Df. 2.3]{MR2369114}.
\end{nul}

\begin{lem}\label{lem:rho} If $\Phi$ is a perfect operator category, then there is an isomorphism of functors
\begin{equation*}
\rho\colon\equivto{\fib}{\fib\circ\chi_{t,!}\circ E_{/TT}},
\end{equation*}
where $\chi_{t,!}\colon\fromto{(\Phi/TT)}{(\Phi/T)}$ is the functor given by composition with $\chi_t$.
\begin{proof} For any morphism $\phi\colon\fromto{J}{T}$, every square of the diagram
\begin{equation*}
\begin{tikzpicture} 
\matrix(m)[matrix of math nodes, 
row sep=6ex, column sep=6ex, 
text height=1.5ex, text depth=0.25ex] 
{J_t&J&TJ\\ 
\{t\}&T&TT\\
\{t\}&&T\\}; 
\path[>=stealth,->,font=\scriptsize] 
(m-1-1) edge[right hook->] (m-1-2) 
edge (m-2-1) 
(m-1-2) edge node[right]{$\phi$} (m-2-2) 
(m-2-1) edge[right hook->] (m-2-2)
(m-1-2) edge[right hook->] (m-1-3) 
(m-2-1) edge[-,double distance=1.5pt] (m-3-1)
(m-2-2) edge[right hook->] node[below]{$\iota_T$} (m-2-3)
(m-3-1) edge[right hook->] (m-3-3) 
(m-1-3) edge node[right]{$T\phi$} (m-2-3)
(m-2-3) edge node[right]{$\chi_t$} (m-3-3); 
\end{tikzpicture}
\end{equation*}
is a pullback square.
\end{proof}
\end{lem}

\begin{nul}\label{nul:explainsigma} By adjunction, we obtain a natural transformation $\sigma\colon\fromto{\chi_{t,!}\circ E_{/T}}{E\circ\fib}$. If $\phi\colon\fromto{I}{T}$ is a morphism of $\Phi$, then one can apply $T$ to this morphism to obtain a morphism $T(\phi)\colon\fromto{TI}{TT}$. One can, alternatively apply $T$ to the special fiber $I_t$ to obtain a morphism $\fromto{T(I_t)}{T}$. The component $\sigma_I$ then fits into a commutative square
\begin{equation*}
\begin{tikzpicture}[baseline]
\matrix(m)[matrix of math nodes, 
row sep=6ex, column sep=6ex, 
text height=1.5ex, text depth=0.25ex] 
{TI&T(I_t)\\ 
TT&T\\}; 
\path[>=stealth,->,font=\scriptsize] 
(m-1-1) edge node[above]{$\sigma_J$} (m-1-2) 
edge node[left]{$T(\phi)$} (m-2-1) 
(m-1-2) edge node[right]{$e_{I_t}$} (m-2-2) 
(m-2-1) edge node[below]{$\chi_t$} (m-2-2); 
\end{tikzpicture}
\end{equation*}
Here, the special fiber of $I$ inside $TI$ is mapped isomorphically to the special fiber of $T(I_t)$ under $\sigma_I$.
\end{nul}

\begin{dfn} Suppose $\Phi$ a perfect operator category. Define a morphism of endofunctors $\mu\colon\fromto{T^2}{T}$ as the composition
\begin{equation*}
T^2=U\circ E\circ U\circ E=U\circ \chi_{t,!}\circ E_{/T}\circ E\ \tikz[baseline]\draw[>=stealth,->,font=\scriptsize](0,0.5ex)--node[above]{$\id\circ\sigma\circ\id$}(1.5,0.5ex);\ U\circ E\circ\fib\circ E\ \tikz[baseline]\draw[>=stealth,->,font=\scriptsize](0,0.5ex)--node[above]{$\id\circ\kappa$}(0.75,0.5ex);\ U\circ E=T,
\end{equation*}
where $\kappa\colon\fromto{\fib\circ E}{\id_{\Phi}}$ is the counit isomorphism (Lm. \ref{lem:Phicolocal}).
\end{dfn}

\begin{nul}\label{nul:explainmu} More explicitly, if $I$ is an object of $\Phi$, then one has, following \ref{nul:explainsigma}, a commutative diagram
\begin{equation}\label{eqn:explainmu}
\begin{tikzpicture}[baseline]
\matrix(m)[matrix of math nodes, 
row sep=5ex, column sep=6ex, 
text height=1.5ex, text depth=0.25ex] 
{T^{\,2}I&T((TI)_t)&TI\\
TT&&T\\}; 
\path[>=stealth,->,font=\scriptsize] 
(m-1-1) edge node[above]{$\sigma_{TI}$} (m-1-2)
(m-1-2) edge node[above]{$\cong$} (m-1-3)
(m-1-1) edge node[left]{$T(e_I)$} (m-2-1)
(m-1-3) edge (m-2-3) 
(m-2-1) edge node[below]{$\chi_t$} (m-2-3); 
\end{tikzpicture}
\end{equation}
The composite $\fromto{T^{\,2}I}{TI}$ is the component $\mu_I$.
\end{nul}

\begin{thm}\label{th:Tisamonad} The endofunctor $T$ on a perfect operator category $\Phi$ is a monad with unit $\iota\colon\fromto{\id_{\Phi}}{T}$ and multiplication $\mu\colon\fromto{T^{\,2}}{T}$.
\end{thm}
\noindent The proof, though quite elementary, is a little involved, so we postpone it (\S \ref{sect:proofthatTismonad}).

\begin{exm} When $\Phi=\FF$, the monad $T$ is the \emph{partial map monad}, so that the set of maps $\fromto{J}{TI}$ is precisely the set of partial maps $\fromto{J\supseteq K}{I}$. We emphasize, however, that this is not the case in general.
\end{exm}

\begin{nul} One notes that the \emph{monad} $T$ is a \emph{local} or \emph{parametric right adjoint} \cite[Df. 2.3]{MR2369114} as well. Work of Mark Weber \cite[Pr. 2.6]{MR2369114} and \cite[Pr. 5.9]{MR2116333} shows that this implies instantly that the Kleisli category of this monad admits a factorization system, which is the inert/active factorization we describe in a more pedestrian way in Lm. \ref{lem:actinert} below. (See also \cite[\S 2.6]{MR2925893}.)
\end{nul}

Let us turn to the functoriality of this monad structure in admissible functors.
\begin{nul} An admissible functor $F\colon\fromto{\Psi}{\Phi}$ between perfect operator categories induces a functor
\begin{equation*}
F_{/T_{\Psi}}\colon\fromto{(\Psi/T_{\Psi})}{(\Phi/FT_{\Psi})},
\end{equation*}
and there is a unique conservative morphism $\chi_{F(t_{\Psi})}\colon\fromto{(FT_{\Psi},F(t_{\Psi}))}{(T_{\Phi},t_{\Phi})}$ of $\Phi_{\ast}$. So let $F_{/T}$ (with no subscript on $T$) denote the composite
\begin{equation*}
(\Psi/T_{\Psi})\ \tikz[baseline]\draw[>=stealth,->,font=\scriptsize](0,0.5ex)--node[above]{$F_{/T_{\Psi}}$}(0.75,0.5ex);\ (\Phi/FT_{\Psi})\ \tikz[baseline]\draw[>=stealth,->,font=\scriptsize](0,0.5ex)--node[above]{$\chi_{F(t_{\Psi}),!}$}(1.25,0.5ex);\ (\Phi/T_{\Phi}).
\end{equation*}
\end{nul}

We now have the following trivial observations.
\begin{lem}\label{lem:heresbeta} For any admissible functor $F\colon\fromto{\Psi}{\Phi}$ between perfect operator categories, there is a natural isomorphism $\beta_F\colon\equivto{\fib_{\Phi}\circ F_{/T}}{F\circ\fib_{\Psi}}$.
\end{lem}

\begin{lem} For any admissible functor $F\colon\fromto{\Psi}{\Phi}$ between perfect operator categories, the natural transformation $\fromto{F}{\fib_{\Phi}\circ F_{/T}\circ E_{\Psi}}$ corresponding to the isomorphism $\beta_F$ by adjunction is itself an isomorphism.
\end{lem}
\noindent Adjoint to the inverse of this isomorphism is a natural transformation
\begin{equation*}
\alpha_F\colon\fromto{F_{/T}\circ E_{\Psi}}{E_{\Phi}\circ F}.
\end{equation*}
It may be characterized as follows.
\begin{prp}\label{prp:alphaFsummary} For any object $I$ of $\Psi$, there is a unique commutative square
\begin{equation*}
\begin{tikzpicture}[baseline]
\matrix(m)[matrix of math nodes, 
row sep=4ex, column sep=4ex, 
text height=1.5ex, text depth=0.25ex] 
{FT_{\Psi}I&T_{\Phi}FI\\ 
FT_{\Psi}&T_{\Phi}\\}; 
\path[>=stealth,->,font=\scriptsize] 
(m-1-1) edge node[above]{$\alpha_{F,I}$} (m-1-2) 
edge (m-2-1) 
(m-1-2) edge (m-2-2) 
(m-2-1) edge node[below]{$\chi_{F(t_{\Psi})}$} (m-2-2); 
\end{tikzpicture}
\end{equation*}
of $\Phi$ whose special fiber is the square
\begin{equation*}
\begin{tikzpicture}[baseline]
\matrix(m)[matrix of math nodes, 
row sep=4ex, column sep=4ex, 
text height=1.5ex, text depth=0.25ex] 
{FI&FI\\ 
F\{t_{\Psi}\}&\{t_{\Phi}\}\\}; 
\path[>=stealth,->,font=\scriptsize] 
(m-1-1) edge[-,double distance=1.5pt] (m-1-2) 
edge (m-2-1) 
(m-1-2) edge (m-2-2) 
(m-2-1) edge[-,double distance=1.5pt]  (m-2-2); 
\end{tikzpicture}
\end{equation*}
\end{prp}

The natural transformation $\alpha_F$ is generally not an isomorphism, but it does behave well with respect to the monad structure.
\begin{nul} More precisely, recall that if $C$ and $D$ are categories equipped with monads $(T_C,\iota_C,\mu_C)$ and $T_D$, then a \emph{colax morphism of monads} $(F,\eta)\colon\fromto{(C,T_C)}{(D,T_D)}$ is a functor $F\colon\fromto{C}{D}$ equipped with a natural transformation $\eta\colon\fromto{FT_C}{T_DF}$ with the property that the diagrams
\begin{equation}\label{eq:colaxmorph}
\begin{tikzpicture}[baseline]
\matrix(m)[matrix of math nodes, 
row sep={5ex,between origins}, column sep={5ex,between origins}, 
text height=1.5ex, text depth=0.25ex] 
{&&FT_C^{\,2}&[3ex]\\
T_DFT_C&&&\\
&&&FT_C\\
T_D^{\,2}F&&&\\
&&T_DF&\\}; 
\path[>=stealth,->,font=\scriptsize] 
(m-1-3) edge node[above left]{$\eta T_C$} (m-2-1) 
edge node[above right]{$F\mu_C$} (m-3-4)
(m-4-1) edge node[below left]{$\mu_DF$} (m-5-3)
(m-3-4) edge node[below right]{$\eta$} (m-5-3)
(m-2-1) edge node[left]{$T_C\eta$} (m-4-1); 
\end{tikzpicture}
\textrm{\qquad and\qquad}
\begin{tikzpicture}[baseline]
\matrix(m)[matrix of math nodes, 
row sep=3ex, column sep=4ex, 
text height=1.5ex, text depth=0.25ex] 
{&FT_C\\
F&\\
&T_DF.\\}; 
\path[>=stealth,->,font=\scriptsize] 
(m-2-1) edge node[above left]{$F\iota$} (m-1-2)
edge node[below left]{$\iota F$} (m-3-2)
(m-1-2) edge node[right]{$\eta$} (m-3-2);
\end{tikzpicture}
\end{equation}
of $\Fun(C,D)$ commute.
\end{nul}
\noindent Thus the interaction of the natural transformation $\alpha_F$ with the monad structures on $\Psi$ and $\Phi$ is summarized by the following result.
\begin{thm}\label{thm:admissiblesinducecolax} An admissible functor $F\colon\fromto{\Psi}{\Phi}$ of perfect operator categories induces a colax morphism of monads $(F,\alpha_F)\colon\fromto{(\Psi,T_{\Psi})}{(\Phi,T_{\Phi})}$.
\end{thm}
\noindent Again we postpone the proof (\S \ref{sect:proofthatadmissiscolax}).

\begin{nul} Observe that the uniqueness of conservative morphisms with target $(T_{\Phi},t_{\Phi})$ implies that if $G\colon\fromto{X}{\Psi}$ is another admissible functor of perfect operator $\XX$-categories, then $(F\circ G)_{/T}=F_{/T}\circ G_{/T}$. Hence the assignment $\goesto{\Phi}{(\Phi/T)}$ defines a functor $\fromto{\Adm^{\mathrm{perf}}}{\Cat}$. On the other hand, there is a forgetful functor $\fromto{\Adm^{\mathrm{perf}}}{\Cat}$.

Now the functor $\beta$ of Lm. \ref{lem:heresbeta} can be regarded as a natural transformation from the functor $\goesto{\Phi}{(\Phi/T)}$ to the forgetful functor $\fromto{\Adm^{\mathrm{perf}}}{\Cat}$ such that for any perfect operator category $\Phi$, one has $\beta_{\Phi}=\fib_{\Phi}$. Simlarly, $\alpha$ can be regarded as a lax natural transformation from the forgetful functor $\fromto{\Adm^{\mathrm{perf}}}{\Cat}$ to the functor $\goesto{\Phi}{(\Phi/T)}$ such that for any perfect operator category $\Phi$, one has $\alpha_{\Phi}=E_{\Phi}$.
\end{nul}

\begin{nul} A $2$-morphism $\xi\colon\fromto{(F,\eta)}{(F',\eta')}$ of colax morphisms of monads is an isomorphism of functors $\xi\colon\fromto{F}{F'}$ such that the square
\begin{equation*}
\begin{tikzpicture}[baseline]
\matrix(m)[matrix of math nodes, 
row sep=6ex, column sep=6ex, 
text height=1.5ex, text depth=0.25ex] 
{FT_C&T_DF\\ 
F'T_C&T_DF'\\}; 
\path[>=stealth,->,font=\scriptsize] 
(m-1-1) edge node[above]{$\eta$} (m-1-2)
(m-1-1) edge node[left]{$\xi T_C$} (m-2-1)
(m-1-2) edge node[right]{$T_D\xi$} (m-2-2) 
(m-2-1) edge node[below]{$\eta$} (m-2-2); 
\end{tikzpicture}
\end{equation*}
commutes.

Composition and identities are defined in the obvious manner; hence this defines a $2$-category $\mathrm{Mnd}_{\mathrm{colax}}$ of small categories with monads and colax functors. By applying the nerve of each $\Mor$-groupoid and taking the simplicial nerve of the resulting simplicial category, we obtain a quasicategory $\Mnd_{\mathrm{colax}}$, which is in fact a $2$-category in the sense of \cite[\S 2.3.4]{HTT}.

We may summarize Th. \ref{th:Tisamonad} and Th. \ref{thm:admissiblesinducecolax} together by stating that the assignment $\goesto{\Phi}{(\Phi,T)}$ defines a functor $\fromto{\Adm^{\mathrm{perf}}}{\mathbf{Mnd}_{\mathrm{colax}}}$.
\end{nul}


\section{Leinster categories}\label{sect:Leinster} The Kleisli category of the canonical monad on a perfect operator category $\Phi$ is the category of free algebras for this monad. It can be thought of as indexing both operations in $\Phi$ as well as projection maps in a coherent manner. In the examples of interest, this Kleisli category recovers a number of combinatorial categories familiar to homotopy theorists.

In particular, we'll find that the Kleisli category of the canonical monad on $\FF$ is Segal's $\Gamma^{\op}$, the Kleisli category of the canonical monad on $\OO$ is $\Delta^{\op}$, and the Kleisli category of the canonical monad on $\OO^{(n)}$ is Joyal's $\Theta_n^{\op}$

\begin{dfn} Suppose $\Phi$ a perfect operator category. Then the \textbf{\emph{Leinster category}} $\Lambda(\Phi)$ of $\Phi$ is the Kleisli category of the monad $T_{\Phi}$.
\end{dfn}

\begin{nul} That is, the objects of the Leinster category of a perfect operator category $\Phi$ are precisely those of $\Phi$ itself, and for any two objects $I$ and $J$,
\begin{equation*}
\Mor_{\Lambda(\Phi)}(J,I)\coloneq \Mor_{\Phi}(J,TI).
\end{equation*}
The identity at an object $I$ is the unit $\iota_I\colon\fromto{I}{TI}$. The composition law is defined by the composite
\begin{equation*}
\begin{tikzpicture}[baseline]
\matrix(m)[matrix of math nodes, 
row sep=2ex, column sep=6ex, 
text height=1.5ex, text depth=0.25ex] 
{\Mor_{\Phi}(K,TJ)\times\Mor_{\Phi}(J,TI)&\Mor_{\Phi}(K,TJ)\times\Mor_{\Phi}(TJ,T^{\,2}I)\\
{\qquad\qquad\qquad\qquad\qquad}&\Mor_{\Phi}(K,T^{\,2}I)\\
{\qquad\qquad\qquad\qquad\qquad\qquad\qquad\qquad}&\Mor_{\Phi}(K,TI)\\ }; 
\path[>=stealth,->,font=\scriptsize]
(m-1-1) edge node[above]{$\id\times T_{\Phi}$} (m-1-2)
(m-2-1) edge node[above]{$\circ$} (m-2-2)
(m-3-1) edge node[above]{$\mu_{I,\ast}$} (m-3-2);
\end{tikzpicture}
\end{equation*}
for objects $I$, $J$, and $K$.
\end{nul}

\begin{nul} Suppose now $F\colon\fromto{\Psi}{\Phi}$ an admissible functor of perfect operator categories. Since colax morphisms of monoids induce functors of their Kleisli categories --- and in fact this is a functor $\fromto{\Mnd_{\mathrm{colax}}}{\Cat}$ ---, it follows that $f$ induces a functor $\Lambda(F)\colon\fromto{\Lambda(\Psi)}{\Lambda(\Phi)}$ of the Leinster categories.

The admissible functor $F$ induces a functor $\Lambda(F)\colon\fromto{\Lambda(\Psi)}{\Lambda(\Phi)}$ such that the map on objects is simply
\begin{equation*}
\Obj\Lambda(F)=\Obj F\colon\fromto{\Obj\Lambda(\Psi)=\Obj\Psi}{\Obj\Phi=\Obj\Lambda(\Phi)},
\end{equation*}
and the map on morphisms is given by
\begin{equation*}
\Mor_{\Psi}(J,T_{\Psi}I)\ \tikz[baseline]\draw[>=stealth,->,font=\scriptsize](0,0.5ex)--(0.5,0.5ex);\ \Mor_{\Phi}(FJ,FT_{\Psi}I)\ \tikz[baseline]\draw[>=stealth,->,font=\scriptsize](0,0.5ex)--node[above]{$\alpha_{F,\ast}$}(0.75,0.5ex);\ \Mor_{\Phi}(FJ,T_{\Phi}FI)
\end{equation*}
(Pr. \ref{prp:alphaFsummary}).

Consequently, the Leinster category construction defines a functor
\begin{equation*}
\Lambda\colon\fromto{\Adm^{\mathrm{perf}}}{\Cat}.
\end{equation*}
\end{nul}

\begin{exm} For any object $I$ of a perfect operator category $\Phi$ and for any point $i\in|I|$, the classifying morphism $\chi_i\colon\fromto{I}{T}$ of $\Phi$ is a morphism $\fromto{I}{\{i\}}$ of $\Lambda(\Phi)$. In particular, of course, $\{i\}$ is not terminal in $\Lambda(\Phi)$.
\end{exm}

\begin{exm}\label{exm:LambdaFisGamma} The Leinster category of $\mathbf{F}$ is Segal's category $\Gamma^{\op}$ of pointed finite sets \cite{MR51:4232,MR50:5782}.
\end{exm}

\begin{exm}\label{exm:LambdaOisDelta} Let us study the Leinster category of $\mathbf{O}$. Denote by $\bot$ and $\top$ the two points of $T(\varnothing)$, so that
\begin{equation*}
T(\varnothing)=\{\bot,\top\};
\end{equation*}
these correspond to two wide morphisms $\fromto{T}{T(\varnothing)}$, and for any object $I$ of $\Lambda(\OO)$, they give rise to a map
\begin{equation*}
c_I:=(\bot_!,\top_!)\colon\fromto{\Mor_{\Lambda(\OO)}(I,\ast)}{\Mor_{\Lambda(\OO)}(I,\varnothing)\times\Mor_{\Lambda(\OO)}(I,\varnothing)}.
\end{equation*}
The map $c_I$ is injective, so for any $\phi,\psi\in\Mor_{\Lambda(\OO)}(I,\varnothing)$, denote by $\phi\star\psi$ the unique element of $\Mor_{\Lambda(\OO)}(I,\ast)$ such that $c_I(\phi\star\psi)=(\phi,\psi)$, if it exists.

Then the subset $c_I\left(\Mor_{\Lambda(\OO)}(I,\ast)\right)\subset\Mor_{\Lambda(\OO)}(I,\varnothing)\times\Mor_{\Lambda(\OO)}(I,\varnothing)$ is a total ordering $\leq$ on $\Mor_{\Lambda(\OO)}(I,\varnothing)$; in particular, for any $\phi,\psi\in\Mor_{\Lambda(\OO)}(I,\varnothing)$, the element $\psi$ is the successor of $\phi$ if and only if $\phi\star\psi=\chi_i$ for some point $i\in|I|$.

The functor $\Mor_{\Lambda(\OO)}(-,\varnothing)$ thus defines a functor $\fromto{\Lambda(\OO)^{\op}}{\Delta}$, where $\Delta$ is defined as the full subcategory of $\OO$ consisting of nonempty objects; a quasiinverse functor is given by
\begin{equation*}
\goesto{\mathbf{n}}{\mathbf{n}^{\vee}\coloneq\Mor_{\Delta}(\mathbf{n},\mathbf{1})}.
\end{equation*}
We deduce that $\Lambda(\OO)$ is equivalent to $\Delta^{\op}$. (This observation goes back at least to Street \cite{MR574662}.)
\end{exm}

The following result is trivial to prove.
\begin{prp} Suppose $p\colon\fromto{\Phi'}{\Phi}$ an admissible functor between perfect operator categories that is also a Grothendieck fibration that classifies a functor $\fromto{\Phi^{\op}}{\Adm^{\mathrm{perf}}}$. Then the functor $\fromto{\Lambda(\Phi')}{\Lambda(\Phi)}$ is a Grothendieck fibration such that for any object $I\in\Phi$, one has $(\Lambda(\Phi'))_{I}\simeq \Lambda(\Phi'_I)$.
\end{prp}

\begin{exm}\label{exm:LambdaOnisThetan} The previous result now immediately implies that the Leinster category of the iterated wreath product $\OO^{(n)}$ coincides with Berger's iterated wreath product $(\Delta\wr\Delta\wr\cdots\wr\Delta)^{\op}$ \cite{MR2331244}. In particular, $\Lambda(\OO^{(n)})$ is equivalent to Joyal's category $\Theta_{n}^{\op}$.
\end{exm}

\section{Quasioperads and their algebras}\label{sect:Luriestyle} When $\Phi$ is a perfect operator category, the theory of complete Segal $\Phi$-operads introduced above can also be codified in a manner similar to Jacob Lurie's theory of $\infty$-operads \cite[Ch. 2]{HA}. Here, we explain how to generalize the basic elements of Lurie's theory to any perfect operator category; in most cases the proofs are trivial extensions of the proofs of loc. cit.

The following terminology was first introduced on $\Gamma^{\op}$ by Jacob Lurie.
\begin{dfn} Suppose $\Phi$ a perfect operator category. We call a morphism $\fromto{J}{I}$ of $\Lambda(\Phi)$ \textbf{\emph{inert}} if, the corresponding morphism $\fromto{J}{TI}$ in $\Phi$ has the property that the natural morphism
\begin{equation*}
\fromto{J\times_{TI}I}{I}
\end{equation*}
is an isomorphism. Denote by $\Lambda^{\dag}(\Phi)$ the collection of inert morphisms of $\Lambda(\Phi)$. We will simply write $N\Lambda(\Phi)$ for the marked simplicial set $(N\Lambda(\Phi),\Lambda^{\dag}(\Phi))$.

Let us call a morphism $\fromto{J}{I}$ of $\Lambda(\Phi)$ \textbf{\emph{active}} if the corresponding morphism $\fromto{J}{TI}$ of $\Phi$ factors as
\begin{equation*}
\into{J}{TJ\to TI}, 
\end{equation*}
where the morphism $\fromto{TJ}{TI}$ is of the form $T\phi$ for some morphism $\phi\colon\fromto{J}{I}$ of $\Phi$. Denote by $\Lambda_{\dag}(\Phi)$ the collection of active morphisms of $\Lambda(\Phi)$.
\end{dfn}

\begin{nul}\label{nul:easyactinert} Suppose $\Phi$ a perfect operator category. We observe the following.
\begin{enumerate}[(\ref{nul:easyactinert}.1)]
\item A morphism $J\to I$ of $\Lambda(\Phi)$ is active just in case the corresponding morphism $\fromto{J}{TI}$ in $\Phi$ has the property that the natural morphism
\begin{equation*}
\fromto{J\times_{TI}I}{J}
\end{equation*}
is an isomorphism.
\item A morphism of $\Lambda(\Phi)$ is both inert and active just in case it is an isomorphism.
\suspend{enumerate}
Suppose $\psi\colon\fromto{K}{J}$ and $\phi\colon\fromto{J}{I}$ morphisms of $\Lambda(\Phi)$.
\resume{enumerate}[{[(\ref{nul:easyactinert}.1)]}]
\item If $\psi$ is inert, then $\phi\circ\psi$ is inert just in case $\phi$ is.
\item Dually, if $\phi$ is active, then $\phi\circ\psi$ is active just in case $\psi$ is.
\end{enumerate}
\end{nul}

It is an observation of Christopher Schommer-Pries that $\Lambda(\Phi)$ always admits an inert-active factorization system. In fact, Mark Weber \cite[Pr. 2.6]{MR2369114} and \cite[Pr. 5.9]{MR2116333} proves this as a natural consequence of the fact that $T$ is a \emph{local} or \emph{parametric right adjoint}. (See also \cite[\S 2.6]{MR2925893}. In their language, what we call ``inert,'' Weber and others call ``generic,'' and what we call ``active,'' they call ``free.'') Clemens Berger has introduced the concept of a \emph{moment category} that codifies the salient features of this factorization system, which we intend to explore elsewhere. For now, we simply record the observation of Schommer-Pries (with a perhaps suboptimal proof).
\begin{lem}[C. Schommer-Pries]\label{lem:actinert} Every morphism $\fromto{J}{I}$ of the Leinster category of a perfect operator category $\Phi$ admits factorization $J\to K\to I$ into an inert morphism $\fromto{J}{K}$ followed by an active morphism $\fromto{K}{I}$. Moreover, this factorization is unique up to unique isomorphism
\begin{proof} Using the structural morphism $e_I\colon\fromto{TI}{I}$, regard the morphism $\fromto{J}{TI}$ as a morphism over $T$. Set $K\coloneq J\times_{TI}I$; the projection $\fromto{K}{I}$ in $\Phi$ induces an active morphism $\fromto{K}{I}$ in $\Lambda(\Phi)$. Now the universal property of $T$ states that the set of maps $\fromto{J}{TK}$ over $T$ is in bijection with the set of maps $\fromto{J_t}{K}$. Hence we may choose the canonical isomorphism $J_t\cong K$, yielding a morphism $\fromto{J}{TK}$. It is obvious from the construction that this now gives the desired factorization $J\to K\to I$ of the original morphism, and the uniqueness follows from the observations above.
\end{proof}
\end{lem}

\begin{exm} When $\Phi=\FF$, a pointed map $\phi\colon\fromto{J_{+}}{I_{+}}$ of $\Gamma^{\op}$ corresponds to an inert morphism in our sense just in case it is inert in the sense of Lurie, that is, just in case the inverse image $\phi^{-1}(i)$ of any point $i\in I$ is a singleton. It corresponds to an active morphism in our sense just in case it is active in the sense of Lurie, that is, just in case the inverse image of the base point is a singleton.
\end{exm}

\begin{exm} When $\Phi=\OO$, a morphism $\phi\colon\fromto{\mathbf{m}}{\mathbf{n}}$ of $\Delta^{\op}$ corresponds to an inert morphism in our sense just in case it corresponds to an injection $\fromto{\mathbf{n}}{\mathbf{m}}$ given by the formula $\goesto{i}{i+k}$ for some fixed integer $k\geq 1$. It corresponds to an active morphism in our sense just in case it corresponds to a map $\fromto{\mathbf{n}}{\mathbf{m}}$ that carries $0$ to $0$ and $n$ to $m$.
\end{exm}

\begin{prp}\label{prp:admspreserveinerts} The functor $\fromto{\Lambda(\Psi)}{\Lambda(\Phi)}$ induced by an admissible functor $F\colon\fromto{\Psi}{\Phi}$ of perfect operator categories preserves inert morphisms.
\begin{proof} Suppose $\fromto{J}{I}$ an inert morphism of $\Lambda(\Phi)$. Consider the rectangle
\begin{equation*}
\begin{tikzpicture} 
\matrix(m)[matrix of math nodes, 
row sep=4ex, column sep=4ex, 
text height=1.5ex, text depth=0.25ex] 
{F(J\times_{TI}I)&FI&FI\\ 
FJ&FTI&TFI.\\}; 
\path[>=stealth,->,font=\scriptsize] 
(m-1-1) edge node[above]{$\cong$} (m-1-2) 
edge[right hook->] (m-2-1) 
(m-1-2) edge[-,double distance=1.5pt] (m-1-3)
edge[right hook->] (m-2-2)
(m-1-3) edge[right hook->] (m-2-3)
(m-2-1) edge (m-2-2)
(m-2-2) edge node[below]{$\alpha_{F,I}$} (m-2-3); 
\end{tikzpicture}
\end{equation*}
The left hand square is a pullback since $F$ is admissible; the morphism
\begin{equation*}
\fromto{F(J\times_{TI}I)}{FI}
\end{equation*}
is an isomorphism because $\fromto{J}{I}$ is inert; and the fact that right hand square is a pullback follows from the characterization of $\alpha_{F,I}$ given in Pr. \ref{prp:alphaFsummary}.
\end{proof}
\end{prp}

\begin{ntn} Suppose $X$ is an $\infty$-category, suppose $S$ a $1$-category, and suppose $q\colon\fromto{X}{S}$ a functor. Suppose $x,y\in X$, and suppose $g\colon\fromto{q(x)}{q(y)}$ a morphism of $S$. Denote by $\Map_X^{g}(x,y)$ the union of the connected components of $\Map_X(x,y)$ lying over the connected component of $g$ in $\Map_S(q(x),q(y))$.
\end{ntn}

We can now define the notion of $\Phi$-quasioperad in exact analogy with Lurie \cite[Df. 2.1.1.10]{HA}.

\begin{dfn}\label{dfn:inftyPhioperad} Suppose $\Phi$ an operator category. Then a \textbf{\emph{$\Phi$-quasioperad}} or \textbf{\emph{$\infty$-operad over $\Phi$}} is an inner fibration $p\colon\fromto{X^{\otimes}}{N\Lambda(\Phi)}$ satisfying the following conditions.
\begin{enumerate}[(\ref{dfn:inftyPhioperad}.1)]
\item For every morphism $\phi\colon\fromto{J}{I}$ of $\Lambda^{\dag}(\Phi)$ and every object $x\in X^{\otimes}_J$, there is a $p$-cocartesian edge $\fromto{x}{y}$ in $X^{\otimes}$ covering $\phi$.
\item For any objects $I,J\in\Phi$, any objects $x\in X^{\otimes}_{I}$ and $y\in X^{\otimes}_{J}$, any morphism $\phi\colon\fromto{J}{I}$ of $\Lambda(\Phi)$, and any $p$-cocartesian edges $\{\fromto{y}{y_i}\ |\ i\in|I|\}$ lying over the inert morphisms $\{\rho_i\colon\fromto{I}{\{i\}}\ |\ i\in|I|\}$, the induced map
\begin{equation*}
\fromto{\Map_{X^{\otimes}}^{\phi}(x,y)}{\prod_{i\in|I|}\Map_{X^{\otimes}}^{\rho_i\circ\phi}(x,y_i)}
\end{equation*}
is an equivalence.
\item For any object $I\in\Phi$, the $p$-cocartesian morphisms lying over the inert morphisms $\{\fromto{I}{\{i\}}\ |\ i\in|I|\}$ together induce an equivalence
\begin{equation*}
\fromto{X^{\otimes}_I}{\prod_{i\in|I|}X^{\otimes}_{\{i\}}}.
\end{equation*}
\end{enumerate}
\end{dfn}

\begin{exm} When $\Phi=\FF$, the above definition coincides with Lurie's definition of $\infty$-operad \cite{HA}. 
\end{exm}

\begin{exm} When $\Phi=\{1\}$, the conditions of the above definition are trivial, and we are simply left with the notion of a quasicategory.
\end{exm}

In light of \cite[Pr. 2.1.2.5]{HA}, the inert-active factorization on the Leinster category of a perfect operator category lifts to any $\infty$-operad over it.
\begin{dfn} Suppose $\Phi$ a perfect operator category, and suppose
\begin{equation*}
p\colon\fromto{X^{\otimes}}{N\Lambda(\Phi)}
\end{equation*}
a $\Phi$-quasioperad. Call a $p$-cocartesian edge of $X^{\otimes}$ that covers an inert morphism of $\Lambda(\Phi)$ \textbf{\emph{inert}}. Dually, call any edge of $X^{\otimes}$ that covers an active morphism of $\Lambda(\Phi)$ \textbf{\emph{active}}.
\end{dfn}

\begin{prp} For any perfect operator category and any $\Phi$-quasioperad $X^{\otimes}$, the inert morphisms and the active morphisms determine a factorization system on $X^{\otimes}$.
\end{prp}

As in \cite[\S 2.1.4]{HA}, we can introduce a model category of $\infty$-preoperads over a perfect operator category $\Phi$ whose fibrant objects are $\infty$-operads over $\Phi$.

\begin{ntn} As in \cite[Nt. 2.1.4.5]{HA}, for any operator category $\Phi$ and any $\infty$-operad $X^{\otimes}$ over $\Phi$, denote by $X^{\otimes,\natural}$ the object $(X^{\otimes},E)$ of $s\Set^{+}_{/N\Lambda(\Phi)}$, where $E$ denotes the collection of inert morphisms of $X^{\otimes}$.
\end{ntn}

\begin{nul}\label{nul:catpattern} For any perfect operator category $\Phi$, consider the following \emph{categorical pattern} (in the sense of Lurie; see \cite[App. B]{HA})
\begin{equation*}
\mathscr{P}=(M,T,\{p_{\alpha}\colon\fromto{\Lambda^2_0}{N\Lambda(\Phi)}\}_{\alpha\in A})
\end{equation*}
on the simplicial set $N\Lambda(\Phi)$. The class $M$ consists of all the inert morphisms of $\Lambda(\Phi)$; the class $T$ is the class of all $2$-simplices; and the set $A$ is the set of diagrams $I\ot J\to I'$ of $\Lambda(\Phi)$ in which both $J\to I$ and $J\to I'$ are inert, and
\begin{equation*}
|J|=|J\times_{TI}I|\sqcup|J\times_{TI'}I'|\cong|I|\sqcup|I'|.
\end{equation*}
\end{nul}

Now applying \cite[Th. B.0.19]{HA} to the categorical pattern $\mathscr{P}$, one deduces the following.
\begin{thm}\label{thm:modstructpreoperad} Suppose $\Phi$ a perfect operator category. There exists a left proper, tractable, simplicial model structure --- called the \textbf{\emph{operadic model structure}} --- on $s\Set^{+}_{/N\Lambda(\Phi)}$ with the following properties.
\begin{enumerate}[(\ref{thm:modstructpreoperad}.1)]
\item A marked map $\fromto{X}{Y}$ over $N\Lambda(\Phi)$ is a cofibration just in case it is a monomorphism.
\item A marked simplicial set over $N\Lambda(\Phi)$ is fibrant just in case it is of the form $Z^{\otimes,\natural}$ for some $\Phi$-quasioperad $Z^{\otimes}$.
\item A marked map $\fromto{X}{Y}$ over $N\Lambda(\Phi)$ is a weak equivalence just in case, for any $\Phi$-quasioperad $Z^{\otimes}$, the induced map
\begin{equation*}
\fromto{\Map_{/N\Lambda(\Phi)}(Y,Z^{\otimes,\natural})}{\Map_{/N\Lambda(\Phi)}(X,Z^{\otimes,\natural})}
\end{equation*}
is a weak equivalence.
\end{enumerate}
\end{thm}

\begin{ntn} For any perfect operator category $\Phi$, denote by $\Operad_{\infty}^{\Phi,\Delta}$ the simplicial subcategory of $s\Set^{+}_{/N\Lambda(\Phi)}$ spanned by the fibrant objects for the operadic model structure. Since this is a fibrant simplicial category, we may apply the nerve to obtain an $\infty$-category $\Operad_{\infty}^{\Phi}$.
\end{ntn}

Using Lurie's characterization \cite[Lm. B.2.4(3)]{HA} of $\mathscr{P}$-equivalences between fibrant objects, we obtain the following.
\begin{prp}\label{prp:DKequivsonquasiops} Suppose $\Phi$ a perfect operator category. Then a morphism of $\Phi$-quasioperads $g\colon\fromto{X}{Y}$ is an equivalence just in case the following conditions are satisfied.
\begin{enumerate}[(\ref{prp:DKequivsonquasiops}.1)]
\item \emph{Essential surjectivity}. The functor $\fromto{X_{\{1\}}}{Y_{\{1\}}}$ is essentially surjective.
\item \emph{Full faithfulness}. For any object $I\in\Phi$, any vertex $x\in X_{I}$, and any vertex $y\in X_{\{1\}}$, the induced map
\begin{equation*}
\fromto{\Map^{\alpha}_{X}(x,y)}{\Map^{\alpha}_{Y}(g(x),g(y))}
\end{equation*}
is an equivalence, where $\alpha$ is the unique active morphism $\fromto{I}{\{1\}}$ of $\Lambda(\Phi)$.
\end{enumerate}
\end{prp}

We may apply \cite[Pr. B.2.9]{HA} to the functors $\fromto{\Lambda(\Psi)}{\Lambda(\Phi)}$ induced by operator morphisms $\fromto{\Psi}{\Phi}$ thanks to Pr. \ref{prp:admspreserveinerts}; it follows that that the operadic model structure on $s\Set^{+}_{/N\Lambda(\Phi)}$ enjoys the same functoriality in $\Phi$ that is enjoyed by the operadic model structure on the category $s\Set_{/N\Delta_{\Phi}^{\op}}$ (Pr. \ref{prp:operadicfunctorial}):
\begin{prp} For any operator morphism $G\colon\fromto{\Psi}{\Phi}$, the adjunction
\begin{equation*}
\adjunct{G_!}{s\Set^+_{/N\Lambda(\Psi)}}{s\Set^+_{/N\Lambda(\Phi)}}{G^{\star}}
\end{equation*}
is a Quillen adjunction for the operadic model structure.
\end{prp}


\section{Boardman--Vogt tensor products and weak algebras}\label{sect:BV} It is well-known that the classical Boardman--Vogt tensor product \cite{MR524181} exhibits better homotopical properties when it is extended to suitably weak operads, as in \cite{MR2366165} and \cite{HA}. It can also be \emph{externalized} over the wreath product construction given the previous section, as we now demonstrate.

\begin{ntn} For any operator categories $\Phi$ and $\Psi$, consider the functor
\begin{equation*}
\wr\colon\fromto{\Psi\times\Phi}{\Psi\wr\Phi}
\end{equation*}
that carries the pair $(J,I)$ to the object $J\wr I\coloneq((J_i)_{i\in|I|};I)$, in which $J_i=J$ for each $i\in|I|$. This induces a functor
\begin{equation*}
W\colon\fromto{\Delta_{\Psi}\times_{\Delta}\Delta_{\Phi}}{\Delta_{\Psi\wr\Phi}}
\end{equation*}
given by the assignment
\begin{equation*}
\goesto{([J_0\to\cdots\to J_m],[I_0\to\cdots\to I_m])}{[J_0\wr I_0\to\cdots\to J_m\wr I_m]}.
\end{equation*}
\end{ntn}

\begin{nul}\label{nul:Wandopmors} It is clear that for any operator morphisms $H\colon\fromto{\Psi'}{\Psi}$ and $G\colon\fromto{\Phi'}{\Phi}$, the square
\begin{equation*}
\begin{tikzpicture} 
\matrix(m)[matrix of math nodes, 
row sep=4ex, column sep=4ex, 
text height=1.5ex, text depth=0.25ex] 
{\Delta_{\Psi'}\times_{\Delta}\Delta_{\Phi'}&\Delta_{\Psi'\wr\Phi'}\\ 
\Delta_{\Psi}\times_{\Delta}\Delta_{\Phi}&\Delta_{\Psi\wr\Phi}\\}; 
\path[>=stealth,->,font=\scriptsize] 
(m-1-1) edge node[above]{$W$} (m-1-2) 
edge (m-2-1) 
(m-1-2) edge (m-2-2) 
(m-2-1) edge node[below]{$W$} (m-2-2); 
\end{tikzpicture}
\end{equation*}
commutes.
\end{nul}

\begin{dfn} Suppose $\Phi$ and $\Psi$ operator categories, and suppose $p\colon\fromto{X}{N\Delta^{\op}_{\Phi}}$ and $q\colon\fromto{Y}{N\Delta^{\op}_{\Psi}}$ left fibrations. For any complete Segal $\Psi\wr\Phi$-operad $Z$, a \textbf{\emph{pairing}} $\fromto{(X,Y)}{Z}$ is a commutative diagram
\begin{equation*}
\begin{tikzpicture} 
\matrix(m)[matrix of math nodes, 
row sep=4ex, column sep=4ex, 
text height=1.5ex, text depth=0.25ex] 
{Y\times_{N\Delta^{\op}}X&Z\\ 
N\Delta_{\Psi}^{\op}\times_{N\Delta^{\op}}N\Delta_{\Phi}^{\op}&N\Delta_{\Psi\wr\Phi}^{\op}\\}; 
\path[>=stealth,->,font=\scriptsize] 
(m-1-1) edge (m-1-2) 
edge node[left]{$(q,p)$} (m-2-1) 
(m-1-2) edge (m-2-2) 
(m-2-1) edge node[below]{$W$} (m-2-2); 
\end{tikzpicture}
\end{equation*}
Write $\Pair^{\Psi,\Phi}(Y,X;Z)$ for the simplicial set $\Map_{/N\Delta_{\Psi\wr\Phi}^{\op}}(Y\times_{N\Delta^{\op}}X,Z)$ of pairings $\fromto{(Y,X)}{Z}$.
\end{dfn}

\begin{prp} Suppose $\Phi$ and $\Psi$ operator categories, and suppose $\fromto{X}{N\Delta^{\op}_{\Phi}}$ and $\fromto{Y}{N\Delta^{\op}_{\Psi}}$ left fibrations. Then the functor
\begin{equation*}
\Pair^{\Psi,\Phi}(Y,X;-)\colon\fromto{\Operad^{\Psi\wr\Phi,\op}}{\Kan}
\end{equation*}
is corepresentable.
\begin{proof} Denote by $\mathbf{LFib}(S)$ the simplicial nerve of the simplicial category of left fibrations to a fixed simplicial set $S$. Denote by $L\colon\fromto{\mathbf{LFib}(N\Delta^{\op}_{\Psi\wr\Phi})}{\Operad^{\Psi\wr\Phi}}$ the left adjoint to the inclusion. Then $\Pair^{\Psi,\Phi}(Y,X;-)$ is corepresented by the object $LW_{!}(Y\times_{N\Delta^{\op}}X)$.
\end{proof}
\end{prp}

\begin{dfn} Suppose $\Phi$ and $\Psi$ operator categories, suppose $\fromto{X}{N\Delta^{\op}_{\Phi}}$ and $\fromto{Y}{N\Delta^{\op}_{\Psi}}$ left fibrations, and suppose $Z$ a complete Segal $\Psi\wr\Phi$-operad. A pairing $\fromto{(Y,X)}{Z}$ will be said to \textbf{\emph{exhibit $Z$ as the Boardman--Vogt tensor product of $Y$ and $X$}} just in case, for every operad complete Segal $\Psi\wr\Phi$-operad $Z'$ it induces an equivalence
\begin{equation*}
\equivto{\Map_{\Operad^{\Psi\wr\Phi}}(Z,Z')}{\Pair^{\Psi,\Phi}((Y,X),Z')}.
\end{equation*}
In the presence of such a pairing, we write $Y{{}^{\Psi}\otimes^{\Phi}}X$ for $Z$. Using \cite{HTT}, one can organize this into a functor
\begin{equation*}
-{{}^{\Psi}\otimes^{\Phi}}-\colon\fromto{\mathbf{LFib}(N\Delta_{\Phi}^{\op})^{\op}\times\mathbf{LFib}(N\Delta_{\Psi}^{\op})^{\op}}{\Operad^{\Psi\wr\Phi}}.
\end{equation*}
\end{dfn}

\begin{wrn} In contrast with the Boardman--Vogt tensor product of ordinary symmetric operads, note that the order matters here: since the wreath product is noncommutative, the objects $Y{{}^{\Psi}\otimes^{\Phi}}X$ and $X{{}^{\Phi}\otimes^{\Psi}}Y$ are not even objects of the same $\infty$-category.
\end{wrn}

It follows immediately from \ref{nul:Wandopmors} that the Boardman--Vogt tensor product is compatible with operator morphisms.
\begin{prp} For any operator morphisms $H\colon\fromto{\Psi'}{\Psi}$ and $G\colon\fromto{\Phi'}{\Phi}$ and for any left fibrations $\fromto{X}{N\Delta^{\op}_{\Phi'}}$ and $\fromto{Y}{N\Delta^{\op}_{\Psi'}}$, there is a canonical equivalence
\begin{equation*}
(H\wr G)_!(Y{{}^{\Psi'}\otimes^{\Phi'}}X)\simeq H_!Y{{}^{\Psi}\otimes^{\Phi}}G_!X
\end{equation*}
of complete Segal $\Psi\wr\Phi$-operads.
\end{prp}

It is easy to see that the Boardman--Vogt tensor product preserves colimits separately in each variable. Consequently, we obtain the following.
\begin{prp} Suppose $\Phi$ and $\Psi$ operator categories, and suppose $\fromto{X}{N\Delta^{\op}_{\Phi}}$ and $\fromto{Y}{N\Delta^{\op}_{\Psi}}$ left fibrations. Then the functors
\begin{equation*}
-{{}^{\Psi}\otimes^{\Phi}}X\colon\fromto{\mathbf{LFib}(N\Delta_{\Psi}^{\op})^{\op}}{\Operad^{\Psi\wr\Phi}}
\end{equation*}
and
\begin{equation*}
Y{{}^{\Psi}\otimes^{\Phi}}-\colon\fromto{\mathbf{LFib}(N\Delta_{\Phi}^{\op})^{\op}}{\Operad^{\Psi\wr\Phi}}
\end{equation*}
each admit a right adjoint, denoted $\Alg^{\Phi,\Psi\wr\Phi}(X,-)$ and $\Alg^{\Psi,\Psi\wr\Phi}(Y,-)$, respectively.
\end{prp}

\begin{nul} Hence for any two operator categories $\Phi$ and $\Psi$, we obtain a pair of functors
\begin{equation*}
\Alg^{\Psi,\Psi\wr\Phi}\colon\fromto{\mathbf{LFib}(N\Delta_{\Psi}^{\op})^{\op}\times\Operad^{\Psi\wr\Phi}}{\mathbf{LFib}(N\Delta_{\Phi}^{\op})}
\end{equation*}
and
\begin{equation*}
\Alg^{\Phi,\Psi\wr\Phi}\colon\fromto{\mathbf{LFib}(N\Delta_{\Phi}^{\op})^{\op}\times\Operad^{\Psi\wr\Phi}}{\mathbf{LFib}(N\Delta_{\Psi}^{\op})}.
\end{equation*}

In fact, these functors are valued in the $\infty$-category of weak operads. That is, suppose $\Phi$ and $\Psi$ operator categories, suppose
\begin{equation*}
\fromto{X}{N\Delta^{\op}_{\Phi}}\textrm{\quad and\quad}\fromto{Y}{N\Delta^{\op}_{\Psi}}
\end{equation*}
left fibrations, and suppose $Z$ a complete Segal $\Psi\wr\Phi$-operad. Then the left fibration $\fromto{\Alg^{\Psi,\Psi\wr\Phi}(Y,Z)}{N\Delta^{\op}_{\Phi}}$ is a complete Segal $\Phi$-operad, and the left fibration $\fromto{\Alg^{\Phi,\Psi\wr\Phi}(X,Z)}{N\Delta^{\op}_{\Psi}}$ is a complete Segal $\Psi$-operad.

This is not an entirely formal matter. It suffices to prove this for $X$ and $Y$ corepresentable. In this case, one shows that forming the Boardman--Vogt tensor product of any corepresentable left fibration with any element of the set $S_{\Phi}$ (or $S_{\Psi}$) of Pr. \ref{prp:wkoperadmodstruct} itself lies in the strongly saturated class generated by $S_{\Psi\wr\Phi}$. We leave these details to the reader.
\end{nul}

Now let's concentrate on the situation in which one of the two complete Segal operads is the terminal operad.

\begin{ntn} Suppose $\Phi$ and $\Psi$ two operator categories. For any complete Segal $\Psi\wr\Phi$-operad $Z$, we will write $\Mon^{\Phi,\Psi\wr\Phi}(Z)$ for $\Alg^{\Phi,\Psi\wr\Phi}(U_{\Phi},Z)$, and we will write $\Mon^{\Phi,\Psi\wr\Phi}(Z)$ for $\Alg^{\Psi,\Psi\wr\Phi}(U_{\Psi},Z)$.
\end{ntn}

\begin{nul} Suppose $\Phi$ an operator category. Note that we have canonical equivalences $\Phi\wr\{1\}\simeq\Phi\simeq\{1\}\wr\Phi$, through which the functor $\Alg^{\Phi,\Phi\wr\{1\}}$ and $\Alg^{\Phi,\{1\}\wr\Phi}$ may be identified. We write
\begin{equation*}
\Alg^{\Phi}\colon\fromto{\mathbf{LFib}(N\Delta^{\op}_{\Phi})^{\op}\times\Operad^{\Phi}}{\mathbf{LFib}(N\Delta^{\op})}
\end{equation*}
for the common functor. For any left fibration $\fromto{Y}{N\Delta^{\op}_{\Phi}}$ and any complete Segal $\Phi$-operad $Z$, we will refer to objects of $\Alg^{\Phi}(Y,Z)$ (i.e., those $0$-simplices that lie over $0\in N\Delta^{\op}$) as \textbf{\emph{$Y$-algebras in $Z$}}. If $Y$ is the terminal $\Phi$-operad $U_{\Phi}$, then we will refer to $Y$-algebras in $Z$ as \textbf{\emph{complete Segal $\Phi$-monoids in $Z$}}.
\end{nul}

\begin{thm}\label{prp:UotimesUisU} For any two operator categories $\Phi$ and $\Psi$, the pairing
\begin{equation*}
\fromto{(U_{\Psi},U_{\Phi})}{U_{\Psi\wr\Phi}}
\end{equation*}
given by $W$ exhibits $U_{\Psi\wr\Phi}$ as the Boardman--Vogt tensor product of $U_{\Psi}$ and $U_{\Phi}$.
\end{thm}
\noindent The proof is a bit involved, so it appears in an appendix; see \S \ref{sect:proofthatUotimesUisU}.

\begin{cor} For any two operator categories $\Phi$ and $\Psi$, one has, for any $\Psi\wr\Phi$-operad $Z$, canonical equivalences
\begin{equation*}
\Mon^{\Psi\wr\Phi}(Z)\simeq\Mon^{\Phi}(\Mon^{\Psi,\Psi\wr\Phi}(Z))\simeq\Mon^{\Psi}(\Mon^{\Phi,\Psi\wr\Phi}(Z)).
\end{equation*}
\end{cor} 
\noindent Roughly speaking, we have shown that complete Segal $\Psi\wr\Phi$-monoids are complete Segal $\Psi$-monoids in complete Segal $\Phi$-monoids, which in turn are complete Segal $\Phi$-monoids in complete Segal $\Psi$-monoids.

\begin{exm} The $(\infty,1)$-category of complete Segal $\OO^{(n)}$-monoids (in some complete Segal $\OO^{(n)}$-operad $Z$) is equivalent to the $(\infty,1)$-category
\begin{equation*}
\Mon^{\OO,\OO}(\Mon^{\OO,\OO\wr\OO}(\cdots\Mon^{\OO,\OO^{(n)}}(Z)\cdots))
\end{equation*}
of complete Segal $\OO$-monoids in complete Segal $\OO$-monoids in \dots\ in complete Segal $\OO$-monoids in $Z$. The assertion that the operad $E_n$ is equivalent to the symmetrization of the terminal $\OO^{(n)}$-operad (Pr. \ref{prp:mainEn}) thus states that if $Z$ is a symmetric operad, then the homotopy theory of $E_n$-algebras in $Z$ is equivalent to the homotopy theory of complete Segal $\OO$-monoids in complete Segal $\OO$-monoids in \dots\ in complete Segal $\OO$-monoids in the underlying complete Segal $\OO^{(n)}$-operad of $Z$.
\end{exm}

\begin{exm} For any integers $n\geq m\geq 0$, we have an inclusion
\begin{equation*}
s\colon\into{\OO^{(m)}\cong\{1\}^{(n-m)}\wr\OO^{(m)}}{\OO^{(n-m)}\wr\OO^{(m)}\cong\OO^{(n)}},
\end{equation*}
which is a section of the coronal fibration $p\colon\fromto{\OO^{(n)}}{\OO^{(m)}}$. We can thus form the colimit $\OO^{(\infty)}\coloneq\colim_{n\geq 0}\OO^{(n)}$. This can be viewed as the category whose objects are sequences
\begin{equation*}
(M_1,M_2,\dots)
\end{equation*}
with $M_n\in\OO^{(n)}$ such that for every for each $n\geq1$, one has $p(M_n)=M_{n-1}$ and for each $n\gg1$, one has $s(M_n)=M_{n+1}$. Then $\OO^{(\infty)}$ is an operator category.
\end{exm}


\section{Boardman--Vogt tensor products and $\infty$-algebras} We can define an analogue of the Boardman--Vogt tensor product introduced in \S \ref{sect:BV} and study its interaction with the model structure introduced in he previous section. Once again, in most cases the proofs are trivial extensions of the proofs of loc. cit.

\begin{ntn} For any perfect operator categories $\Phi$ and $\Psi$, define a natural transformation
\begin{equation*}
\omega\colon\fromto{\wr\circ(T_{\Psi}\times T_{\Phi})}{T_{\Psi\wr\Phi}\circ\wr}
\end{equation*}
as follows. For any pair $(K,I)\in\Psi\times\Phi$, let
\begin{equation*}
\omega_{(K,I)}=((\phi_j)_{j\in|T_{\Phi}I|},\id_{T_{\Phi}I})\colon\fromto{T_{\Psi}K\wr T_{\Phi}I}{T_{\Psi\wr\Phi}(K\wr I)}
\end{equation*}
be the morphism in which
\begin{equation*}
\phi_j=\begin{cases}
\id_{T_{\Psi}K}&\textrm{if }j\in|I|\subset|T_{\Phi}I|;\\
!&\textrm{if }j\notin|I|.
\end{cases}
\end{equation*}
Using this, we obtain an induced functor $W\colon\fromto{\Lambda(\Psi)\times\Lambda(\Phi)}{\Lambda(\Psi\wr\Phi)}$ on the Leinster categories given by the assignment $\goesto{(K,I)}{K\wr I}$.
\end{ntn}

\begin{dfn} Suppose $\Phi$ and $\Psi$ two perfect operator categories, suppose $X\in s\Set^{+}_{/N\Lambda(\Phi)}$, and suppose $Y\in s\Set^{+}_{/N\Lambda(\Psi)}$. Write $Y{{}^{\Psi}\otimes^{\Phi}}X$ for the product $Y\times X$, regarded as a marked simplicial set over $N\Lambda(\Psi\wr\Phi)$ via
\begin{equation*}
Y\times X\to N\Lambda(\Psi)\times N\Lambda(\Phi)\ \tikz[baseline]\draw[>=stealth,->,font=\scriptsize](0,0.5ex)--node[above]{$W$}(0.5,0.5ex);\ N\Lambda(\Psi\wr\Phi).
\end{equation*}
We call this the \textbf{\emph{Boardman--Vogt tensor product}} of $Y$ and $X$.
\end{dfn}

\begin{nul}\label{nul:BVotimesFF} Consider the operator category $\FF$. To relate our Boardman--Vogt tensor product to the monoidal structure $\odot$ constructed in \cite{HA}, we need only note that the functor $\wedge\colon\fromto{\Gamma^{\op}\times\Gamma^{\op}}{\Gamma^{\op}}$ of \cite{HA} is isomorphic to the composition of the functor $W\colon\fromto{\Lambda(\FF)\times\Lambda(\FF)}{\Lambda(\FF\wr\FF)}$ with the functor $\fromto{\Lambda(\FF\wr\FF)}{\Lambda(\FF)}$ induced by the unique operator morphism $U\colon\fromto{\FF\wr\FF}{\FF}$. Consequently, we obtain an isomorphism $U_{!}(Y{{}^{\FF}\otimes^{\FF}}X)\cong Y\odot X$.
\end{nul}

\begin{prp} For any operator morphisms $H\colon\fromto{\Psi'}{\Psi}$ and $G\colon\fromto{\Phi'}{\Phi}$ between perfect operator categories, and for any objects $Y\in s\Set_{/N\Lambda(\Psi')}$ and $X\in s\Set_{/N\Lambda(\Phi')}$, there is a canonical isomorphism
\begin{equation*}
(H\wr G)_!(Y{{}^{\Psi'}\otimes^{\Phi'}}X)\cong H_!Y{{}^{\Psi}\otimes^{\Phi}}G_!X
\end{equation*}
of simplicial sets over $N\Lambda(\Psi\wr\Phi)$.
\end{prp}

The Boardman--Vogt tensor product preserves colimits separately in each variable. Consequently, we have the following.
\begin{prp} Suppose $\Phi$ and $\Psi$ two perfect operator categories, suppose $X\in s\Set_{/N\Lambda(\Phi)}$, and suppose $Y\in s\Set_{/N\Lambda(\Psi)}$. Then the functors
\begin{equation*}
Y{{}^{\Psi}\otimes^{\Phi}}-\colon\fromto{s\Set^{+}_{/N\Lambda(\Phi)}}{s\Set^{+}_{N\Lambda(\Psi\wr\Phi)}}\textrm{\; and\;}-{{}^{\Psi}\otimes^{\Phi}}X\colon\fromto{s\Set^{+}_{/N\Lambda(\Psi)}}{s\Set^{+}_{N\Lambda(\Psi\wr\Phi)}}
\end{equation*}
both admits right adjoints.
\end{prp}

\begin{ntn} For any two perfect operator categories $\Phi$ and $\Psi$, we obtain a pair of functors
\begin{equation*}
\Alg_{\infty}^{\Psi,\Psi\wr\Phi}\colon\fromto{(s\Set^{+}_{/N\Lambda(\Psi)})^{\op}\times s\Set^{+}_{/N\Lambda(\Psi\wr\Phi)}}{s\Set^{+}_{/N\Lambda(\Phi)}}
\end{equation*}
and
\begin{equation*}
\Alg_{\infty}^{\Phi,\Psi\wr\Phi}\colon\fromto{(s\Set^{+}_{/N\Lambda(\Phi)})^{\op}\times s\Set^{+}_{/N\Lambda(\Psi\wr\Phi)}}{s\Set^{+}_{/N\Lambda(\Psi)}}
\end{equation*}
such that $\Alg_{\infty}^{\Psi,\Psi\wr\Phi}(Y,-)$ is right adjoint to $Y{{}^{\Psi}\otimes^{\Phi}}-$ and $\Alg_{\infty}^{\Phi,\Psi\wr\Phi}(X,-)$ is right adjoint to $-{{}^{\Psi}\otimes^{\Phi}}X$.

It is clear that these functors, along with the Boardman--Vogt tensor product, comprise an adjunction of two variables $\fromto{s\Set^{+}_{/N\Lambda(\Psi)}\times s\Set^{+}_{/N\Lambda(\Phi)}}{s\Set^{+}_{/N\Lambda(\Psi\wr\Phi)}}$.
\end{ntn}

The interaction between this variant of the Boardman--Vogt tensor product and the operadic model structure is the most one could hope for. We apply \cite[Rk. B.2.5 and Pr. B.2.9]{HA} to the functors $W\colon\fromto{\Lambda(\Psi)\times\Lambda(\Phi)}{\Lambda(\Psi\wr\Phi)}$ to deduce the following.
\begin{thm} For any two perfect operator categories $\Phi$ and $\Psi$, the functors ${}^{\Psi}\otimes^{\Phi}$, $\Alg_{\infty}^{\Psi,\Psi\wr\Phi}$, and $\Alg_{\infty}^{\Phi,\Psi\wr\Phi}$ form a Quillen adjuction of two variables for the operadic model structures.
\end{thm}

\begin{ntn} Suppose $\Phi$ and $\Psi$ two perfect operator categories. For any $(\Psi\wr\Phi)$- quasioperad $Z$, we will write $\Mon_{\infty}^{\Phi,\Psi\wr\Phi}(Z)$ for $\Alg_{\infty}^{\Phi,\Psi\wr\Phi}(U_{\Phi},Z)$, and we will write $\Mon_{\infty}^{\Phi,\Psi\wr\Phi}(Z)$ for $\Alg_{\infty}^{\Psi,\Psi\wr\Phi}(U_{\Psi},Z)$.
\end{ntn}

\begin{nul} Suppose $\Phi$ a perfect operator category. Note that we have canonical equivalences $\Phi\wr\{1\}\simeq\Phi\simeq\{1\}\wr\Phi$, through which the functor $\Alg_{\infty}^{\Phi,\Phi\wr\{1\}}$ and $\Alg_{\infty}^{\Phi,\{1\}\wr\Phi}$ may be identified. We write
\begin{equation*}
\Alg_{\infty}^{\Phi}\colon\fromto{(s\Set^{+}_{/N\Lambda(\Phi)})^{\op}\times s\Set^{+}_{/N\Lambda(\Phi)}}{s\Set^{+}}
\end{equation*}
for the common functor. For any marked map $\fromto{Y}{N\Lambda(\Phi)}$ and any $\infty$-operad $Z$ over $\Phi$, we will refer to objects of the quasicategory $\Alg_{\infty}^{\Phi}(Y,Z)$ as \textbf{\emph{$\infty$-algebras over $Y$ in $Z$}}. If $Y$ is the terminal $\Phi$-operad $U_{\Phi}$, then we will refer to $Y$-algebras in $Z$ as \textbf{\emph{$\infty$-monoids over $\Phi$ in $Z$}}.
\end{nul}

The following theorem is proved exactly as in \cite[Th. 2.4.4.3]{HA}.
\begin{thm}\label{prp:UotimesUisUredux} For any two operator categories $\Phi$ and $\Psi$, the functor $W$ induces an equivalence $\equivto{U_{\Psi}{{}^{\Psi}\otimes^{\Phi}}U_{\Phi}}{U_{\Psi\wr\Phi}}$.
\end{thm}

\begin{cor} For any two perfect operator categories $\Phi$ and $\Psi$, one has, for any $\Psi\wr\Phi$-operad $Z$, canonical equivalences
\begin{equation*}
\Mon_{\infty}^{\Psi\wr\Phi}(Z)\simeq\Mon_{\infty}^{\Phi}(\Mon_{\infty}^{\Psi,\Psi\wr\Phi}(Z))\simeq\Mon_{\infty}^{\Psi}(\Mon_{\infty}^{\Phi,\Psi\wr\Phi}(Z)).
\end{equation*}
\end{cor} 


\section{Complete Segal operads versus quasioperads}\label{sect:equivalence} We now compare our two approaches to the theory of complete Segal $\Phi$-operads. Here is the main theorem.

\begin{thm*} For any perfect operator category $\Phi$, there exist inverse equivalences of $\infty$-categories
\begin{equation*}
P^{\otimes}\colon\Operad^{\Phi}_{\mathrm{CSS}}\simeq\Operad^{\Phi}_{\infty}\colon C.
\end{equation*}
\end{thm*}
\noindent We do this by obtaining a homotopy equivalence of relative categories (Th. \ref{thm:sameasLurie}).

We begin with the functor $P^{\otimes}$. To construct this functor, we must first relate the Leinster category $\Lambda(\Phi)$ to the category of $\Phi$-sequences $\Delta_{\Phi}$.

\begin{ntn} Suppose $\Phi$ a perfect operator category. For any object $\mathbf{m}\in\Delta$, denote by $\widetilde{\mathscr{O}}(\mathbf{m})$ the twisted arrow category of $\mathbf{m}$; this is a poset whose objects are pairs of integers $(i,j)$ such that $0\leq i\leq j\leq m$, and $(i',j')\leq(i,j)$ just in case $i\leq i'\leq j'\leq j$.

Now for any $m$-simplex $\sigma\colon\fromto{\mathbf{m}}{N\Lambda(\Phi)}$ corresponding to a sequence of morphisms
\begin{equation*}
I_0\to I_1\to\cdots\to I_m
\end{equation*}
in $\Lambda(\Phi)$, we obtain a functor $F_{\sigma}\colon\fromto{\widetilde{\mathscr{O}}(\mathbf{m})^{\op}}{\Delta_{\Phi}}$, given by the formula
\begin{equation*}
F_{\sigma}(i,j)\coloneq[I_0\times_{TI_j}I_j\to I_1\times_{TI_j}I_j\to\cdots\to I_i\times_{TI_j}I_j].
\end{equation*}
We write $\Delta_{\Phi}^{\sigma}$ for the colimit of the diagram $j\circ F_{\sigma}$, where $j$ denotes the Yoneda embedding $\into{\Delta_{\Phi}}{\Fun(\Delta_{\Phi}^{\op},\Set)}$.

The assignment $\goesto{\sigma}{\Delta_{\Phi}^{\sigma}}$ is functorial with respect to the category $\mathrm{Simp}(\Lambda(\Phi))$, whose objects are pairs $(\mathbf{m},\sigma)$ consisting of an object $\mathbf{m}\in\Delta$ and an $m$-simplex $\sigma\colon\fromto{\mathbf{m}}{N\Lambda(\Phi)}$ of $\Lambda(\Phi)$ and whose morphisms $\fromto{(\mathbf{m},\sigma)}{(\mathbf{n},\tau)}$ are morphisms $\eta\colon\fromto{\mathbf{m}}{\mathbf{n}}$ of $\Delta$ and a natural isomorphism $\sigma\cong\tau\circ\eta$; in particular, we obtain a functor
\begin{equation*}
\Delta_{\Phi}^{\ast}\colon\fromto{\iota\Fun(\mathbf{m},\Lambda(\Phi))}{\iota\Fun(\Delta_{\Phi}^{\op},\Set)}.
\end{equation*}
\end{ntn}

\begin{cnstr} Suppose $\Phi$ a perfect operator category, and suppose that $X\colon\fromto{\Delta_{\Phi}^{\op}}{s\Set}$ a functor. For any integer $m\geq 0$ and any $m$-simplex
\begin{equation*}
\sigma\colon\fromto{\mathbf{m}}{N\Lambda(\Phi)},
\end{equation*}
one may define a simplicial set $X(\Delta_{\Phi}^{\sigma})$ via right Kan extension, so that
\begin{equation*}
X(\Delta_{\Phi}^{\sigma})\cong \lim_{(i,j)\in\widetilde{\mathscr{O}}(\mathbf{m})}X(F_{\sigma}(i,j)).
\end{equation*}
The assignment $\goesto{\sigma}{X(\Delta_{\Phi}^{\sigma})}$ defines a functor $\fromto{\iota\Fun(\mathbf{m},\Lambda(\Phi))^{\op}}{s\Set}$. Apply the construction of \cite[Df. 3.2.5.2]{HTT} to this functor to obtain a map
\begin{equation*}
\fromto{\PP^{\otimes}(X)_m}{N\iota\Fun(\mathbf{m},\Lambda(\Phi))^{\op}}.
\end{equation*}
This map is functorial in $\mathbf{m}\in\Delta^{\op}$, whence we obtain a morphism of simplicial spaces
\begin{equation*}
\fromto{\PP^{\otimes}(X)}{\NN\Lambda(\Phi)},
\end{equation*}
where $\NN\Lambda(\Phi)$ is the classifying diagram of the category $\Lambda(\Phi)$ in the sense of Rezk \cite[3.5]{MR1804411} (except that in each degree we are taking the nerve of the opposite groupoid).
\end{cnstr}

\begin{prp}\label{prp:PPisCSS} Suppose $\Phi$ a perfect operator category, and suppose that
\begin{equation*}
X\colon\fromto{\Delta_{\Phi}^{\op}}{\Kan}
\end{equation*}
a functor. Assume that $X$ that is fibrant for the injective model structure on the category $\Fun(\Delta_{\Phi}^{\op},s\Set)$, and that $X$ classifies a complete Segal $\Phi$-operad. Then $\PP^{\otimes}(X)$ is a complete Segal space.
\begin{proof} To see that $\PP^{\otimes}(X)$ is Reedy fibrant, consider any square
\begin{equation*}
\begin{tikzpicture} 
\matrix(m)[matrix of math nodes, 
row sep=4ex, column sep=4ex, 
text height=1.5ex, text depth=0.25ex] 
{\Lambda^n_{k}&\PP^{\otimes}(X)_m\\ 
\Delta^n&\lim_{\mathbf{k}\subsetneqq\mathbf{m}}\PP^{\otimes}(X)_k.\\}; 
\path[>=stealth,->,font=\scriptsize] 
(m-1-1) edge (m-1-2) 
edge (m-2-1) 
(m-1-2) edge (m-2-2) 
(m-2-1) edge (m-2-2); 
\end{tikzpicture}
\end{equation*}
Unwinding the definitions, we see that a lift $\fromto{\Delta^n}{\PP^{\otimes}(X)_m}$ amounts to a lift $\sigma$ of the diagram
\begin{equation*}
\begin{tikzpicture} 
\matrix(m)[matrix of math nodes, 
row sep=4ex, column sep=4ex, 
text height=1.5ex, text depth=0.25ex] 
{\Lambda^n_{k}&N\iota\Fun(\mathbf{m},\Lambda(\Phi))^{\op}\\ 
\Delta^n&\lim_{\mathbf{k}\subsetneqq\mathbf{m}}N\iota\Fun(\mathbf{k},\Lambda(\Phi))^{\op}\\}; 
\path[>=stealth,->,font=\scriptsize] 
(m-1-1) edge (m-1-2) 
edge (m-2-1) 
(m-1-2) edge (m-2-2) 
(m-2-1) edge (m-2-2); 
\end{tikzpicture}
\end{equation*}
along with a compatible collection of maps $\fromto{\Delta^{J}}{X(\Delta^{\sigma(\max J)}_{\Phi})}$, one for each nonempty subset $J\subset\mathbf{m}$. The existence of the lift $\sigma$ follows from \cite[Lm. 3.9]{MR1804411}, and the compatible collection of maps follows from the injective fibrancy.

The Segal conditions on $\PP^{\otimes}(X)$ reduce to showing that for any $m$-simplex $\sigma\colon\fromto{\mathbf{m}}{\Lambda(\Phi)}$, the induced map
\begin{equation*}
\fromto{X(\Delta^{\sigma}_{\Phi})}{X(\Delta^{\sigma|\{0,1\}}_{\Phi})\times_{X(\Delta^{\sigma|\{1\}}_{\Phi})}\cdots\times_{X(\Delta^{\sigma|\{m-1\}}_{\Phi})}X(\Delta^{\sigma|\{m-1,m\}}_{\Phi})}
\end{equation*}
is a weak equivalence. But this follows from the fact that, for $X$ itself, the maps $s_{\mathbf{n},J}$ are all equivalences.

The completeness conditions on $\PP^{\otimes}(X)$ reduces to the assertion that for any object $I\in\Phi$, the natural map
\begin{equation*}
\fromto{X[I]}{X[I=I=I=I]\times_{X[I=I]\times X[I=I]}X[I]}
\end{equation*}
is a weak equivalence, where the maps $\fromto{X[I=I=I=I]}{X[I=I]}$ are given by the inclusions $\into{\{0,2\}}{\{0,1,2,3\}}$ and $\into{\{1,3\}}{\{0,1,2,3\}}$. But this follows from the case when $I=\{1\}$ along with the decomposition
\begin{equation*}
X[I]\simeq\prod_{i\in|I|}X[\{i\}].\qedhere
\end{equation*}
\end{proof}
\end{prp}

The results of \cite{MR2342834} thus motivate the following.
\begin{dfn} Suppose $\Phi$ a perfect operator category, and suppose
\begin{equation*}
X\colon\fromto{\Delta_{\Phi}^{\op}}{s\Set}
\end{equation*}
a functor. Then we define a simplicial set $P^{\otimes}(X)$ in the following manner. For any $\mathbf{m}\in\Delta^{\op}$, an $m$-simplex $(\sigma,x)$ of $P^{\otimes}(X)$ consists of a functor $\sigma\colon\fromto{\mathbf{m}}{\Lambda(\Phi)}$ and a vertex $x\in X(\Delta_{\Phi}^{\sigma})_0$. This is obviously functorial in $\mathbf{m}$, and the assignment $\goesto{(\sigma,x)}{\sigma}$ defines a projection $\fromto{P^{\otimes}(X)}{N\Lambda(\Phi)}$. Hence we obtain a functor $P^{\otimes}\colon\fromto{\Fun(\Delta_{\Phi}^{\op},s\Set)}{s\Set_{/N\Lambda(\Phi)}}$.
\end{dfn}

\begin{nul} Suppose $\Phi$ a perfect operator category, and suppose
\begin{equation*}
X\colon\fromto{\Delta_{\Phi}^{\op}}{s\Set}
\end{equation*}
a functor. Note that the fiber of $P^{\otimes}(X)$ over an object $I\in\Lambda(\Phi)$ is the simpicial set whose $m$-simplices are the vertices of $X(m,I)$, where $(m,I)$ denotes the constant sequence $I=I=\cdots=I$ of length $m$. In particular, the fiber over $\{1\}$ is the simplicial set whose $m$-simplices are the vertices of the simplicial set $(X|\Delta^{\op})_m$.
\end{nul}

\begin{prp} Suppose $\Phi$ a perfect operator category, and suppose 
\begin{equation*}
X\colon\fromto{\Delta_{\Phi}^{\op}}{s\Set}
\end{equation*}
a functor. Assume that $X$ is fibrant for the operadic model structure on the category $\Fun(\Delta_{\Phi}^{\op},s\Set)$. Then $P^{\otimes}(X)$ is a $\Phi$-quasioperad.
\begin{proof} Since the $m$-simplices of $P^{\otimes}(X)$ can be identified with the $0$-simplices of $\PP^{\otimes}(X)_m$, Pr \ref{prp:PPisCSS} and \cite[Th. 4.11]{MR2342834} together imply that $P^{\otimes}(X)$ is a quasicategory, whence by \cite[Pr. 2.3.1.5]{HTT}, $\fromto{P^{\otimes}(X)}{\Lambda(\Phi)}$ is an inner fibration.

Suppose $\fromto{J}{I}$ an inert morphism of $\Lambda(\Phi)$, and suppose that $x\in P^{\otimes}(X)_J$ an object; hence $x\in X[J]_0$. Let $x'\in X[J\times_{TI}I]$ be its image under the map induced by the interval inclusion $\into{J\times_{TI}I}{J}$. Now we may use the inverse to the isomorphism $\fromto{J\times_{TI}I}{I}$ to define a morphism of $\Phi$-sequences
\begin{equation*}
\fromto{[J\times_{TI}I\to I]}{[J\times_{TI}I]},
\end{equation*}
and one may consider the image $x''$ of $x'$ in $X[J\times_{TI}I\to I]$ under the induced map. The pair $(x,x'')\in X[J]\times_{X[J\times_{TI}I]}X[J\times_{TI}I\to I]$ is now a cocartesian edge covering $\phi$.

Suppose that the following are given:
\begin{itemize}
\item[---] two objects $I,J\in\Phi$,
\item[---] vertices $x\in P^{\otimes}(X)_{I}$ and $y\in P^{\otimes}(X)_{J}$,
\item[---] a morphism $\phi\colon\fromto{J}{I}$ of $\Lambda(\Phi)$, and
\item[---] $p$-cocartesian edges $\{\fromto{y}{y_i}\ |\ i\in|I|\}$ lying over the inert morphisms $\{\rho_i\colon\fromto{I}{\{i\}}\ |\ i\in|I|\}$.
\end{itemize}
Now the simplicial set $\Map_{P^{\otimes}(X)}^{\phi}(x,y)$ can be identified with the fiber of the natural map
\begin{equation*}
\fromto{X[J]\times_{X[J\times_{TI}I]}X[J\times_{TI}I\to I]}{X[J]\times X[I]}
\end{equation*}
over the point $(x,y)$. The fact that the induced map
\begin{equation*}
\fromto{\Map_{X^{\otimes}}^{\phi}(x,y)}{\prod_{i\in|I|}\Map_{X^{\otimes}}^{\rho_i\circ\phi}(x,y_i)}
\end{equation*}
is an equivalence thus follows from the decomposition
\begin{eqnarray}
X[J]\times_{X[J\times_{TI}I]}X[J\times_{TI}I\to I]&\simeq&X[J]\times_{\prod_{i\in|I|}X[J_i]}\prod_{i\in|I|}X[J_i\to\{i\}]\nonumber\\
&\to&X[J]\times \prod_{i\in|I|}X[\{i\}]\nonumber\\
&\simeq&X[J]\times X[I].\nonumber
\end{eqnarray}

Finally, we must show that for any object $I\in\Phi$ and for any $I$-tuple
\begin{equation*}
(y_i)_{i\in|I|}\in\prod_{i\in|I|}P^{\otimes}(X)_{\{i\}},
\end{equation*}
there exists an object $y\in P^{\otimes}(X)_{I}$ a collection of cocartesian edges $\fromto{y}{y_i}$ lying over the inert edges $\fromto{I}{\{i\}}$. For this, choose $y\in X[I]$ that maps to $(y_i)_{i\in|I|}$. (This is possible since the map $\fromto{X[I]}{\prod_{i\in|I|}X[\{i\}]}$ is a trivial fibration.) For any point $i\in|I|$, the image of $y$ under the natural map
\begin{equation*}
\fromto{X[I]\cong X[I]\times_{X[\{i\}]}X[\{i\}]}{X[I]\times_{X[\{i\}]}X[\{i\}=\{i\}]}
\end{equation*}
is the desired cocartesian edge $\fromto{y}{y_i}$.
\end{proof} 
\end{prp}

\begin{exm} When $\Phi=\{1\}$, the functor $P^{\otimes}$ carries a complete Segal space $X$ to the quasicategory whose $m$-simplices are the $0$-simplices of $X_m$. By a theorem of Joyal and Tierney \cite[Th. 4.11]{MR2342834}, this is known to be an equivalence of homotopy theories.
\end{exm}

It is obvious from the construction given here that $P^{\otimes}$ is compatible with changes of operator category, in the following sense.
\begin{prp} Suppose $G\colon\fromto{\Psi}{\Phi}$ is an operator morphism between perfect operator categories. Then we have a natural isomorphism of functors
\begin{equation*}
P^{\otimes}\circ G^{\star}\cong G^{\star}\circ P^{\otimes}.
\end{equation*}
\end{prp}

\begin{prp} Suppose $\Phi$ a perfect operator category. A morphism $\fromto{X}{Y}$ between injectively fibrant complete Segal $\Phi$-operads is a weak equivalence just in case the induced morphism $\fromto{P^{\otimes}(X)}{P^{\otimes}(Y)}$ is so.
\begin{proof} Let us use the characterizations of Prs. \ref{prp:DKequivs} and \ref{prp:DKequivsonquasiops}. Combining the example and the proposition above, we deduce that a morphism $\fromto{X}{Y}$ of complete Segal $\Phi$-operads induces an essentially surjective functor on underlying complete Segal spaces just in case the morphism $\fromto{P^{\otimes}(X)}{P^{\otimes}(Y)}$ induces an essentially surjective functor on underlying quasicategories.

On the other hand, if $I\in\Phi$, and if $\alpha\colon\fromto{I}{\{1\}}$ is the morphism of $\Lambda(\Phi)$ given by $I\to\{t\}\to T$, then for any complete Segal $\Phi$-operad $X$, the mapping space $\Map_X^{\alpha}(x,y)$ is equivalent to the fiber of the map
\begin{equation*}
\fromto{X[I\to\{1\}]}{X[I]\times X[\{1\}]}.
\end{equation*}
In particular, a morphism $\fromto{X}{Y}$ of complete Segal $\Phi$-operads is fully faithful just in case the induced morphism $\fromto{P^{\otimes}(X)}{P^{\otimes}(Y)}$ of $\Phi$-quasicategories is so.
\end{proof}
\end{prp}
\noindent In other words, $P^{\otimes}$ is a relative functor in the sense of \cite{MR2877401}, and it reflects weak equivalences.

To define a functor in the opposite direction, we introduce the following.
\begin{cnstr} Suppose $\Phi$ a perfect operator category, and suppose $(\mathbf{m},I)$ a $\Phi$-sequence. Define a poset $A(\mathbf{m},I)$ as follows: the objects are triples $(r,s,i)$, where $r$ and $s$ are integers such that $0\leq r\leq s\leq m$, and $i\in|I_s|$; we declare that $(r,s,i)\leq (r',s',i')$ just in case $r\leq r'\leq s'\leq s$, and $\goesto{i'}{i}$ under the map $\fromto{|I_{s'}|}{|I_{s}|}$.

Let $A(\mathbf{m},I)^{\dag}$ be the set of morphisms of $A(\mathbf{m},I)$ of the form $\fromto{(r,s,i)}{(r,s',i')}$. Using the factorization system on $\Lambda(\Phi)$ of Lm. \ref{lem:actinert}, we deduce that there is a functor $\fromto{A(\mathbf{m},I)}{\Lambda(\Phi)}$ given by the assignment
\begin{equation*}
\goesto{(r,s,i)}{I_{r,i}\coloneq I_r\times_{I_s}\{i\}}
\end{equation*}
under which morphisms of $A(\mathbf{m},I)^{\dag}$ are carried to inert morphisms of $\Lambda(\Phi)$.

Now the assignment $\fromto{(\mathbf{m},I)}{(NA(\mathbf{m},I),A(\mathbf{m},I)^{\dag})}$ defines a functor
\begin{equation*}
NA\colon\fromto{\Delta_{\Phi}}{s\Set^{+}_{/N\Lambda(\Phi)}}.
\end{equation*}
For any object $X^{\otimes}\in s\Set^{+}_{/N\Lambda(\Phi)}$, write $C(X)$ for the functor $\fromto{\Delta_{\Phi}^{\op}}{s\Set}$ given by the assignment
\begin{equation*}
\goesto{(\mathbf{m},I)}{\Map^{\sharp}_{N\Lambda(\Phi)}(NA(\mathbf{m},I),X^{\otimes})}.
\end{equation*}
\end{cnstr}

\begin{exm}\label{exm:CwhenPhiisone} When $\Phi=\{1\}$, the category $A(\mathbf{m},\{1\})$ is the opposite of the twisted arrow category $\widetilde{\mathscr{O}}(\mathbf{m})$, and in $N\widetilde{\mathscr{O}}(\mathbf{m})^{\op}$, the marked edges are precisely the morphisms $\fromto{(r,s)}{(r,s')}$. Now the inclusion $\into{\mathbf{m}}{N\widetilde{\mathscr{O}}(\mathbf{m})^{\op}}$ given by $\goesto{r}{(r,m)}$ induces a marked anodyne morphism on nerves. Consequently, for any quasicategory $X$, one has $C(X)_m\simeq\iota\Fun(\Delta^m,X)$, which is an inverse homotopy equivalence of $P^{\otimes}$ by the theorem of Joyal and Tierney \cite[Th. 4.11]{MR2342834}.
\end{exm}

The construction $\goesto{(m,I)}{NA(m,I)}$ is compatible with changes of operator category. That is, for any operator morphism $G\colon\fromto{\Psi}{\Phi}$ between perfect operator categories, one has a canonical isomorphism $G_{!}NA(m,I)\cong NA(m,GI)$. Consequently, we have the following.
\begin{lem}\label{lem:GCCG} For any operator morphism $G\colon\fromto{\Psi}{\Phi}$ between perfect operator categories, there is a natural isomorphism
\begin{equation*}
C\circ G^{\star}\cong G^{\star}\circ C.
\end{equation*}
\end{lem}

\begin{prp} Suppose $\Phi$ a perfect operator category, and suppose $X^{\otimes}$ a $\Phi$-quasioperad. Then $C(X^{\otimes})$ is fibrant for the operadic model structure.
\begin{proof} The fact that $C(X^{\otimes})$ is fibrant for the injective model structure follows from \cite[Rk. A.2.9.27]{HTT} and the observation that $\goesto{(m,I)}{(NA(m,I),A(m,I)^{\dag})}$ is cofibrant for the projective model structure on $\Fun(\Delta_{\Phi},s\Set^{+}_{/N\Lambda(\Phi)})$.

For any $\Phi$-sequence $[I_0\to\cdots\to I_m]$, the canonical map
\begin{equation*}
\fromto{\coprod_{i\in|I_m|}NA[I_{0,i}\to\cdots\to\{i\}]}{NA[I_{0}\to\cdots\to I_m]}
\end{equation*}
is an isomorphism by definition. Moreover, for any integer $0\leq k\leq m$, an elementary computation shows that the inclusion
\begin{equation*}
\into{NA[I_0\to\cdots\to I_k]\cup^{NA[I_k]}NA[I_k\to\cdots\to I_m]}{NA[I_{0}\to\cdots\to I_m]}
\end{equation*}
is $\mathscr{P}$-anodyne.

It thus remains only to observe that the underlying simplicial set $p^{\star}C(X^{\otimes})$ is a complete Segal space, and this follows from Lm. \ref{lem:GCCG} and Ex. \ref{exm:CwhenPhiisone}.
\end{proof}
\end{prp}

\begin{prp}\label{prp:Cisrelfun} Suppose $\Phi$ a perfect operator category. A morphism $\fromto{X^{\otimes}}{Y^{\otimes}}$ between $\Phi$-quasioperads is a weak equivalence just in case the induced morphism $\fromto{C(X^{\otimes})}{C(Y^{\otimes})}$ of complete Segal $\Phi$-operads is so.
\begin{proof} Let us use the characterizations of Prs. \ref{prp:DKequivs} and \ref{prp:DKequivsonquasiops}. Applying Lm. \ref{lem:GCCG} and Ex. \ref{exm:CwhenPhiisone}, we see that $\fromto{X^{\otimes}}{Y^{\otimes}}$ is essentially surjective just in case the morphism $\fromto{C(X^{\otimes})}{C(Y^{\otimes})}$ is so.

To complete the proof, let us note that for any object $I\in\Phi$, the inclusion $\into{\Delta^1}{NA[I\to\{1\}]}$ given by $\goesto{r}{(r,m,1)}$ is $\mathscr{P}$-anodyne. Consequently, we obtain an equivalence
\begin{equation*}
[I/\{1\}]\Map_{C(X^{\otimes})}((x_i)_{i\in|I|},y)\simeq\Map^{\alpha}_{X^{\otimes}}(x,y)
\end{equation*}
for any objects $x\in X^{\otimes}_I$, $y\in X^{\otimes}_{\{1\}}$, and any collection $\{\fromto{x}{x_i}\ |\ i\in|I|\}$ of inert morphisms.
\end{proof}
\end{prp}
\noindent In other words, $C$, when restricted to $\Phi$-quasioperads is a relative functor in the sense of \cite{MR2877401}, and it reflects weak equivalences. In fact, it is part of a Quillen pair.
\begin{cor} Suppose $\Phi$ a perfect operator category. Then the functors
\begin{equation*}
\adjunct{B}{\Fun(\Delta^{\op}_{\Phi},s\Set)}{s\Set^{+}_{/N\Lambda(\Phi)}}{C}
\end{equation*}
given by the formulas
\begin{equation*}
B(X)\coloneq\int^{(\mathbf{m},I)}NA(\mathbf{m},I)\times X(\mathbf{m},I)^{\sharp}
\end{equation*}
and
\begin{equation*}
C(Y^{\otimes})(\mathbf{m},I)\coloneq\Map_{N\Lambda(\Phi)}^{\sharp}(NA(\mathbf{m},I),Y^{\otimes})
\end{equation*}
form a Quillen pair for the operadic model structures.
\end{cor}

\begin{nul} Note that, using the functor $C$ along with the construction of \cite[Cnstr. 2.1.1.7]{HA}, we can now form the \textbf{\emph{nerve}} of any symmetric operad enriched in Kan complexes as a complete Segal $\FF$-operad. 
\end{nul}

We are now ready to state the main result of this section.
\begin{thm}\label{thm:sameasLurie} For any perfect operator category $\Phi$, the functors
\begin{equation*}
P^{\otimes}\colon\fromto{\Operad_{\mathrm{CSS}}^{\Phi}}{\Operad_{\infty}^{\Phi}}\textrm{\quad and\quad}C\colon\fromto{\Operad_{\infty}^{\Phi}}{\Operad_{\mathrm{CSS}}^{\Phi}}
\end{equation*}
are inverse equivalences of $\infty$-categories.
\begin{proof} Since we have shown that each functor is conservative, it is sufficient for us to furnish a natural transformation $\eta\colon\fromto{\id}{P^{\otimes}\circ C}$ that is objectwise essentially surjective and fully faithful.

For any object $\mathbf{m}\in\Delta$, any $m$-simplex $\sigma\colon\fromto{\mathbf{m}}{\Lambda(\Phi)}$, and any integers $0\leq i\leq j\leq m$, consider the functor $f_{\sigma,i,j}\colon\fromto{AF_{\sigma}(i,j)}{\mathbf{m}}$ over $\Lambda(\Phi)$ given by the assignment $\goesto{(r,s,k)}{r}$; these are compatible with one another and hence define a map
\begin{equation*}
f_{\sigma}\colon\colim_{(i,j)\in\widetilde{\mathscr{O}}(\mathbf{m})}NAF_{\sigma}(i,j)\to\Delta^m
\end{equation*}
in $s\Set^{+}_{/N\Lambda(\Phi)}$. This induces, for any $\Phi$-quasioperad $X^{\otimes}$, a map
\begin{equation*}
X^{\otimes}_m\to\Mor_{N\Lambda}(\colim_{(i,j)\in\widetilde{\mathscr{O}}(\mathbf{m})}NAF_{\sigma}(i,j),X^{\otimes})\cong P^{\otimes}(C(X^{\otimes}))_m,
\end{equation*}
natural in $X^{\otimes}$, hence a natural transformation $\eta\colon\fromto{\id}{P^{\otimes}\circ C}$.

When $\Phi=\{1\}$, the functor $f_{\sigma,i,j}$ is the projection $\fromto{\widetilde{\mathscr{O}}(\mathbf{m})^{\op}}{\mathbf{m}}$ given by $\fromto{(r,s)}{r}$, which is a retract of the map considered in Ex. \ref{exm:CwhenPhiisone}, which thus must induce an equivalence. Hence by functoriality in $\Phi$, it follows that $\eta$ is objectwise an equivalence on underlying quasicategories, and in particular is objectwise essentially surjective.

Suppose $I\in\Phi$, and consider the $1$-simplex of $\Lambda(\Phi)$ given by the unique active morphism $\alpha\colon\fromto{I}{\{1\}}$. The map $f_{\alpha}\colon\fromto{NA[I\to\{1\}]}{\Delta^1}$ is a retract of the inclusion $\into{\Delta^1}{NA[I\to\{1\}]}$ considered in the proof of Pr. \ref{prp:Cisrelfun}. Consequently, it induces equivalences
\begin{equation*}
\equivto{\Map^{\alpha}_{X^{\otimes}}(x,y)}{\Map^{\alpha}_{P^{\otimes}(C(X^{\otimes}))}(x,y)},
\end{equation*}
whence $\eta$ is objectwise fully faithful.
\end{proof}
\end{thm}

\begin{cnj} We expect that the axiom system given in \cite{BSP} should have an analogue for homotopy theories of weak operads over a fixed operator category $\Phi$, and that a corresponding unicity theorem should hold.
\end{cnj}

We conclude with a remark on the compatibility between the equivalences
\begin{equation*}
\Operad_{\mathrm{CSS}}^{\Phi}\simeq\Operad_{\infty}^{\Phi}
\end{equation*}
and the Boardman--Vogt tensor product.
\begin{prp}\label{prp:compareBV} Suppose $\Phi$ a perfect operator category, and suppose $X$ and $Z$ two $\Phi$-quasioperads. Then one has a canonical equivalence
\begin{equation*}
C(\Alg_{\infty}^{\Phi}(X,Z))\simeq\Alg^{\Phi}(CX,CZ).
\end{equation*}
\begin{proof} The claim is equivalent to the claim that for any simplicial space $Y$, one has a projection formula
\begin{equation*}
B(Y){{}^{\{1\}}\otimes^{\Phi}}X\simeq B(Y{{}^{\{1\}}\otimes^{\Phi}}CX).
\end{equation*}
Note that $B$ preserves all weak equivalences (as it is left Quillen for the \emph{injective} model structure on $\Fun(\Delta_{\Phi}^{\op},s\Set)$) and all colimits; hence it preserves all homotopy colimits. Consequently, we may assume that $Y$ is either $\Delta^0$ or $\Delta^1$. In the former case, the result is obvious; in the latter, it follows from a computation.
\end{proof}
\end{prp}

In general, we expect that the equivalence $C$ is fully compatible with the Boardman--Vogt tensor products for any perfect operator categories, so that
\begin{equation*}
C(Y{{}^{\Psi}\otimes^{\Phi}}X)\simeq C(Y){{}^{\Psi}\otimes^{\Phi}}C(X).
\end{equation*}
This is true when both $X$ and $Y$ are terminal quasioperads by Th. \ref{prp:UotimesUisU} and Th. \ref{prp:UotimesUisUredux}.


\section{Some examples of complete Segal $\Phi$-monoids}\label{sect:examples} We may now prove the assertions stated in the introduction. First, we note that \cite[Cor 5.1.1.5]{HA} states the following.
\begin{prp} The $\infty$-operad $E_{\infty}$ over $\FF$ is equivalent to the terminal $\infty$-operad over $\FF$.
\end{prp}
\noindent Applying the functor $C$, we may also state this result in the context of complete Segal $\FF$-operads.
\begin{cor} The complete Segal $\FF$-operad $C(E_{\infty})$ is equivalent to the terminal complete Segal $\FF$-operad $U_{\FF}$.
\end{cor}
\noindent We may equally well state this result from the perspective of the $\infty$-categories of algebras.
\begin{cor} The $\infty$-category of $E_{\infty}$-algebras in an $\infty$-operad $Z$ over $\FF$ is equivalent to the quasicategory associated with the complete Segal space $\Mon^{\FF}(Z)$.
\end{cor}

Similarly, \cite[Pr. 4.1.2.10 and Ex. 5.1.0.7]{HA}, in light of the theory developed here, state the following.
\begin{prp} The associative $\infty$-operad over $\FF$ and the $\infty$-operad $E_1$ over $\FF$ are each equivalent to the symmetrization of the terminal $\infty$-operad over $\OO$.
\end{prp}
\begin{cor} The complete Segal $\FF$-operad $C(E_1)$ is equivalent to the symmetrization of the terminal complete Segal $\OO$-operad $U_{\OO}$.
\end{cor}
\begin{cor} The $\infty$-category of $E_1$-algebras in an $\infty$-operad $Z$ over $\FF$ is equivalent to the quasicategory associated with the complete Segal space $\Mon^{\OO}(u^{\star}Z)$.
\end{cor}
 
Since the operad $A_n$ is the suboperad of $E_1$ generated by the operations of arity $\leq n$, we deduce the following.
\begin{prp} The complete Segal operad $C(A_n)$ is equivalent to the symmetrization of the terminal complete Segal $\OO_{\leq n}$-operad.
\end{prp}
\begin{cor} The $\infty$-category of $A_n$-algebras in an $\infty$-operad $Z$ over $\FF$ is equivalent to the quasicategory associated with the complete Segal space $\Mon^{\OO_{\leq n}}(u^{\star}Z)$.
\end{cor}

Finally, we have the following, which follows from \cite[Th. 5.1.2.2]{HA} and \ref{nul:BVotimesFF}.
\begin{prp}\label{prp:mainEn} For any integer $k\geq 0$, the $\infty$-operad $E_k$ over $\FF$ is equivalent to the symmetrization of the terminal $\infty$-operad over $\OO^{(k)}$.
\end{prp}

\begin{cor} For any integer $k\geq 0$, the $\infty$-category of $E_k$-algebras in an $\infty$-operad $Z$ over $\FF$ is equivalent to the quasicategory associated with the complete Segal space $\Mon^{\OO^{(k)}}(U^{\star}Z)$.
\end{cor}

\begin{cor} The $\infty$-category of $E_k$-algebras (in spaces) is equivalent to the $\infty$-category of left fibrations
\begin{equation*}
\fromto{X}{N\Lambda(\OO^{(k)})\simeq N\Theta_k^{\op}}
\end{equation*}
satisfying the Segal condition, so that for any object $I\in\OO^{(k)}$, the inert morphisms $\{\rho_i\colon\fromto{I}{\{i\}}\ |\ i\in|I|\}$ induce an equivalence
\begin{equation*}
\equivto{X_I}{\prod_{i\in|I|}X_{\{i\}}}.
\end{equation*}
\end{cor}

\begin{exm}\label{exm:bifiltrobinson} The formalism we have introduced invites the study of a wide range of new examples as well. In particular, let us contemplate the operator categories $\OO^{(n)}_{\leq m}\coloneq(\OO^{(n)})_{\leq m}$. They fit together in a diagram
\begin{equation*}
\begin{tikzpicture}
\matrix(m)[matrix of math nodes, 
row sep=4ex, column sep=4ex, 
text height=1.5ex, text depth=0.25ex] 
{\{1\}&\OO_{\leq 1}&(\OO\wr\OO)_{\leq 1}&\cdots&\OO^{(n)}_{\leq 1}&\cdots&\FF_{\leq 1}\\
\{1\}&\OO_{\leq 2}&(\OO\wr\OO)_{\leq 2}&\cdots&\OO^{(n)}_{\leq 2}&\cdots&\FF_{\leq 2}\\
\vdots&\vdots&\vdots&&\vdots&&\vdots\\
\{1\}&\OO_{\leq m}&(\OO\wr\OO)_{\leq m}&\cdots&\OO^{(n)}_{\leq m}&\cdots&\FF_{\leq m}\\
\vdots&\vdots&\vdots&&\vdots&&\vdots\\
\{1\}&\OO&\OO\wr\OO&\cdots&\OO^{(n)}&\cdots&\FF.\\}; 
\path[>=stealth,->,font=\scriptsize] 
(m-1-1) edge (m-1-2)
edge[-,double distance=1.5pt] (m-2-1)
(m-1-2) edge (m-1-3)
edge (m-2-2)
(m-1-3) edge (m-1-4)
edge (m-2-3)
(m-1-4) edge (m-1-5)
(m-1-5) edge (m-1-6)
edge (m-2-5)
(m-1-6) edge (m-1-7)
(m-1-7) edge (m-2-7)
(m-2-1) edge (m-2-2)
edge[-,double distance=1.5pt] (m-3-1)
(m-2-2) edge (m-2-3)
edge (m-3-2)
(m-2-3) edge (m-2-4)
edge (m-3-3)
(m-2-4) edge (m-2-5)
(m-2-5) edge (m-2-6)
edge (m-3-5)
(m-2-6) edge (m-2-7)
(m-2-7) edge (m-3-7)
(m-3-1) edge[-,double distance=1.5pt] (m-4-1)
(m-3-2) edge (m-4-2)
(m-3-3) edge (m-4-3)
(m-3-5) edge (m-4-5)
(m-3-7) edge (m-4-7)
(m-4-1) edge (m-4-2)
edge[-,double distance=1.5pt] (m-5-1)
(m-4-2) edge (m-4-3)
edge (m-5-2)
(m-4-3) edge (m-4-4)
edge (m-5-3)
(m-4-4) edge (m-4-5)
(m-4-5) edge (m-4-6)
edge (m-5-5)
(m-4-6) edge (m-4-7)
(m-4-7) edge (m-5-7)
(m-5-1) edge[-,double distance=1.5pt] (m-6-1)
(m-5-2) edge (m-6-2)
(m-5-3) edge (m-6-3)
(m-5-5) edge (m-6-5)
(m-5-7) edge (m-6-7)
(m-6-1) edge (m-6-2)
(m-6-2) edge (m-6-3)
(m-6-3) edge (m-6-4)
(m-6-4) edge (m-6-5)
(m-6-5) edge (m-6-6)
(m-6-6) edge (m-6-7); 
\end{tikzpicture}
\end{equation*}
By forming the symmetrization $A^n_m$ of the terminal complete Segal $\OO^{(n)}_{\leq m}$-operad, we obtain a diagram of weak symmetric operads
\begin{equation*}
\begin{tikzpicture}
\matrix(m)[matrix of math nodes, 
row sep=4ex, column sep=4ex, 
text height=1.5ex, text depth=0.25ex] 
{E_0&A_1&A^2_1&\cdots&A^n_1&\cdots&A^{\infty}_1\\
E_0&A_2&A^2_2&\cdots&A^n_2&\cdots&A^{\infty}_2\\
\vdots&\vdots&\vdots&&\vdots&&\vdots\\
E_0&A_m&A^2_m&\cdots&A^n_m&\cdots&A^{\infty}_m\\
\vdots&\vdots&\vdots&&\vdots&&\vdots\\
E_0&E_1&E_2&\cdots&E_n&\cdots&E_{\infty}.\\}; 
\path[>=stealth,->,font=\scriptsize] 
(m-1-1) edge (m-1-2)
edge[-,double distance=1.5pt] (m-2-1)
(m-1-2) edge (m-1-3)
edge (m-2-2)
(m-1-3) edge (m-1-4)
edge (m-2-3)
(m-1-4) edge (m-1-5)
(m-1-5) edge (m-1-6)
edge (m-2-5)
(m-1-6) edge (m-1-7)
(m-1-7) edge (m-2-7)
(m-2-1) edge (m-2-2)
edge[-,double distance=1.5pt] (m-3-1)
(m-2-2) edge (m-2-3)
edge (m-3-2)
(m-2-3) edge (m-2-4)
edge (m-3-3)
(m-2-4) edge (m-2-5)
(m-2-5) edge (m-2-6)
edge (m-3-5)
(m-2-6) edge (m-2-7)
(m-2-7) edge (m-3-7)
(m-3-1) edge[-,double distance=1.5pt] (m-4-1)
(m-3-2) edge (m-4-2)
(m-3-3) edge (m-4-3)
(m-3-5) edge (m-4-5)
(m-3-7) edge (m-4-7)
(m-4-1) edge (m-4-2)
edge[-,double distance=1.5pt] (m-5-1)
(m-4-2) edge (m-4-3)
edge (m-5-2)
(m-4-3) edge (m-4-4)
edge (m-5-3)
(m-4-4) edge (m-4-5)
(m-4-5) edge (m-4-6)
edge (m-5-5)
(m-4-6) edge (m-4-7)
(m-4-7) edge (m-5-7)
(m-5-1) edge[-,double distance=1.5pt] (m-6-1)
(m-5-2) edge (m-6-2)
(m-5-3) edge (m-6-3)
(m-5-5) edge (m-6-5)
(m-5-7) edge (m-6-7)
(m-6-1) edge (m-6-2)
(m-6-2) edge (m-6-3)
(m-6-3) edge (m-6-4)
(m-6-4) edge (m-6-5)
(m-6-5) edge (m-6-6)
(m-6-6) edge (m-6-7); 
\end{tikzpicture}
\end{equation*}

This is an interesting bifiltration $\{A^n_m\}_{m,n}$ of the $E_{\infty}$ operad that incorporates both the $E_n$ operads as well as the $A_n$ operads. This filtration appears to include much that is already known about obstruction theories for finding $E_{\infty}$ structures on spectra. In particular, when $n=1$, we are simply looking at the filtration of $E_1\simeq A_{\infty}$ by the operads $A_n$. When $n=\infty$, we are filtering $E_{\infty}$ by the suboperads generated by the operations of arity $\leq m$, When $m=\infty$, we are looking at the filtration of $E_{\infty}$ by the operads $E_n$. And finally, when when $m=n+1$, we expect that an algebra over $A^{n}_{n+1}$ is the same thing as an $n$-stage $E_{\infty}$ structure in the sense of Robinson \cite[\S 5.2]{MR1974890}, thereby giving his ``diagonal'' filtration on $E_{\infty}$, though we have not checked this.
\end{exm}


\appendix

\section{A proof of Th. \protect{\ref{prp:UotimesUisU}}}\label{sect:proofthatUotimesUisU} Suppose $\Phi$ and $\Psi$ two operator categories. The claim of the theorem is that the functor
\begin{equation*}
W\colon\fromto{N\Delta_{\Psi}^{\op}\times_{N\Delta^{\op}}N\Delta_{\Phi}^{\op}}{N\Delta_{\Psi\wr\Phi}^{\op}}
\end{equation*}
is a weak equivalence in the operadic model structure on $s\Set_{/N\Delta_{\Psi\wr\Phi}^{\op}}$. Note that this functor is faithful (in fact, pseudomonic), but not full. Let us call any object or morphism in its image \textbf{\emph{rectangular}}.

The proof proceeds as follows: first we find a full subcategory $\Delta_{\Psi}^{\op}\times_{\Delta^{\op}}\Delta_{\Phi}^{\op}\subset V^{\op}\subset\Delta^{\op}_{\Psi\wr\Phi}$, and we show that the inclusion $\into{N\Delta_{\Psi}^{\op}\times_{N\Delta^{\op}}N\Delta_{\Phi}^{\op}}{NV^{\op}}$ is an operadic weak equivalence. Then we replace the inclusion $\into{NV^{\op}}{N\Delta^{\op}_{\Psi\wr\Phi}}$ by a suitable left fibration $\fromto{Y}{N\Delta^{\op}_{\Psi\wr\Phi}}$, which we then show has contractible fibers by showing that every map in from a finite simplicial set can, up to homotopy, be ``coned off.''

\begin{nul} Now let $V$ denote the full subcategory of $\Delta^{\op}_{\Psi\wr\Phi}$ spanned by those objects $[K_0\to\cdots\to K_m]$ such that for every integer $1\leq i\leq m$, the fibers of the map $\fromto{K_{i-1}}{K_{i}}$ are all rectangular. Note that if $Z$ is a $\Psi\wr\Phi$-operad, then for any object $[K_0\to\cdots\to K_m]\in V$, we have a natural equivalence from $Z[K_0\to\cdots\to K_m]$ to the (homotopy) limit of a diagram
\begin{equation*}
\begin{tikzpicture} 
\matrix(m)[matrix of math nodes, 
row sep=5ex, column sep=0ex, 
text height=1.5ex, text depth=0.25ex] 
{\prod_{k_1\in|K_1|}Z[K_{0,k_1}\to\{k_1\}]&[-12ex]&\cdots&[-3ex]&[-16ex]\prod_{k_{m}\in|K_m|}Z[K_{m-1,k_m}\to\{k_m\}],\\ 
&[-2ex]\prod_{k_1\in|K_i|}Z[\{k_1\}]&&\prod_{k_{m-1}\in|K_{m-1}|}Z[\{k_{m-1}\}]&\\}; 
\path[>=stealth,->,font=\scriptsize] 
(m-1-1) edge (m-2-2) 
(m-1-3) edge (m-2-2)
edge (m-2-4)
(m-1-5) edge (m-2-4); 
\end{tikzpicture}
\end{equation*}
in which all the terms are the fiber of $Z$ over rectangular objects, and all the maps that appear are induced by rectangular morphisms. Consequently, one may extend any map of simplicial sets $\fromto{N\Delta^{\op}_{\Psi}\times_{N\Delta^{\op}}N\Delta^{\op}_{\Phi}}{Z}$ over $N\Delta^{\op}_{\Psi\wr\Phi}$, in a unique way up to homotopy, to a map $\fromto{NV^{\op}}{Z}$ over $N\Delta^{\op}_{\Psi\wr\Phi}$; in other words, the map
\begin{equation*}
\fromto{[NV^{\op},Z]_{N\Delta^{\op}_{\Psi\wr\Phi}}}{[N\Delta^{\op}_{\Psi}\times_{N\Delta^{\op}}N\Delta^{\op}_{\Phi},Z]_{N\Delta^{\op}_{\Psi\wr\Phi}}}
\end{equation*}
is a bijection. Hence the functor $W$ induces a weak equivalence
\begin{equation*}
\equivto{N\Delta^{\op}_{\Psi}\times_{N\Delta^{\op}}N\Delta^{\op}_{\Phi}}{NV^{\op}}
\end{equation*}
in the operadic model structure on $s\Set_{/N\Delta^{\op}_{\Psi\wr\Phi}}$. 
\end{nul}

\begin{nul} Now the inclusion $\into{NV^{\op}}{(N\Delta_{\Psi\wr\Phi}^{\op})_{NV^{\op}/}}$ is a trivial cofibration for the covariant model structure, and the target is fibrant \cite[Cor. 2.1.2.2]{HA}. Hence it suffices now to prove that the natural map of left fibrations $\fromto{(N\Delta_{\Psi\wr\Phi}^{\op})_{NV^{\op}/}}{N\Delta_{\Psi\wr\Phi}^{\op}}$ is a pointwise weak equivalence. In other words, we aim to show that for any object $(\mathbf{m},K)=[K_0\to\cdots\to K_m]$, the nerve of the ordinary category $V_{(\mathbf{m},K)/}$ is contractible.
\end{nul}

\begin{nul}
The objects of $V_{(\mathbf{m},K)/}$ can be described as triples
\begin{equation*}
(\eta\colon\fromto{\mathbf{m}}{\mathbf{n}},L\colon\fromto{\mathbf{n}}{\Psi\wr\Phi},\phi\colon\fromto{K}{L\circ\eta})
\end{equation*}
consisting of a morphism $\eta$ of $\Delta$, a functor $L\colon\fromto{\mathbf{n}}{\Psi\wr\Phi}$ such that every fiber of each morphism $\fromto{L_{j-1}}{L_j}$ is rectangular, and a map $\phi\colon\fromto{K}{L\circ\eta}$ such that each morphism $\phi_i\colon\into{K_i}{L_{\eta(i)}}$ is an interval inclusion and the squares
\begin{equation*}
\begin{tikzpicture} 
\matrix(m)[matrix of math nodes, 
row sep=4ex, column sep=4ex, 
text height=1.5ex, text depth=0.25ex] 
{K_{i-1}&L_{\eta(i-1)}\\ 
K_i&L_{\eta(i)}\\}; 
\path[>=stealth,->,font=\scriptsize] 
(m-1-1) edge[right hook->] (m-1-2) 
edge (m-2-1) 
(m-1-2) edge (m-2-2) 
(m-2-1) edge[right hook->] (m-2-2); 
\end{tikzpicture}
\end{equation*}
are all pullbacks. Now denote by $R_{(\mathbf{m},K)}\subset V_{(\mathbf{m},K)/}$ the full subcategory spanned by those objects $(\eta\colon\fromto{\mathbf{m}}{\mathbf{n}},L,\phi)$, such that:
\begin{itemize}
\item[---] $\eta$ has the property that $\eta(0)=0$ and $\eta(m)=n$,
\item[---] $L_{\eta(i)}=K_i$ for any $0\leq i\leq m$, and
\item[---] $\phi_i$ is the identity for any $0\leq i\leq m$.
\end{itemize}
The inclusion $i$ of this subcategory admits a retraction $r$, given by the functor that carries any triple $(\theta\colon\fromto{\mathbf{m}}{\mathbf{p}},M,\psi)\in V_{(\mathbf{m},K)/}$ to the triple $(\eta\colon\fromto{\mathbf{m}}{\mathbf{n}},L,\phi)\in R_{(\mathbf{m},K)}$ in which:
\begin{itemize}
\item[---] $\eta\colon\fromto{\mathbf{m}}{\mathbf{n}=\mathbf{p}_{\theta(0)/\ /\theta(m)}}$ is the morphism induced by $\theta$,
\item[---] $L_{j}=M_j\times_{M_{\theta(m)}}K_m$ for any $\theta(0)\leq j\leq \theta(m)$, and
\item[---] $\phi_i$ is the identity for any $0\leq i\leq m$.
\end{itemize}
The factorization of $\psi$ through $\phi$ provides a natural transformation $\fromto{i\circ r}{\id}$; hence we have a weak homotopy equivalence $NV_{(\mathbf{m},K)/}\simeq NR_{(\mathbf{m},K)}$. It thus remains to show that the nerve of $R_{(\mathbf{m},K)}$ is contractible. 
\end{nul}

\begin{nul} If there exists an integer $1\leq i\leq m$ such that the morphism $\fromto{K_{i-1}}{K_{i}}$ is an isomorphism, then it is easy to see that the natural functor
\begin{equation*}
\fromto{R_{[K_0\to\cdots\to K_{i-1}\to K_{i+1}\to\cdots\to K_m]}}{R_{[K_0\to\cdots\to K_m]}}
\end{equation*}
induces an equivalence on nerves. Consequently, we may assume that none of the morphisms $\fromto{K_{i-1}}{K_{i}}$ is an isomorphism, and consequently, that for any object $(\eta\colon\fromto{\mathbf{m}}{\mathbf{n}},L,\phi)\in R_{(\mathbf{m},K)}$, the map $\eta$ is injective.

Now for any integer $m\geq 0$, we note that the natural functor
\begin{equation*}
\fromto{R_{[K_0\to\cdots\to K_m]}}{R_{[K_0\to K_1]}\times\cdots\times R_{[K_{m-1}\to K_m]}}
\end{equation*}
is an equivalence of categories. Hence we reduce to the case in which $m=1$.
\end{nul}

\begin{nul} So suppose that $m=1$, and let us set about showing that the nerve of the category $R_{[K_0\to K_1]}$ is contractible. By definition any morphism $g\colon\fromto{L}{M}$ of $\Psi\wr\Phi$ can be factored, in an essentially unique manner, as a map $\fromto{L}{M(g)}$ covering the identity morphism in $\Phi$, followed by a morphism $\fromto{M(g)}{M}$ that is cocartesian for the coronal fibration $\fromto{\Psi\wr\Phi}{\Phi}$. Let us call this a \textbf{\emph{coronal factorization}}. The resulting sequence $[L\to M(g)\to M]$ clearly lies in $V$.

Fix a coronal factorization $[K_0\to K_{01}\to K_{1}]\in R_{[K_0\to K_1]}$. Our aim is now to produce, for any finite simplicial set $X$ and any map $g\colon\fromto{X}{NR_{[K_0\to K_1]}}$, a weakly contractible simplicial set $P_g^{\ast}\supset X$ with a distinguished vertex $v$ and an extension of $g$ to a map $G\colon\fromto{P_g^{\ast}}{NR_{[K_0\to K_1]}}$ such that $G|X=g$ and $G(v)=[K_0\to K_{01}\to K_{1}]\in R_{[K_0\to K_1]}$. This will complete the proof.
\end{nul}

\begin{nul} We now set about constructing the simplicial set $P_g^{\ast}$. Denote by $\II$ the category whose objects are objects $\mathbf{n}\in\Delta$ such that $n\geq 1$ and whose morphisms $\fromto{\mathbf{n}}{\mathbf{n}'}$ are morphisms $\fromto{\mathbf{n}}{\mathbf{n}'}$ of $\Delta$ that carry $0$ to $0$ and $n$ to $n'$. We have an obvious projection $\fromto{R_{[K_0\to K_1]}}{\II}$ given by $\goesto{(\eta\colon\fromto{\mathbf{1}}{\mathbf{n}},L,\phi)}{\mathbf{n}}$.

Let $\PP$ be the following category. The objects will be triples $(\mathbf{n},\mathbf{p},h)$, where $\mathbf{m},\mathbf{p}\in\II$, and $h$ is a functor $\fromto{\mathbf{p}}{\Fun(\mathbf{1},\mathbf{n})}$ such that $h(0)=00$, $h(m)=nn$, and for any $1\leq i\leq n$, the map $h(i-1)\to h(i)$ is of the form of one of the following:
\begin{equation*}
ij\to ik,\qquad ik\to jk,\textrm{\quad or\quad}ik\to(i+1)(k+1).
\end{equation*}
A morphism $\fromto{(\mathbf{n},\mathbf{p},h)}{(\mathbf{n}',\mathbf{p}',h')}$ of $\PP$ is a morphism $\sigma\colon\fromto{\mathbf{n}}{\mathbf{n}'}$ of $\II$ such that one has an inclusion
\begin{equation*}
\{h(i)\ |\ 0\leq i\leq n\}\subset\{\sigma(h'(i'))\ |\ 0\leq i'\leq n'\}.
\end{equation*}
We have the obvious projection $q\colon\fromto{\PP}{\II}$ given by $\goesto{(\mathbf{n},\mathbf{p},h)}{\mathbf{n}}$.

For any integer $n\geq 1$, one has a map $\fromto{\mathbf{n}}{\Fun(\mathbf{1},\mathbf{n})}$ that carries each $i$ to the identity map at $i$. This defines a section $i_0\colon\fromto{\II}{\PP}$ of $q$. There is another section $i_1\colon\fromto{\II}{\PP}$ of $q$, which carries an object $\mathbf{n}$ to the triple $(\mathbf{n},\mathbf{2},[00\to 0n\to nn])$.

For any integer $N\geq 1$, we can consider the full subcategory $\II_{\leq N}\subset\II$ spanned by those objects $\mathbf{n}$ such that $n\leq N$, and we can consider the pullback $\PP_{\leq N}\coloneq\II_{\leq N}\times_{\II}\PP$. One shows easily that when restricted to $\II_N$, there exist a zigzag of natural transformations connecting $i_0$ and $i_1$ and a zigzag of natural transformations connecting $i_0\circ q$ and the identity map on $\PP_{\leq N}$.

Now for any finite simplicial set $X$ and any map $g\colon\fromto{X}{NR_{[K_0\to K_1]}}$, let $P_g$ denote the fiber product
\begin{equation*}
\begin{tikzpicture} 
\matrix(m)[matrix of math nodes, 
row sep=4ex, column sep=4ex, 
text height=1.5ex, text depth=0.25ex] 
{P_g&&N\PP\\ 
X&NR_{[K_0\to K_1]}&N\II.\\}; 
\path[>=stealth,->,font=\scriptsize] 
(m-1-1) edge (m-1-3) 
edge (m-2-1) 
(m-1-3) edge (m-2-3) 
(m-2-1) edge (m-2-2)
(m-2-2) edge (m-2-3); 
\end{tikzpicture}
\end{equation*}
The sections $i_0$ and $i_1$ pull back to sections of $\fromto{P_g}{X}$. Since $X$ is finite, there exists an integer $N\geq 1$ such that the composite $X\to NR_{[K_0\to K_1]}\to N\II$ factors through $N\II_{\leq N}$. Hence the zigzag of homotopies between $i_0$ and $i_1$ and between $i_0\circ q$ and the identity map on $\PP_{\leq N}$ lift, and we deduce that $\fromto{P_g}{X}$ is an equivalence.

Now set $P_g^{\ast}=P_g/i_1(X)$; let $v$ be the vertex corresponding to $i_1(X)$, and regard $X$ as a simplicial subset of $P_g^{\ast}$ via $i_0$. It is thus clear that $P_g^{\ast}$ is weakly contractible.
\end{nul}

\begin{nul} Now we define the extension of $g$ to a map $G\colon\fromto{P_g^{\ast}}{NR_{[K_0\to K_1]}}$. Applying the coronal factorization repeatedly, one obtains, functorially, for any object
\begin{equation*}
[K_0=L_0\to L_1\to\cdots\to L_n=K_1]\in R_{[K_0\to K_1]},
\end{equation*}
a functor $\overline{L}\colon\fromto{\Fun(\mathbf{1},\mathbf{n})}{\Psi\wr\Phi}$ such that:
\begin{itemize}
\item[---] for any integer $0\leq i\leq n$, one has $\overline{L}_{ii}=L_i$,
\item[---] the factorization $\overline{L}_{ii}\to\overline{L}_{ij}\to\overline{L}_{jj}$ is a coronal factorization, and
\item[---] the factorization $K_0=L_0\to L_{0n}\to L_n=K_1$ is equal to our chosen coronal factorization $K_0\to K_{01}\to K_{1}$.
\end{itemize}
It is easy to check that the morphisms
\begin{equation*}
\fromto{\overline{L}_{ij}}{\overline{L}_{ik}},\qquad\fromto{\overline{L}_{ik}}{\overline{L}_{jk}},\textrm{\quad and\quad}\fromto{\overline{L}_{ik}}{\overline{L}_{(i+1)(k+1)}}
\end{equation*}
all have rectangular fibers. This defines an extension of $g$ to a map
\begin{equation*}
\fromto{P_g}{NR_{[K_0\to K_1]}},
\end{equation*}
which then factors through $P_g^{\ast}$, as desired.
\end{nul}


\section{A proof of Th. \protect{\ref{th:Tisamonad}}}\label{sect:proofthatTismonad} Suppose $\Phi$ a perfect operator category. There are two unit axioms and an associativity axiom to be checked for the endofunctor $T$ equipped with the natural transformations $\iota$ and $\mu$.

\begin{nul}\label{nul:unitaxiom} The first unit axiom is the assertion that $(T\iota)\circ\mu=\id_T$. To verify this, we first claim that for any morphism $\phi\colon\fromto{I}{T}$ of $\Phi$, one has $\sigma_I\circ\iota_I=U\lambda_I$, where $\lambda\colon\fromto{\id_{(\Phi/T)}}{E\circ\fib}$ is the unit natural transformation. Indeed, one has a commutative diagram
\begin{equation*}
\begin{tikzpicture} 
\matrix(m)[matrix of math nodes, 
row sep={8ex,between origins}, column sep={8ex,between origins}, 
text height=1.5ex, text depth=0.25ex] 
{&TJ&\\
J&TT&T(J_t)\\
&T&\\}; 
\path[>=stealth,->,font=\scriptsize]  
(m-2-1) edge node[above left]{$\iota_J$} (m-1-2)
(m-1-2) edge node[above right]{$\sigma_J$} (m-2-3)
edge (m-2-2)
(m-2-2) edge node[right]{$\chi_t$} (m-3-2)
(m-2-1) edge (m-3-2)
(m-2-3) edge (m-3-2); 
\end{tikzpicture}
\end{equation*}
and the special fiber of the composite is the identity on $J_t$. Hence one has
\begin{equation*}
(T\iota)\circ\mu=(UE\kappa)\circ(U\lambda E).
\end{equation*}
The triangle identity for the adjunction $(\fib,E)$ now implies that
\begin{equation*}
(T\iota)\circ\mu=U\id_E=\id_T,
\end{equation*}
as desired. The second unit axiom is analogous.
\end{nul}

\begin{nul} It remains to prove the associativity condition; that is, that for any object $I\in\Phi$, the diagram
\begin{equation}\label{eqn:assocsq}
\begin{tikzpicture}[baseline]
\matrix(m)[matrix of math nodes, 
row sep=4ex, column sep=6ex, 
text height=1.5ex, text depth=0.25ex] 
{T^{\,3}I&T^{\,2}I\\ 
T^{\,2}I&TI\\}; 
\path[>=stealth,->,font=\scriptsize] 
(m-1-1) edge node[above]{$T\mu_I$} (m-1-2)
(m-1-1) edge node[left]{$\mu_{TI}$} (m-2-1)
(m-1-2) edge node[right]{$\mu_I$} (m-2-2) 
(m-2-1) edge node[below]{$\mu_I$} (m-2-2); 
\end{tikzpicture}
\end{equation}
commutes.
\end{nul}
\noindent We begin with the following key technical lemma.
\begin{lem}\label{lem:iotainvofmuisequal} For any object $I$ of $\Phi$, the rectangle
\begin{equation*}
\begin{tikzpicture}[baseline]
\matrix(m)[matrix of math nodes, 
row sep=4ex, column sep=4ex, 
text height=1.5ex, text depth=0.25ex] 
{I&TI&T^{\,2}I\\
I&&TI\\}; 
\path[>=stealth,->,font=\scriptsize] 
(m-1-1) edge node[above]{$\iota_I$} (m-1-2)
(m-1-2) edge node[above]{$\iota_{TI}$} (m-1-3)
(m-1-1) edge[-,double distance=1.5pt] (m-2-1)
(m-1-3) edge node[right]{$\mu_I$} (m-2-3) 
(m-2-1) edge node[below]{$\iota_I$} (m-2-3); 
\end{tikzpicture}
\end{equation*}
is a pullback.
\begin{proof} Form the larger diagram
\begin{equation*}
\begin{tikzpicture}[baseline]
\matrix(m)[matrix of math nodes, 
row sep=4ex, column sep=4ex, 
text height=1.5ex, text depth=0.25ex] 
{I&TI&T^{\,2}I\\
I&&TI\\
\ast&&T.\\}; 
\path[>=stealth,->,font=\scriptsize] 
(m-1-1) edge node[above]{$\iota_I$} (m-1-2)
(m-1-2) edge node[above]{$\iota_{TI}$} (m-1-3)
(m-1-1) edge[-,double distance=1.5pt] (m-2-1)
(m-1-3) edge node[right]{$\mu_I$} (m-2-3) 
(m-2-1) edge node[below]{$\iota_I$} (m-2-3)
edge (m-3-1)
(m-2-3) edge node[right]{$e_I$} (m-3-3)
(m-3-1) edge node[below]{$t$} (m-3-3);
\end{tikzpicture}
\end{equation*}
The lower diagram is a pullback diagram; hence it suffices to show that the exterior rectangle is a pullback. To verify this, observe that $\mu_{\ast}=\chi_t$; hence the naturality of $\mu$ implies that the square
\begin{equation*}
\begin{tikzpicture}[baseline]
\matrix(m)[matrix of math nodes, 
row sep=4ex, column sep=6ex, 
text height=1.5ex, text depth=0.25ex] 
{T^{\,2}I&TT\\ 
TI&T\\}; 
\path[>=stealth,->,font=\scriptsize] 
(m-1-1) edge node[above]{$Te_I$} (m-1-2)
(m-1-1) edge node[left]{$\mu_I$} (m-2-1)
(m-1-2) edge node[right]{$\chi_t$} (m-2-2) 
(m-2-1) edge node[below]{$e_I$} (m-2-2); 
\end{tikzpicture}
\end{equation*}
commutes, since $e_I=T(!)$. All the rectangles of the diagram
\begin{equation*}
\begin{tikzpicture}[baseline]
\matrix(m)[matrix of math nodes, 
row sep=4ex, column sep=4ex, 
text height=1.5ex, text depth=0.25ex] 
{I&TI&T^{\,2}I\\
\ast&T&TT\\
\ast&&T.\\}; 
\path[>=stealth,->,font=\scriptsize] 
(m-1-1) edge node[above]{$\iota_I$} (m-1-2)
(m-1-2) edge node[above]{$\iota_{TI}$} (m-1-3)
(m-1-1) edge (m-2-1)
(m-1-3) edge node[right]{$Te_I$} (m-2-3) 
(m-2-1) edge node[below]{$t$} (m-2-2)
edge[-,double distance=1.5pt] (m-3-1)
(m-2-2) edge node[below]{$\iota_T$} (m-2-3)
(m-2-3) edge node[right]{$\chi_t$} (m-3-3)
(m-3-1) edge node[below]{$t$} (m-3-3);
\end{tikzpicture}
\end{equation*}
are pullback squares, whence the desired result.
\end{proof}
\end{lem}
\noindent Now we can prove the associativity.
\begin{nul}\label{thm:muworks} The first claim is that, using the structure morphism $e_I\colon\fromto{TI}{T}$, one may view the square \eqref{eqn:assocsq} as a square of $(\Phi/T)$; that is, we claim that the diagram
\begin{equation}\label{eqn:assocoverT}
\begin{tikzpicture}[baseline]
\matrix(m)[matrix of math nodes, 
row sep=4ex, column sep={6ex,between origins}, 
text height=1.5ex, text depth=0.25ex] 
{&T^2I&[4ex]TI&\\
T^3I&&&T\\ 
&T^2I&TI&\\}; 
\path[>=stealth,->,font=\scriptsize] 
(m-1-2) edge node[above]{$\mu_I$} (m-1-3)
(m-2-1) edge node[above left]{$T\mu_I$} (m-1-2)
(m-2-1) edge node[below left]{$\mu_{TI}$} (m-3-2)
(m-1-3) edge node[above right]{$e_I$} (m-2-4)
(m-3-3) edge node[below right]{$e_I$} (m-2-4) 
(m-3-2) edge node[below]{$\mu_I$} (m-3-3); 
\end{tikzpicture}
\end{equation}
commutes. To see this, consider the cube
\begin{equation*}
\begin{tikzpicture}[cross line/.style={preaction={draw=white, -, 
line width=6pt}}]
\matrix(m)[matrix of math nodes, 
row sep=2ex, column sep=3ex, 
text height=1.5ex, text depth=0.25ex]
{&T^{\,2}T&[5ex]&TT\\
T^{\,3}I&&T^{\,2}I&\\
[5ex]&TT&&T\\
T^{\,2}I&&TI&\\
};
\path[>=stealth,->,font=\scriptsize]
(m-2-1) edge node[above left]{$T^{\,2}e_I$} (m-1-2)
(m-3-2) edge node[below]{$\chi_t$} (m-3-4)
(m-1-2) edge node[above]{$\mu_T$} (m-1-4)
edge node[left]{$T\chi_t$} (m-3-2)
(m-1-4) edge node[right]{$\chi_t$} (m-3-4)
(m-2-3) edge node[below right]{$Te_I$} (m-1-4)
(m-4-3) edge node[below right]{$e_I$} (m-3-4)
(m-4-1) edge node[above left]{$Te_I$} (m-3-2)
edge node[below]{$\mu_I$} (m-4-3)
(m-2-1) edge[cross line] node[above]{$\mu_{TI}$} (m-2-3)
edge node[left]{$T\mu_I$} (m-4-1)
(m-2-3) edge[cross line] node[right]{$\mu_I$} (m-4-3);
\end{tikzpicture}
\end{equation*}
The outer square commutes since each composite $\fromto{T^{\,2}T}{T}$ is a conservative morphism $\fromto{(T^{\,2}T,\iota_{TT}\iota_T(t))}{(T,t)}$. The top square commutes by naturality, and the commutativity of the two side faces and the bottom face follows from the commutativity of \eqref{eqn:explainmu}. Thus the outer rectangle
\begin{equation*}
\begin{tikzpicture}[baseline]
\matrix(m)[matrix of math nodes, 
row sep=6ex, column sep=6ex, 
text height=1.5ex, text depth=0.25ex] 
{T^{\,3}I&T^{\,2}I\\ 
T^{\,2}T&TT\\
TT&T\\}; 
\path[>=stealth,->,font=\scriptsize] 
(m-1-1) edge node[above]{$\mu_{TI}$} (m-1-2)
(m-1-1) edge node[left]{$T^{\,2}e_I$} (m-2-1)
(m-1-2) edge node[right]{$Te_I$} (m-2-2) 
(m-2-1) edge node[above]{$\mu_T$} (m-2-2)
(m-2-1) edge node[left]{$T\chi_t$} (m-3-1)
(m-2-2) edge node[right]{$\chi_t$} (m-3-2)
(m-3-1) edge node[below]{$\chi_t$} (m-3-2); 
\end{tikzpicture}
\end{equation*}
commutes, and since the bottom square of the diagram
\begin{equation*}
\begin{tikzpicture}[baseline]
\matrix(m)[matrix of math nodes, 
row sep=6ex, column sep=6ex, 
text height=1.5ex, text depth=0.25ex] 
{T^{\,3}I&T^{\,2}I\\ 
T^{\,2}I&TI\\
TT&T\\}; 
\path[>=stealth,->,font=\scriptsize] 
(m-1-1) edge node[above]{$\mu_{TI}$} (m-1-2)
(m-1-1) edge node[left]{$T\mu_I$} (m-2-1)
(m-1-2) edge node[right]{$\mu_I$} (m-2-2) 
(m-2-1) edge node[above]{$\mu_I$} (m-2-2)
(m-2-1) edge node[left]{$Te_I$} (m-3-1)
(m-2-2) edge node[right]{$e_I$} (m-3-2)
(m-3-1) edge node[below]{$\chi_t$} (m-3-2); 
\end{tikzpicture}
\end{equation*}
commutes, it follows that \eqref{eqn:assocoverT} does as well.

Now write $K$ for the special fiber of the composite
\begin{equation*}
T^{\,3}I\ \tikz[baseline]\draw[>=stealth,->,font=\scriptsize](0,0.5ex)--node[above]{$\mu_{TI}$}(0.75,0.5ex);\ T^{\,2}I\ \tikz[baseline]\draw[>=stealth,->,font=\scriptsize](0,0.5ex)--node[above]{$\mu_I$}(0.75,0.5ex);\ TI\ \tikz[baseline]\draw[>=stealth,->,font=\scriptsize](0,0.5ex)--node[above]{$e_I$}(0.75,0.5ex);\ T,
\end{equation*}
and write $L$ for the special fiber of the composite
\begin{equation*}
T^{\,3}I\ \tikz[baseline]\draw[>=stealth,->,font=\scriptsize](0,0.5ex)--node[above]{$T\mu_I$}(0.75,0.5ex);\ T^{\,2}I\ \tikz[baseline]\draw[>=stealth,->,font=\scriptsize](0,0.5ex)--node[above]{$\mu_I$}(0.75,0.5ex);\ TI\ \tikz[baseline]\draw[>=stealth,->,font=\scriptsize](0,0.5ex)--node[above]{$e_I$}(0.75,0.5ex);\ T.
\end{equation*}
By adjunction it suffices to show that the two morphisms of special fibers $\fromto{K}{I}$ and $\fromto{L}{I}$ are equal. To compute $K$, consider the following diagram
\begin{equation*}
\begin{tikzpicture}[baseline]
\matrix(m)[matrix of math nodes, 
row sep=4ex, column sep=4ex, 
text height=1.5ex, text depth=0.25ex] 
{I&TI&T^{\,2}I&T^{\,3}I\\
I&TI&&T^{\,2}I\\
I&&&TI\\
\ast&&&T\\}; 
\path[>=stealth,->,font=\scriptsize] 
(m-1-1) edge node[above]{$\iota_I$} (m-1-2)
(m-1-2) edge node[above]{$\iota_{TI}$} (m-1-3)
(m-1-3) edge node[above]{$\iota_{T^{\,2}I}$} (m-1-4)
(m-2-1) edge node[below]{$\iota_I$} (m-2-2)
(m-2-2) edge node[below]{$\iota_{TI}$} (m-2-4)
(m-3-1) edge node[below]{$\iota_I$} (m-3-4)
(m-4-1) edge node[below]{$t$} (m-4-4)
(m-1-1) edge[-,double distance=1.5pt] (m-2-1)
(m-2-1) edge[-,double distance=1.5pt] (m-3-1)
(m-3-1) edge (m-4-1)
(m-1-2) edge[-,double distance=1.5pt] (m-2-2)
(m-1-4) edge node[right]{$\mu_{TI}$} (m-2-4)
(m-2-4) edge node[right]{$\mu_I$} (m-3-4)
(m-3-4) edge node[right]{$e_I$} (m-4-4);
\end{tikzpicture}
\end{equation*}
By \ref{lem:iotainvofmuisequal}, every rectangle of this diagram is a pullback. To compute $L$, consider the following diagram
\begin{equation*}
\begin{tikzpicture}[baseline]
\matrix(m)[matrix of math nodes, 
row sep=4ex, column sep=4ex, 
text height=1.5ex, text depth=0.25ex] 
{I&TI&T^{\,2}I&T^{\,3}I\\
I&TI&&T^{\,2}I\\
I&&&TI\\
\ast&&&T\\}; 
\path[>=stealth,->,font=\scriptsize] 
(m-1-1) edge node[above]{$\iota_I$} (m-1-2)
(m-1-2) edge node[above]{$\iota_{TI}$} (m-1-3)
(m-1-3) edge node[above]{$\iota_{T^{\,2}I}$} (m-1-4)
(m-2-1) edge node[below]{$\iota_I$} (m-2-2)
(m-2-2) edge node[below]{$\iota_{TI}$} (m-2-4)
(m-3-1) edge node[below]{$\iota_I$} (m-3-4)
(m-4-1) edge node[below]{$t$} (m-4-4)
(m-1-1) edge[-,double distance=1.5pt] (m-2-1)
(m-2-1) edge[-,double distance=1.5pt] (m-3-1)
(m-3-1) edge (m-4-1)
(m-1-2) edge[-,double distance=1.5pt] (m-2-2)
(m-1-4) edge node[right]{$T\mu_I$} (m-2-4)
(m-2-4) edge node[right]{$\mu_I$} (m-3-4)
(m-3-4) edge node[right]{$e_I$} (m-4-4);
\end{tikzpicture}
\end{equation*}
Again by \ref{lem:iotainvofmuisequal}, and since $T$ preserves all pullbacks, every rectangle of this diagram is a pullback. Hence $K=L=I$, and the morphisms $\fromto{K}{I}$ and $\fromto{L}{I}$ are each simply the identity.
\end{nul}

This completes the proof of the theorem.


\section{A proof of Th. \protect{\ref{thm:admissiblesinducecolax}}}\label{sect:proofthatadmissiscolax} Fix an admissible functor $F\colon\fromto{\Psi}{\Phi}$ between perfect operator categories. We have to show that the two diagrams \eqref{eq:colaxmorph} commute with $C=\Psi$, $D=\Phi$, and $\eta=\alpha_F$. 

\begin{nul} Note that $F$ induces a functor
\begin{equation*}
F_{/TT}\colon\fromto{(\Psi/T_{\Psi}T_{\Psi})}{(\Phi/T_{\Phi}T_{\Phi})}
\end{equation*}
that assigns to any object $\fromto{I}{T_{\Psi}T_{\Psi}}$ the composite
\begin{equation*}
FI\ \tikz[baseline]\draw[>=stealth,->,font=\scriptsize](0,0.5ex)--(0.5,0.5ex);\ FT_{\Psi}T_{\Psi}\ \tikz[baseline]\draw[>=stealth,->,font=\scriptsize](0,0.5ex)--node[above]{$\alpha_{F,T_{\Psi}}$}(1,0.5ex);\ T_{\Phi}FT_{\Psi}\ \tikz[baseline]\draw[>=stealth,->,font=\scriptsize](0,0.5ex)--node[above]{$T_{\Phi}\chi_{F(t_{\Psi})}$}(1.25,0.5ex);\ T_{\Phi}T_{\Phi}.
\end{equation*}
\end{nul}

The commutativity of the first diagram of \eqref{eq:colaxmorph} follows directly from the following.
\begin{lem}\label{lem:protocommutemu} The following diagram of natural transformations of functors
\begin{equation*}
\fromto{(\Phi/T_{\Phi})}{(\Psi/T_{\Psi})}
\end{equation*}
commutes:
\begin{equation*}
\begin{tikzpicture}[baseline]
\matrix(m)[matrix of math nodes, 
row sep=6ex, column sep=6ex, 
text height=1.5ex, text depth=0.25ex] 
{F_{/T}\circ\chi_{t_{\Psi},!}\circ E_{\Psi,/T_{\Psi}}&\chi_{t_{\Phi},!}\circ F_{/TT}\circ E_{\Psi,/T_{\Psi}}&\chi_{t_{\Phi},!}\circ E_{\Phi,/T_{\Phi}}\circ F_{/T}\\
F_{/T}\circ E_{\Psi}\circ \fib_{\Psi}&E_{\Phi}\circ F\circ\fib_{\Psi}&E_{\Phi}\circ \fib_{\Phi}\circ F_{/T}.\\}; 
\path[>=stealth,->,font=\scriptsize] 
(m-1-1) edge[-,double distance=1.5pt] (m-1-2)
edge node[left]{$F_{/T}\circ\sigma_{\Psi}$} (m-2-1)
(m-1-2) edge node[above]{$\chi_{t_{\Phi},!}\circ \alpha_F$} (m-1-3)
(m-1-3) edge node[right]{$\sigma_{\Phi}\circ F_{/T}$} (m-2-3)
(m-2-1) edge node[below]{$\alpha_F\circ\fib_{\Psi}$} (m-2-2)
(m-2-2) edge[-,double distance=1.5pt] (m-2-3);
\end{tikzpicture}
\end{equation*}
\begin{proof} Suppose $\psi\colon\fromto{J}{T_{\Psi}}$ a morphism. The claim is that the diagram
\begin{equation*}
\begin{tikzpicture}[baseline]
\matrix(m)[matrix of math nodes, 
row sep=4ex, column sep=4ex, 
text height=1.5ex, text depth=0.25ex] 
{FT_{\Psi}J&T_{\Phi}FJ\\
FT_{\Psi}J_{t_{\Psi}}&T_{\Phi}(FJ)_{F(t_{\Psi})}\\}; 
\path[>=stealth,->,font=\scriptsize] 
(m-1-1) edge (m-1-2)
edge (m-2-1)
(m-1-2) edge (m-2-2)
(m-2-1) edge (m-2-2);
\end{tikzpicture}
\end{equation*}
of $(\Phi/T_{\Phi})$ commutes, where the structure morphisms in question are
\begin{equation*}
\begin{tikzpicture}[baseline]
\matrix(m)[matrix of math nodes, 
row sep=2ex, column sep=3ex, 
text height=1.5ex, text depth=0.25ex] 
{FT_{\Psi}J&&FT_{\Psi}T_{\Psi}&&FT_{\Psi}&&T_{\Phi},\\
T_{\Phi}FJ&&T_{\Phi}FT_{\Psi}&&T_{\Phi}T_{\Phi}&&T_{\Phi},\\
FT_{\Psi}J_{t_{\Psi}}&&&FT_{\Psi}&&&T_{\Phi},\\}; 
\path[>=stealth,->,font=\scriptsize] 
(m-1-1) edge node[above]{$FT_{\Psi}\psi$} (m-1-3)
(m-1-3) edge node[above]{$F(\chi_{t_{\Psi}})$} (m-1-5)
(m-1-5) edge node[above]{$\chi_{F(t_{\Psi})}$} (m-1-7)
(m-2-1) edge node[above]{$T_{\Phi}F\psi$} (m-2-3)
(m-2-3) edge node[above]{$T_{\Phi}\chi_{F(t_{\Psi})}$} (m-2-5)
(m-2-5) edge node[above]{$\chi_{t_{\Phi}}$} (m-2-7)
(m-3-1) edge node[above]{$Fe_{J_{t_{\Psi}}}$} (m-3-4)
(m-3-4) edge node[above]{$\chi_{F(t_{\Psi})}$} (m-3-7);
\end{tikzpicture}
\end{equation*}
and of course the usual structure morphism $e_{F(J_{t_{\Psi}})}\colon\fromto{T_{\Phi}F(J_{t_{\Psi}})}{T_{\Phi}}$. By adjunction it suffices to show that the square of special fibers
\begin{equation*}
\begin{tikzpicture}[baseline]
\matrix(m)[matrix of math nodes, 
row sep=4ex, column sep=4ex, 
text height=1.5ex, text depth=0.25ex] 
{\left(FT_{\Psi}J\right)_{t_{\Phi}}&\left(T_{\Phi}FJ\right)_{t_{\Phi}}\\
\left(FT_{\Psi}J_{t_{\Psi}}\right)_{t_{\Phi}}&\left(T_{\Phi}(FJ)_{F(t_{\Psi})}\right)_{t_{\Phi}}\\}; 
\path[>=stealth,->,font=\scriptsize] 
(m-1-1) edge (m-1-2)
edge (m-2-1)
(m-1-2) edge (m-2-2)
(m-2-1) edge (m-2-2);
\end{tikzpicture}
\end{equation*}
commutes. But since the special fiber of the morphisms $\alpha_{F,J}$ and $\alpha_{F,J_{t_{\Psi}}}$ are each the identity, the special fiber square is in particular commutative.
\end{proof}
\end{lem}

\begin{nul} We now prove the commutativity of the second diagram of \eqref{eq:colaxmorph}. The claim is that for any object $J$ of $\Psi$, the diagram
\begin{equation*}
\begin{tikzpicture}[baseline]
\matrix(m)[matrix of math nodes, 
row sep=2ex, column sep=3ex, 
text height=1.5ex, text depth=0.25ex] 
{&FT_{\Psi}J\\
FJ&\\
&T_{\Phi}FJ.\\}; 
\path[>=stealth,->,font=\scriptsize] 
(m-2-1) edge node[above left]{$F\iota_J$} (m-1-2)
edge node[below left]{$\iota_{FJ}F$} (m-3-2)
(m-1-2) edge node[right]{$\alpha_{F,J}$} (m-3-2);
\end{tikzpicture}
\end{equation*}
commutes. By composing with the structure map $e_{FJ}\colon\fromto{T_{\Phi}FJ}{T_{\Phi}}$, this can be regarded as a diagram of $(\Phi/T_{\Phi})$. Hence it suffices to verify that the special fiber triangle
\begin{equation*}
\begin{tikzpicture}[baseline]
\matrix(m)[matrix of math nodes, 
row sep=2ex, column sep=3ex, 
text height=1.5ex, text depth=0.25ex] 
{&\left(FT_{\Psi}J\right)_{t_{\Phi}}\\
\left(FJ\right)_{t_{\Phi}}&\\
&\left(T_{\Phi}FJ\right)_{t_{\Phi}}\\}; 
\path[>=stealth,->,font=\scriptsize] 
(m-2-1) edge (m-1-2)
edge (m-3-2)
(m-1-2) edge (m-3-2);
\end{tikzpicture}
\end{equation*}
commutes. But since the special fiber of $\alpha_{F,J}\colon\fromto{FT_{\Psi}J}{T_{\Phi}FJ}$ is the identity on $FJ$, the special fiber triangle commutes.
\end{nul}

This completes the proof of the theorem.


\bibliographystyle{amsplain}
\bibliography{hca}

\end{document}